
\input amssym.def
\input amssym.tex 
\magnification 1100 

\catcode`\@=12

\hsize=34pc\vsize=47pc
\voffset=0pc  \hoffset=.3pc

\font\ninerm=cmr9 at 9pt
\font\ninebf=cmbx10 at 9pt

\tolerance=6000

\parskip=0pt
\parindent=18pt
\baselineskip=13pt
\def\myitem{\noindent\hangindent=1.5em\hangafter=1}

\def\theo#1#2{\medskip
  {\noindent\bf   Theorem #1.} {\it #2}\medskip}
\def\lem#1#2{\medskip{\noindent \bf  Lemma #1.}
  {\it #2}\medskip}
\def\cor#1#2{\medskip{\noindent \bf  Corollary #1.}
  {\it #2}\medskip}
\def\prop#1#2{\medskip{\noindent 
\bf  Proposition #1.} {\it #2}\medskip}

\def\rem#1#2{\medskip {\noindent \bf Remark #1.} {\rm #2}
\medskip}

\def\proof#1{\medskip {\noindent \bf  Proof.} {\rm #1} \medskip}
\def\prooff#1#2{\medskip {\noindent \bf  Proof #1.} 
{\rm #2} \medskip}
\def\ref#1{\vskip 1.5pc{\centerline 
{\bf References}}\vskip 4pt \noindent}
\def\affil#1
{\vskip 1truepc{\noindent {\ninerm #1}}}
\def\qed{{{\hfill $\blacksquare$}}\medskip}

\def\mypar{\vskip 0pc}

 \abovedisplayskip=12pt plus3pt minus2pt
 \belowdisplayskip=12pt plus3pt minus2pt
 \abovedisplayshortskip=12pt plus3pt minus2pt
 \belowdisplayshortskip=12pt plus3pt minus2pt

 \def\lskipamount{12pt}
 \def\lskip{\vskip\lskipamount plus3pt minus2pt}
 \def\lbreak{\par \ifdim\lastskip<\lskipamount
  \removelastskip \penalty-200 \lskip \fi}

 \def\lnobreak{\par \ifdim\lastskip<\lskipamount
  \removelastskip \penalty200 \lskip \fi}

\font\titlefont=cmr12 at 14pt
\font\sectionfont=cmbx10 at 10pt
\font\autit=cmr12 


\def\Sec#1{\vskip 1.5truepc\centerline{\hbox {{\sectionfont #1}}}
\vskip1truepc\noindent}
\def\<{{\langle}}
\def\>{{\rangle}}
\def\pdb#1{{\partial\phantom{x}\over\partial#1}}
\def\pa{{\partial}}
\def\ppd#1#2{{\partial#1\over\partial#2}}

\def\ovl{\overline}
\def\wt{\widetilde}
\def\dimR{\dim_{\Bbb R}}
\def\dimC{\dim_{\Bbb C}}
\def\d{{\mathord{\rm d}}}
\def\hol{{\mathord{\rm hol}}}
\def\ee{{\mathord{\rm e}}}

\def\({{\rm(}}
 \def\){{\rm)}}
\def\gr{\hbox{\rm gr}\,}
\def\ad{\hbox{\rm ad}}

\def\Cent{\hbox{\rm Cent}\,}

\def\End{{\hbox{\rm End}}\,}

\def\Ker{{\hbox{\rm Ker}}}
\def\Unit{{\hbox{\rm Unit}\,}}

\def\Herm{{\hbox{\rm Herm}}}

\def\i{^{-1}}
\def\half{{1\over 2}}
\def\Cinf{C^{\infty}}

\def\la{\lambda}  
\def\ga{\gamma}
\def\om{\omega}
\def\al{\alpha}
\def\be{\beta}  
\def\ep{\varepsilon}
\def\kap{\kappa}
\def\sig{\sigma}   \def\Sig{\Sigma}
\def\th{\theta}
\def\ze{\zeta}
\def\a{{\frak a}}
\def\b{{\frak b}}
\def\g{{\frak g}}
\def\h{{\frak h}}
\def\j{{\frak j}}
\def\k{{\frak k}}
\def\l{{\frak l}}
\def\m{{\frak m}}

\def\p{{\frak p}}
\def\q{{\frak q}}
\def\r{{\frak r}}
\def\s{{\frak s}}
\def\t{{\frak t}}
\def\u{{\frak u}}
\def\oY{{\ovl{Y}}}\def\oX{{\ovl{X}}}
\def\oz{{\ovl{z}}}  \def\ow{{\ovl{w}}}
\def\ophi{{\ovl{\phi}}}
\def\oe{{\ovl{e}}}
\def\og{{\ovl{g}}}

\def\so{{\frak s\frak o}}
\def\sp{{\frak s\frak p}}
\def\sl{{\frak s\frak l}}
\def\gl{{\frak g\frak l}}
\def\Ug{{{\cal U}({\frak g})}}

\def\Vect{{\frak Vect}}
\def\mapright#1{\smash
{\mathop{
\longrightarrow}\limits^{#1}}}

\def\iso{{\,\mapright{\sim}\,}}

\def\C{{\Bbb C}}
\def\Z{{\Bbb Z}}
\def\R{{\Bbb R}}

\def\OR{O_{\Bbb R}}

\def\Omin{O_{\rm min}}
\def\H{{\cal H}}

\def\A{{\cal A}}
\def\D{{\cal D}} 
\def\E{{\cal E}}

\def\I{{\cal I}}
\def\J{{\cal J}}
\def\L{{\cal L}}
\def\Q{{\cal Q}} 
\def\T{{\cal T}} 
\def\U{{\cal U}} 
\def\NH{{\bf N}^{1\over 2}}
\def\NN{{\bf N}}
\def\fo{f_{0}}
\def\ofo{\ovl{f}_0}
\def\fa{f_{1}}
\def\fm{f_{m}} 

\def\ks{\k^{\s}} \def\Ks{K^{\s}}
\def\kl{\k^{\l}} \def\Kl{K^{\l}}

\def\Yo{Y^{o}}
\def\LE{{\cal L}_{E}}

\def\sL{\sqrt{\Lambda}}

\def\vep{\varepsilon}
\def\kap{\kappa}

\def\norm{{|\hskip -1.5pt |}}

\def\biggnorm{{\bigg|\hskip -1.5pt \bigg|}}

\def\ns#1{\norm #1\norm^2}
\def\biggns#1{\biggnorm {#1}\biggnorm^2}

\def\of{\ovl{f}}
\def\mdeg{{\mathord{\bf deg}}}

\def\RP{\R^{+}}
\def\ZP{\Z_{+}}
\def\gR{\frak g_{\Bbb R}}
\def\GR{G_{\Bbb R}}
\def\KR{K_{\Bbb R}}
\def\pR{\frak p_{\Bbb R}}
\def\kR{\frak k_{\Bbb R}}

\def\V{{\cal V}}
\def\ee{\mathord{\rm e}} 
\def\del{\partial} \def\odel{\ovl{\partial}} 
\def\om{\omega}   \def\Om{\Omega}
\def\Re{\mathord{\rm Re}\,} \def\Im{\mathord{\rm Im}\,}
\def\symbol{\mathord{\rm symbol}\;}
 \def\JJ{{\bf J}}

\def\YoY{Y\times\oY}
\def\ov{\ovl{v}} \def\ou{\ovl{u}} 
\def\oH{\ovl{H}}
\def\oNH{\ovl{{\bf N}}^{1\over 2}} 

\def\K{{\cal K}}
 \def\oV{\ovl{V}}
\def\AND{\qquad\hbox{and}\qquad} 
\def\and{\quad\hbox{and}\quad} 
\def\pd{\partial}

\def\ek{e_{\k}} \def\oek{\oe_{\k}}
\def\Hs{H^{\sharp}}
\def\mm{^{\m_0^-}}
\def\soe{_{\oe}}

\def\bt{{\bf t}} \def\bv{{\bf v}}
\def\adh{\ad\, h}

\def\rbot{\,\hbox to 5pt{\leaders\hrule\hfil}
   \vbox to 5pt{\leaders\vrule\vfil}\thinspace}

 
\headline={\hfill} 
\centerline{\titlefont Geometric Quantization of Real Minimal 
Nilpotent Orbits. }     
\footnote{}
\vskip2pc
\centerline{\autit Ranee Brylinski\footnote{\rm *}
{\noindent \rm  Research 
supported in part  by a Sloan Foundation fellowship and
NSF  Grant No. DMS-9505055}}

\vskip 2pc 
{\ninerm {\ninebf Abstract:}
In this paper, we begin a quantization program for nilpotent orbits $\OR$ of a
real semisimple Lie group $\GR$. These orbits arise naturally as  the coadjoint
orbits of $\GR$ which  are stable under scaling, and thus  they have a 
canonical symplectic structure $\om$ where the  $\GR$-action is  
Hamiltonian. These orbits and their covers    generalize the oscillator phase
space $T^*\R^n$, which occurs   here when   $\GR=Sp(2n,\R)$ and $\OR$ is
minimal.

A complex  structure $\JJ$ polarizing $\OR$ and invariant under a maximal 
compact subgroup $\KR$ of $\GR$ is provided by the Kronheimer-Vergne
Kaehler structure $(\JJ,\om)$. We argue that the Kaehler potential serves as 
the Hamiltonian. Using this setup, we  realize the Lie algebra $\gR$ of $\GR$ as
a  Lie algebra of  rational functions on the holomorphic cotangent bundle
$T^*Y$ where $Y=(\OR,\JJ)$.

Thus we transform the quantization problem on $\OR$ into a quantization 
problem on $T^*Y$. We explain this in detail and solve the new quantization
problem on $T^*Y$ in a uniform manner for     minimal nilpotent orbits in the 
non-Hermitian case.  The Hilbert space of quantization consists of 
holomorphic   half-forms on $Y$. We construct the reproducing kernel.   The
Lie algebra $\gR$ acts by explicit pseudo-differential operators on half-forms
where the energy operator quantizing the Hamiltonian is inverted.  The Lie
algebra representation exponentiates to give a minimal unitary ladder 
representation  of a cover of $\GR$. Jordan algebras play a key role in the
geometry and the quantization.}

\Sec{\S1. Introduction.}  

{\bf I. Quantization of Phase Space.}
Quantization of a classical phase space $M$ with symplectic form
$\om$ is a process whereby observables $\phi$ are converted into
self-adjoint operators ${\Q}(\phi)$ on a Hilbert space $\H$  of states. 
The observables are simply the smooth functions on $M$.

The Hilbert space $\H$ should arise, according to the philosophy of
Geometric Quantization, as a space of polarized sections of a suitable
complex line bundle over $M$.   A real (complex) polarization of $M$
consists     of a integrable Lagrangian distribution inside the
(complexified) tangent bundle. A polarized section, of a bundle with
connection,  is  a section annihilated by all vector fields lying in the
polarization; in the real case, this means that  the section is covariantly
constant along the leaves of the  corresponding Lagrangian foliation.

We require that the quantization  satisfies Dirac's axioms (see e.g.,
[Ki], [A-M]) in some form. Dirac's consistency axiom is that the Poisson
bracket of functions on $M$ goes over into the commutator of
operators  so that
$$\Q(\{\phi,\psi\})=i[\Q(\phi),\Q(\psi)]\eqno(1.1)$$ 
(We have set $\hbar=1$.) 
Additional axioms mandate that the constant
function $1$ quantizes to the identity operator, and a complete set of
observables quantizes to give a complete set of operators.

In Hamiltonian mechanics, the physics of the system in encoded in a
single observable $F$ (usually written as  $E$ or $H$)  called the {\it
Hamiltonian}. Often $F$ is the total energy.  

Any observable $\phi$
generates a Hamiltonian flow: this is the  flow of the Hamiltonian 
vector field $\xi_\phi$ defined by the equation
$$\xi_\phi\rbot\om+\d\phi=0\eqno(1.2)$$  
The Poisson bracket on $\Cinf(M)$ is given by
$\{\phi,\psi\}=\xi_\phi(\psi)=\om(\xi_\phi,\xi_\psi)$.

The Hamiltonian flow of  the
$F$ gives   the time evolution of the physical system.  For any
observable $\phi$, the time derivative $\dot{\phi}$ of $\phi$ as the
system evolves is given by $\dot{\phi}=\{F,\phi\}$. This is a concise
version of Hamilton's equations. On physical grounds, 
in certain circumstances, $F$ should be a
positive function on $M$.

In quantization of a Hamiltonian mechanical system,  
$F$ should be promoted  to a self-adjoint operator $\Q(F)$ on $\H$ 
with positive spectrum. When $F$ is the total classical energy,  the
spectrum of $\Q(F)$  should be discrete and give the possible quantized
energy levels of the quantum system.

{\bf  II. Quantization of the $n$-dimensional Harmonic Oscillator.}
The most familiar model situation is  the case where  $M$ is
the cotangent bundle of some (configuration) manifold $X$ and $\om$ 
is the  canonical symplectic form so that  $\om=\d\th$ where $\th$ is 
the Liouville $1$-form on $T^*X$. In this case we have the manifest
cotangent polarization where the leaves are the cotangent spaces of
$X$ .  We expect $\H$ to be a space  of square
integrable half-forms on $X$ (see  \S2 and below starting around (1.5)).
 A smooth function 
$f$  on $X$ quantizes to a give a multiplication operator on $\H$.  If
$\eta$ is a vector field on $X$, then the  symbol $\sig^\eta$ quantizes to
the Lie derivative 
$\L_\eta$ operator on half-forms.   Consistent  quantization of
additional observables is problematic, as we see already in the
oscillator example below.

A second model situation is the case where $M$ is a Kaehler manifold
and $\om$ is the Kaehler form.  Then the complex structure $\JJ$ of
$M$ gives a complex polarization.   Now ``polarized" simply means
``holomorphic". Thus $\H$ should be a space of holomorphic
square-integrable sections of a suitable holomorphic complex line
bundle over $M$.

The most familiar  example of   a Hamiltonian mechanical system, the 
oscillator phase space, admits both cotangent and Kaehler
polarizations.   The oscillator phase space is
$M=T^*\R^n$.
The canonical coordinates on $T^*\R^n$ are the position coordinates
$q_1,\dots,q_n$  together with the momentum coordinates 
$p_1,\dots,p_n$. The canonical symplectic form is  
$\om=\sum_{k=1}^n\d p_k\wedge\d q_k$. The  Poisson bracket
satisfies $\{p_j,p_k\}=\{q_j,q_k\}=0$ and $\{p_j,q_k\}=\delta_{jk}$. For 
general observables we have the classical  formula
$$\{\phi,\psi\}=\sum_{k=1}^n\left(
{\pd{\phi}\over\pd{p_k}}{\pd{\psi}\over\pd{q_k}}-
{\pd{\psi}\over\pd{p_k}}{\pd{\phi}\over\pd{q_k}}\right)\eqno(1.3) $$

In physics, $T^*\R^n$ arises as the phase space of $n$ uncoupled 
harmonic oscillators with Hamiltonian equal to  the total energy
(kinetic plus potential)
$$F=\half\sum_{k=1}^n(p_k^2+q_k^2)\eqno(1.4)$$

We also have a natural Kaehler structure.
We  identify $T^*\R^n=\R^{2n}=\C^n$ so that the  complex-valued
observables $z_k=(p_k+iq_k)/\sqrt{2}$ are holomorphic coordinates.
Now $\C^n$ is a Kaehler manifold with  Kaehler  form $\om$ and
Kaehler metric    $g=\sum_{k=1}^n(\d p_k^2+\d q_k^2)$.
In the $z_j,\oz_k$ coordinates we have 
$\om=i\sum_{k=1}^n\d\oz_k\wedge\d z_k$  and the Poisson bracket
satisfies $\{z_j,z_k\}=\{\oz_j,\oz_k\}=0$ and
$\{\oz_j,z_k\}=i\delta_{jk}$. Also (1.4) becomes 
$$F=\sum_{k=1}^n|z_k|^2 \eqno(1.5)$$

The quantization of  the Kaehler phase space $M=\C^n$ gives the
Fock-Bargmann model of the quantum mechanical oscillator.
(Quantization by means of  the real cotangent polarization gives the 
Schroedinger model.)  In this model, $\H$ is a space of holomorphic
functions $f(z_1,\dots,z_n)$ on $\C^n$.
The Hamiltonian $F$ quantizes into the energy operator
$$\Q(F)=\sum_{k=1}^n\left(z_k\pdb{z_k}+\half
\right)\eqno(1.6)$$
The functions $z_k$ and $\oz_k$ quantize into the creation and
annihilation operators
$$\Q(z_k)=z_k\AND   \Q(\oz_k)=\pdb{z_k}     \eqno(1.7)$$
Then $\Q(F)$ is a grading operator on the quantum space and 
$\Q(z_k)=z_k$ and $\Q(\oz_k)$  are raising and lowering operators 
moving  the eigenspaces of   $\Q(F)$.

One way to ``explain" the $\half$-shift  in (1.6) (a
quantum correction)   is to  adopt the symmetrization procedure of
canonical quantization so that  
$$\Q(z_k\oz_k)=\half\left(\Q(z_k)\Q(\oz_k)+\Q(\oz_k)\Q(z_k)\right)
=\half\left(z_k\pdb{z_k}+\pdb{z_k}z_k\right)=z_k\pdb{z_k}+\half
\eqno(1.8)$$
  
There is a unique Hermitian inner product $\<f|g\>$ on the space 
$H=\C[z_1,\dots,z_n]$ of polynomial functions  such that the operators
$\Q(z_k)$ and $\Q(\oz_k)$ in (1.7) are mutually adjoint. (The condition
that
$\Q(\phi)$ is self-adjoint for real $\phi$ amounts to the condition that 
$\Q(\phi)$ and $\Q(\ophi)$ are mutually adjoint for complex
$\phi$.) This   inner product is positive definite with 
$$\ns{z_{1}^{a_1}\cdots z_{m}^{a_m}}=a_{1}!\cdots a_{m}!  \eqno(1.9)$$
 
The inner product (1.9) is  given by the integral formula
$$\<f|g\>=\int_{\C^m}f(z)\ovl{g(z)}\ee^{-|z|^2 }|\d z\d\oz|
\eqno(1.10)$$
and this expression  defines the inner product on the 
Hilbert space completion  $\H$   of $H$. Thus $\H$ consists of  all the
holomorphic functions $f(z_1,\dots,z_n)$ on $\C^n$ which are  ``square
integrable" in the sense that $\ns{f}=\<f|f\>$ is finite.

The reproducing kernel (see \S8) of $\H$ is the holomorphic function
$\K(z,\ow)$ on $\C^n\times\ovl{\C^n}$
$$\K(z,\ow)=\exp(z_1\ow_1+\cdots+z_n\ow_n)\eqno(1.11)$$ 
Here $\oX$ denotes the complex conjugate manifold to a complex 
manifold $X$, so that holomorphic functions on $\oX$ identify with
anti-holomorphic functions on $X$.  $\oX$ is obtained from $X$ by
reversing the sign of  the complex structure.

The Hamiltonian flow of $F$ lies inside a larger symmetry.
The   Hamiltonian $F$  sits inside the
space  $\g$  of all  homogeneous quadratic 
polynomials $z_1,\dots,z_n,\oz_1,\dots,\oz_n$. The space $\g$ is a 
finite-dimensional  Lie subalgebra of complex-valued observables 
under Poisson bracket. The Lie algebra $\g$ breaks naturally into $3$
pieces: $\g=\k\oplus\p^+\oplus\p^-$ where
$$\k=\hbox{span  of } \,z_j\oz_k,\qquad 
\p^+=\hbox{span  of } \,z_j z_k,\qquad 
\p^-=\hbox{span  of } \,\oz_j\oz_k  \eqno(1.12)$$
Here $\k$ arises as the subspace of $\phi\in\g$ which Poisson commute
with $F$ so that $\phi$ is a conserved quantity. Then $\p^+$ and $\p^-$
are the irreducible $\k$-representations in $\g$ complementary to 
$\k$.

The subspace $\gR\subset\g$ of real-valued observables  is a 
Lie algebra real form  of $\g$. We have $\gR=\kR\oplus\pR$ where
$$\eqalign{
\kR=&\hbox{span  of } \,z_j\oz_k+\oz_jz_k
\and i(z_j\oz_k-\oz_jz_k)\cr
\pR=&\hbox{span  of } \,z_j z_k+ \oz_j\oz_k\and  i(z_j z_k-\oz_j\oz_k) }
\eqno(1.13)$$
As Lie algebras, $\kR\simeq\u(n)$, $\k\simeq\gl(n,\C)$,
$\gR\simeq\sp(2n,\R)$ and $\g\simeq\sp(2n,\C)$.

The Hamiltonian flow of $\kR$ on $M=\C^n$ is the natural linear  
representation of  the unitary group $U(n)$. The   Hamiltonian flow of
$\gR$ is the natural linear  representation of  the non-compact
symplectic group $Sp(2n,\R)$. Clearly $U(n)$ is exactly the subgroup of 
$Sp(2n,\R)$ which preserves the Hamiltonian  $F$ in  (1.5). 

We can quantize all the observables in $\g$, in a way consistent with
(1.7)  and $\Q(1)=1$, by
$$\Q(z_jz_k)=z_jz_k,\qquad 
\Q(z_j\oz_k)=z_j\pdb{z_k}+{\delta_{jk}\over 2},\qquad
\Q(\oz_j\oz_k)={\pa^2\phantom{Xi}\over\pa z_j\pa z_k}\eqno(1.14)$$
These operators obey  (1.1) and
$\Q(\phi)^\dagger=\Q(\ophi)$ for $\phi\in\g$.
Moreover this condition by itself determines the inner product
$\<f|g\>$ uniquely. The No-Go Theorem (see e.g., [A-M]) shows that we
cannot extend  the quantization to all polynomial observables.

A benefit of looking at this large Lie algebra of symmetry  is  that
we  can see another  source for the  $\half$-shift in (1.6).
Indeed, the eminently reasonable values of  $\Q(z_jz_k)$ and
$\Q(\oz_j\oz_k)$ in (1.14) imply the value of  $\Q(z_j\oz_k)$ because of 
the Dirac axiom (1.1). So the term involving $\half$ is created exactly
because   $\pdb{z_k}$ and $z_k$ do not commute but   instead
$[\pdb{z_k},z_k]=1$.

 \def\snu{\sqrt{\nu}}

The most convincing  way to understand the $\half$-shift is to introduce
half-forms.  This means that we replace our Hilbert space $\H$ of
holomorphic functions on $\C^n$ by a new Hilbert space $\H'$ of
holomorphic half-forms $s=f\snu$ where $f=f(z_1,\cdots,z_n)$ is still a
holomorphic function and 
$$\nu=\d z_1\wedge\cdots\wedge\d z_n\eqno(1.15)$$ 
is a holomorphic $n$-form.
Then every holomorphic vector  field $\eta$ acts naturally on
half-forms by the  Lie derivative $\L_\eta$  (see \S5). 

On half-forms, $z_j$ and $\oz_k$ quantize into the operators
$$\Q'(z_j)=z_j\AND \Q'(\oz_j)=\L_{\pa_j}\eqno(1.16)$$
Let $\pa_k=\pdb{z_k}$. We compute
$\L_{\pa_k}(\snu)=0$ and 
$\L_{z_j\pa_k}(\snu)=\half\delta_{jk}\snu$.  This gives
$$ \L_{\pa_k}(f\snu)=\ppd{f}{z_k}\snu\AND  
\L_{z_j\pa_k}(f\snu)=
\left(z_j\ppd{f}{z_k}+\half\delta_{jk}f\right)\snu
\eqno(1.17)$$

On half-forms, the observables in $\g$ quantize into the operators
$$\Q'(z_jz_k)=z_jz_k,\qquad 
\Q'(z_j\oz_k)=\L_{z_j\pa_k},\qquad
\Q'(\oz_j\oz_k)=\L_{\pa_j}\L_{\pa_k}\eqno(1.18)$$
These operators in (1.16) and (1.18)  obey (1.1) and 
$\Q'(\phi)^\dagger=\Q'(\ophi)$   where the inner product  
$\<f\snu|g\snu\>$ is again given by the RHS of  (1.10). 
In particular we get
$$\Q(F)=\L_E\quad\hbox{where} E=\sum_{k=1}^nz_k\pa_k
\eqno(1.19)$$
so that $\Q(F)$ is the Lie derivative of the holomorphic Euler vector
field on   $\C^n$.

The operators $i\Q(\phi)$, $\phi\in\g$, give a Lie algebra
representation of $\g$ on $H$ by skew-adjoint operators. THis
integrates to the unitary oscillator representation
$$Mp(2m,\R)\to \Unit L^{2}_{hol}(\C^m)\eqno(1.20)$$ 
where $Mp(2m,\R)$ is the
metaplectic group which doubly covers the symplectic group
$Sp(2m,\R)$. This representation splits into exactly two irreducible
pieces.

There is one more thing we can learn from the oscillator example.
This is that Kaehler polarizations can turn out to be related to 
cotangent bundle geometry. Indeed, we gave no geometric reason for
the  assignments in (1.16) and (1.18).  In quantizing observables on
cotangent bundles $T^*Q$, we have the guiding philosophy that the
principal symbol of $\Q(\phi)$ should be $\phi$ if $\phi$ is 
homogeneous on the fibers of the projection $T^*Q\to Q$. On a Kaehler
manifold we a priori have no notion like this.

However, if $(M,\om)$ is  Kaehler with complex structure $\JJ$, then
we can  ask if   $M$ is a symplectic real form of  the cotangent bundle
$T^*Z$ of some complex manifold $Z$. An obvious choice is  for  $Z$ to 
be   $(M,\JJ)$ (so $Z$ forgets $\om$). Then the ``good" observables on
$M$ would be those that extend to holomorphic (or maybe
rational) functions on  $T^*Z$ which are homogeneous  on the fibers of
$T^*Z\to Z$. The good observables correspond to bona  fide symbols.
See \S2,3 and [B3] for a way to work this out based on the Hamiltonian
$F$. The result of this  is easy to describe directly for the oscillator.
 
We put $Z=\C^n$. Let $\ze_1,...,\ze_n$ be the holomorphic momentum
functions on $T^*Z$ so that $z_1,\dots,z_n,\ze_1,...,\ze_n$ are
holomorphic coordinates on $T^*Z$ and the canonical holomorphic
symplectic form on $T^*Z$ is 
$\Om=\sum_{k=1}^n\d\ze_k\wedge\d z_k$. Then $\Om$ defines a
Poisson bracket   $\{\Phi,\Psi\}_{\Om}$
on the algebra of holomorphic functions on $T^*Z$.  We have
$\{z_j,z_k\}_{\Om}=\{\ze_j,\ze_k\}_{\Om}=0$ and
$\{\ze_j,z_k\}_{\Om}=\delta_{jk}$.

We have an  obvious complex Poisson algebra isomorphism
$$\al:\C[z_1,\dots,z_n,\oz_1,...,\oz_n]\to
\C[z_1,\dots,z_n,\ze_1,...,\ze_n]\eqno(1.21)$$
where $\al(z_k)=z_k$ and $\al(\oz_k)=i\ze_k$.
Then $\al(\phi)$ is the unique extension of $\phi$ to a holomorphic
function $\Phi$ on $T^*Z$ with respect to the  symplectic embedding 
$b$ of $M=\C^n$ into $T^*Z=\C^{2n}$ where $b(w)=(w,\ow)$.
Then
$$\al(z_jz_k)=z_jz_k,\qquad \al(z_j\oz_k)=iz_j\ze_k,\qquad 
\al(\oz_j\oz_k)=-\ze_j\ze_k\eqno(1.22)$$
Now the formulas in (1.16) and (1.18) make sense as $i\ze_k$ is the
symbol of $\pdb{z_k}$.

The quantization  of the oscillator has  manifold
applications  in physics -- in quantum mechanics,
quantum field theory, supersymmetry, etc. It also of course  occupies a
central place in mathematics. 

{\bf III. Quantization of Hamiltonian Symmetry. }
To formulate  a  
mathematical quantization problem generalizing the oscillator case, 
we suppress  (for the time being)  
the Hamiltonian $F$ and focus instead on the large
finite-dimensional symmetry algebra $\g$.  This brings us to  the 
notion of  Hamiltonian symmetry.

Suppose we have an action of a connected Lie group $G$ on a
symplectic manifold $(M,\om)$.  We regard $M$ as a phase space. 
Assume the  action is   {\it symplectic}, i.e.,   $G$  preserves   $\om$. Let
$\g$ be the Lie algebra of   
$G$ . For each $x\in\g$,  we have the  $1$-psg ($1$-parameter
subgroup) $\ga_x:\R\to G$,  $\ga_x(t)=\exp(tx)$, generated by  $x$. 
By Noether's Theorem,  there  is a smooth function $\mu^x$
(defined at least  locally about every point of $M$), unique up addition
of a constant, such that the Hamiltonian flow of $\mu^x$ is the 
action of $\ga_x$.  Then $\mu^x$ is conserved under the action
of $\ga_x$. If   $\mu^x$  exists globally on $M$, then $\mu^x$ is called a
{\it first integral} or  {\it momentum function} for  $\ga_x$.

The symplectic $G$-action is called {\it Hamiltonian} if there exists a 
map $$\mu^*: \g\to\Cinf(M)\eqno(1.23)$$ 
$x\mapsto\mu^x$, such that $\mu^x$
is a first integral for $\ga_x$ and $\{\mu^x,\mu^y\}=\mu^{[x,y]}$
for all $x,y\in\g$, i.e., $\mu^*$ is a Lie algebra homomorphism.
Then the functions $\mu^x$ define a {\it moment map}
$$\mu:M\to\g^*\eqno(1.24)$$   
by $\mu^x(m)=\<\mu(m),x\>$. 
If $\g$ is semisimple then we often identify $\g$ with its
dual by means of the Killing form so that moment maps take values in
$\g$.
 
The moment map $\mu$ obtained in this way is $G$-equivariant and
Poisson. Consequently  the image of  $\mu$ in $\g^*$  is a union of
coadjoint orbits. The image of the moment map is an important
invariant of the action.  It is easy to prove that $\mu$ is a covering onto
a single  coadjoint  orbit if and only if the  Hamiltonian action of $G$ on
$M$ is transitive; then $\mu$ is symplectic.   Such an action is called   
{\it elementary}. 

Thus, symplectically and equivariantly, the elementary
Hamiltonian $G$-spaces are, up to covering, just the coadjoint orbits   
of  $G$.

Going back to our oscillator phase space,  
we see that the action of $Sp(2n,\R)$ on our  manifold
$M=T^*\R^n=\C^n$, with the origin of $\C^n$ deleted, is an
elementary Hamiltonian action.  The moment map
$\C^n-\{0\}\to\sp(2n,\R)$    is a $2$-fold covering on the smallest 
(non-zero)  adjoint orbit $\OR$ of $Sp(2n,\R)$.  This orbit $\OR$ is stable
under scaling and so consists of nilpotent elements.

The quantization problem on the  Hamiltonian $G$-space $(M,\om)$ is
to quantize the momentum functions $\mu^x$ into operators in a
manner agreeable with Dirac's axioms.  It is natural to study the
elementary case first, as here the symmetry is largest. Thus one seeks a
quantization of the functions $\mu^x$, $x\in\g$, for  coadjoint orbits
and their covers.  

In analogy with the oscillator, we consider the  case where the
symmetry group $G$ is a real semisimple Lie group $\GR$ (with finite
center)  and $M$ is  an adjoint orbit $\OR$ stable under
scaling.  Then $\OR$ is a ``nilpotent orbit" of $\GR$ -- see \S2.

Quantization of coadjoint orbits has traditionally been considered
as part of the Orbit Method in representation theory. In the  Orbit
method, one uses  polarizations invariant under the whole
symmetry group and obtains unitary representations by induction.
The theory incorporates metaplectic covers and the Mackey machine.
Much more can be said about the Orbit Method.
We note that unitary representations attached to nilpotent orbits
are called {\it unipotent}   in representation theory.
  
On the other hand, coming into this problem from geometry, we have
found
different methods  which apply (at least) to nilpotent orbits. The main
idea is to transform the quantization problem on $\OR$ into a
quantization problem on a cotangent bundle, and then solve that
problem.

{\bf IV. Outline of this Paper.}

In this paper, we quantize the nilpotent orbit  $\OR$ of  $\GR$  
in the case where $\OR$ is   strongly  minimal (see \S3).
The  oscillator phase space is the double cover of the strongly minimal
nilpotent orbit of $\GR=Sp(2n,\R)$.
 
We assume that $\gR$ is simple, the maximal compact subgroup
$\KR$ of $\GR$ has finite center, and $\GR$ is  simply-connected.  
(Thus we exclude the oscillator case as there $\KR=U(n)$.)

We obtain the  analogs of the Fock space model
of the quantum mechanical oscillator. We find  analogs of all the 
features of the oscillator quantization described above in II. 
This is worked out in detail in this
paper, with the exception of  the integral formula (1.10) for the inner
product which will be written up elsewhere.
We  work from scratch and assume
no prior knowledge on existence of unitary representations.

This completes the work from [B-K4]. In [B-K4]   we worked out
with Kostant the results covered in  \S4-7 of this paper for 
the three  cases where $\GR$ is a 
a  split group of type $E_6$, $E_7$, $E_8$. 

We start from the fact, a product of the work of Kronheimer ([Kr])  and
Vergne ([Ve]),   that  
$\OR$ admits a $\KR$-invariant complex structure
$\JJ$ which together with the KKS symplectic form $\sig$ gives a
(positive) Kaehler structure on $\OR$.  
The  Vergne diffeomorphism $\V:\OR\to Y$ identifies
the complex manifold $(\OR,\JJ)$    with a complex 
homogeneous space $Y$  of the complexification $K$ of $\KR$.
This is a general theory  that applies to every nilpotent orbit
for $\GR$ semisimple.
For   the oscillator, this recovers the $U(n)$-invariant Kaehler
structure and the  identification $T^*\R^n=\C^n$   used in II.

We outline this theory in \S2 and we explain how it gives rise to an 
embedding of $\OR$ into $T^*Y$ as a totally real symplectic
submanifold ([B1]). This enables us to transform the quantization
problem on $\OR$ into a quantization problem on $T^*Y$, as long as  
the Hamiltonian functions $\phi^w$, $w\in\gR$ extend from $\OR$
to $T^*Y$.

An important aspect is that the Kaehler structure on $\OR$ 
possesses a global Kaehler potential $\rho$ which we argue plays the
role of the Hamiltonian $F$.  The  Hamiltonian flow of  $\rho$  
is the action of the center of $\KR$ in the oscillator case.
In our cases, the Hamiltonian flow of  $\rho$ lies outside the
$\GR$-action.

In \S3, we specialize to the case where $\OR$ is strongly minimal
and $\KR$ has finite center. We explain how to
convert the Hamiltonian  functions $\phi^w$,
$w\in\gR$,  on $\OR$ into rational meromorphic functions $\Phi^w$ on
the  cotangent bundle of  $Y$. We interpret the $\Phi^w$ as
``pseudo-differential symbols". 

To describe the symbols, we  consider the Cartan decomposition
$\gR=\kR\oplus\pR$  (cf. (1.13)).
For $x\in\kR$, $\Phi^x$ is just the usual
symbol of the  holomorphic vector field $\eta^x$ on $Y$ defined by
differentiating the $K$-action. But for $v\in\pR$, $\Phi^v$ is a sum of
two terms, each homogeneous under the fiberwise scaling action of
$\C^*$ on the leaves of the cotangent polarization of $T^*Y$.
The passage from the observable function $\phi^w$ to the
symbol $\Phi^w$ preserves Poisson brackets.

The middle part \S4-\S7 of the paper is devoted to  quantizing  the 
symbols $\Phi^w$, $w\in\gR$, into skew-adjoint operators on a 
holomorphic half-form line  bundle $\NH$ over $Y$.  
In  \S5, we construct all such bundles. We find the 
space $H$ of  global algebraic holomorphic sections of $\NH$  
is a     multiplicity free ladder representation of $K$.  We get a simple
geometric description  of the sections which are  the highest weight
vectors. 

In \S4, we set up the Jordan structure that is used throughout 
the  paper (explicity in \S5 and \S7).  A main point is that the
polynomial function $P$ constructed in \S3  is realized in terms of
Jordan norms.

We construct, in  Corollary 6.2 and Theorem 6.3
the  pseudo-differential operators 
$\Q(\Phi^w)$ on half-forms  which  quantize the  symbols $\Phi^w$,  or
equivalently, the functions $\phi^w$.   Theorem 6.3 says that
these operators satisfy (1.1), i.e.,  the operators
$\pi^w=i\Q(\Phi^w)$ give  a representation of $\g$.
In Theorem 6.6 we construct  the $\gR$-invariant
inner product $B$ on $H$.  In
Theorem 6.8 we compute $B$ by giving the analog (6.30) of (1.9).   

 Our operators are pseudo-differential
(not purely differential) in that they involve inverting the
positive-spectrum ``energy" operator $E'$
which is the quantization of $\rho$.  In fact, instead of the
order two operators $\L_{\pa_j}\L_{\pa_k}$  
 from (1.18) we obtain  order $4$ differential operators
divided by $E'(E'+1)$; these are ``formally" of order $2$.
The action of the maximal compact group $\KR$ on $H$ is 
just the natural one defined by the action of  $\KR$ on $Y$
and $\NH$.   

Theorem 6.4 says that  our representation $\pi$  of $\g$ on $H$ is
irreducible. Also we describe the algebra   
generated by the operators $\pi^w$ on $H$.  
It  follows in Theorem 6.6 that $\pi$ integrates to give an irreducible 
minimal unitary representation of $\GR$ on the Hilbert space
completion $\H$ of $H$.

Next \S7 is devoted to proving   the results of  \S6. 
We show that our  pseudo-differential operators  
satisfy the bracket relations of $\gR$ by  reformulating the problem
and applying  the generalized Capelli Identity of Kostant and Sahi 
([K-S]). An important aspect of their work is that
Jordan algebras provide a natural setting for generalizing
the classical Capelli identity involving square matrices. 
The complex Jordan algebra $\k_{-1}$ occurring here is  semisimple
(while in [B-K4] it was simple). It turns out  that the simple   
components of   $\k_{-1}$ become coupled together  in our 
calculations in a  subtle way reflected by  Proposition 7.8.

In \S8, we compute the reproducing kernel $\K$ of  the Hilbert
space completion $\H$ of $H$. We find that $\K$ is a holomorphic
function on $\YoY$ and hence $\H$ consists entirely of holomorphic
sections of $\NH$.  
Finally, in \S9 we give some examples.

Different models, or proofs of existence, 
for  most of the unitary  representations we construct
have been obtained by other authors. These include
Binegar,   Gross, Howe, Kazhdan,  Kostant, Li, Oersted,
Rawnsley, Savin, Sijacki, Sternberg, Sabourin, Torasso,
Vogan, Wallach, Wolf, and Zierau.
Moreover in [T], Torasso  constructs in a uniform manner
by  the Orbit Method  Schroedinger type models of all minimal unitary 
representations. Precisely, Torasso  constructs  unitary irreducible
representations attached to  all minimal admissible nilpotent orbits
of simple groups of relative rank at least three over a local field
of zero characteristic. It would be very interesting to construct
intertwining operators between our models.

There is a rich literature on geometric
models of unitary  highest weight representations, and there are
many interesting ties here with our work.

This  paper builds on several years of  joint work with Bert Kostant on
the algebraic   holomorphic symplectic geometry of    nilpotent
orbits of a complex semisimple Lie group. This   work  includes
[B-K1-5].  In addition \S4  of this paper is  joint work.

I   thank Alex Astashkevich, Olivier Biquard, Murat Gunaydin, Bert
Kostant, Michele Vergne, and Francois Ziegler  for  useful
conversations relating to this work.
Parts of this work were carried out during visits to Harvard
(1993-94, summers of 1995 and 1996), the Institute for Advanced Study
(Spring 1995) and Brown University (summer  of 1997). I thank all these
departments   for their hospitality.
I thank Mark Gotay for putting together this volume and for his comments
on my paper.

I am delighted to dedicate this paper to Victor Guillemin and to be
able to contribute it to this volume in his honor.  In my graduate
student days at MIT I was ensconced in algebraic  geometry and
algebraic group actions.  I was symplectically agnostic.   But since
my symplectic conversion in the end of the last decade, I have had
the opportunity to talk to Victor a lot and learn from him and his
many books and  papers.  I thank him for warmly welcoming me as a
visitor  into his symplectic group.

\Sec{\S2.   The Quantization Problem for  Real Nilpotent Orbits.}

The phase spaces we wish to   quantize are the so-called 
``nilpotent orbits" of  $\GR$ where $\GR$ is a connected   
non-compact real semisimple Lie group with finite center.
Then $\GR$ is a finite cover of   the adjoint group of its  Lie
algebra $\gR$, and $\gR$ is semisimple. To define the nilpotent orbits
we consider the coadjoint action of $\GR$
on the dual $\gR^*$ of $\gR$.  

Each coadjoint orbit  $\OR$
carries   a natural $\GR$-invariant symplectic form $\sig$, often called 
the KKS or Lie-Poisson form. The form $\sig$  is uniquely characterized  by
the following property: let
$$\phi:\gR\to\Cinf(\OR),\qquad w\mapsto\phi^w\eqno(2.1)$$
be the   pullback map on functions defined by the embedding
$\OR\subset\gR^*$. Then $\phi$ is a Lie algebra  homomorphism 
with respect to Poisson bracket on $\Cinf(\OR)$ defined by $\sig$.

In analogy with the cotangent bundle, we wish to single out those
coadjoint orbits  which are {\it conical} in the sense that they are stable
under the Euler scaling action of
$\R^+$ (positive reals). There is a nice Lie theoretic characterization of
these orbits.  To get this, we first use the Killing form $(\,,\,)_{\gR}$
to identify $\gR^*$ with  $\gR$; we do this throughout the paper routinely.
Then (conical) coadjoint orbits get identified with (conical) adjoint orbits. 

An adjoint orbit is conical  if and
only if it consists of nilpotent elements in $\gR$. Such orbits are called
``nilpotent orbits".  It is well-known in Lie theory that there are only
finitely many nilpotent orbits in $\gR$.  

From now on,  we take $\OR$ to be a nilpotent orbit in $\gR$.
The quantization problem on $\OR$ is to quantize into operators the
functions $\phi^w$, $w\in\gR$.  This is a reasonable goal. Ideally 
quantization would convert all smooth functions on $\OR$ into operators
in a manner satisfying  Dirac's  axioms. 
See, e.g., [Ki,\S2.1] for a complete axiom list. 
But full quantization  is  impossible even for
polynomial functions on $\R^2$ (the infamous No-Go Theorem).
We are left hoping that, except for anomalies, finite-dimensional 
Hamiltonian symmetry  will quantize.

In analogy with  the Fock space quantization of the 
oscillator, we look for a  Kaehler polarization of  our phase space 
$(\OR,\sig)$ which is invariant under a fixed maximal compact subgroup
$\KR$ of $\GR$.  This means that we look for a $\KR$-invariant integrable
complex structure $\JJ$ on $\OR$ such that $\JJ$ and $\sig$ together give
a (positive) Kaehler structure on $\OR$.

Fortunately,   such a complex structure $\JJ$ on $\OR$  arises
from the works of Kronheimer ([Kr]) and Vergne ([Ve]) on  instantons
and nilpotent orbits.  
This gives the  $\KR$-invariant {\it instanton Kaehler structure} 
$(\JJ,\sig)$ on $\OR$.   This  structure is discussed and studied in 
detail in [B1].  We recall two main points.  

The first point  is the Vergne diffeomorphism ([Ve]).
To set this up, we introduce the Cartan decomposition
$$\gR=\kR\oplus\pR\eqno(2.2)$$
where $\kR\subset\gR$ is the Lie algebra of $\KR$ and $\pR$ is its
orthogonal complement with respect to the Killing form.
The natural  action of $\KR$ on $\pR$ complexifies to a complex  
algebraic action of  $K$ on $\p$ where $K$ is the complexification of 
$\KR$ (so that $K$ is a complex reductive algebraic group) and 
$\p=\pR\oplus i\pR$. 

Now the Vergne diffeomorphism
$$\V:\OR\to Y\eqno(2.3)$$
is a   $(\KR\times\RP)$-equivariant diffeomorphism of
real manifolds   which maps $\OR$ onto a  $K$-orbit $Y$ in 
$\p$. $Y$, being a $K$-orbit,  is  manifestly a complex submanifold of
$\p$. Moreover $\JJ$
is the pullback through $\V$ of the complex structure on $Y$.

An important feature is that $Y$ is stable under the Euler scaling action 
of  $\C^*$ on $\p$. This follows since $\OR$ is $\RP$-stable and $\V$ is
$\RP$-equivariant.
Let $E$ be the infinitesimal generator of  the Euler $\C^*$-action so that
$E$ is the  algebraic holomorphic  Euler vector field on $Y$.

 In general, the target $Y$ of the Vergne diffeomorphism is known
(by the Kostant-Sekiguchi correspondence [Sek]) but not the actual map
giving $\V$. A little insight into $\V$ comes from Lie theory.

To explain this, we introduce the complexified Lie
algebra $\g=\gR\oplus i\gR$. 
$\g$ is a complex semisimple Lie algebra and   carries the complex
conjugation map $x+iy\mapsto\ovl{x+iy}=x-iy$.

An {\it S-triple} in $\g$ is a basis $(e,h,f)$  of a 
subalgebra isomorphic to $\sl(2,\C)$ which  satisfies the bracket
relations 
$[h,e]=2e$, $[h,f]=-2f$, $[e,f]=h$.  The S-triple is {\it adapted}
to $(\gR,\kR)$ if  $e$ and $f$ are complex conjugates and
$h\in i\kR$.
Given $\OR$, we can find an S-triple $(e,h,\oe)$ adapted to
$(\gR,\kR)$ such that $e+ih+\oe$ lies in $\OR$. Then Vergne's  
construction gives
$$\V(e+ih+\oe)=e$$

The second point is that  the Kaehler structure
$(\JJ,\sig)$ on $\OR$  admits a global Kaehler potential $\rho$. This
means that $\rho$ is a smooth real valued function on $\OR$ such
that $i\del\odel\rho=\om$.  Moreover $\rho$ is uniquely
determined   by the added condition that $\rho$ transforms
homogeneously under the Euler
$\RP$-action on $\OR$. Then $\rho$ is $\KR$-invariant and 
Euler homogeneous of degree $1$.  

Next we examine how to use this Kaehler structure in quantization.
The Vergne diffeomorphism identifies $Y$ as  $\OR$ equipped with
a complex polarization. The philosophy of Geometric  Quantization 
now predicts that we can quantize suitably nice
real-valued functions $\phi$ on
$\OR$ into self-adjoint operators on a Hilbert space consisting of
holomorphic sections of a suitable holomorphic vector bundle over $Y$. 

Our quantization program for $\OR$   becomes:  ``quantize"
each function $\phi^w$, $w\in\gR$, into a self-adjoint operator  
$\Q(\phi^w)$ on a Hilbert space $\H$ of  square integrable holomorphic
sections   of a holomorphic half-form  complex line
bundle $\NH$ over $Y$ in such a way that the Dirac axiom
$$\Q(\phi^{[w,w']})=i[\Q(\phi^w),\Q(\phi^{w'})]\eqno(2.4)$$ 
is satisfied. In the course of  doing this, we will end up quantizing  one
additional function on $\OR$.

There are additional axioms  which should also be satisfied, but these
are somewhat hidden as we are only dealing with the functions
$\phi^w$. E.g., the axiom that the constant function $1$ quantizes to the
identity operator is ``hidden". These ``hidden axioms" are basically 
incorporated by our methodology developed below using symbols.

If  the Hamiltonian flow of $\phi$ preserves $\JJ$ and $\phi$ is
homogeneous of degree $1$, then we    mandate that the
quantized operator is simply 
$$\Q(\phi)=-i\L_{\widehat{\xi_\phi}}\eqno(2.5)$$
Here $\widehat{\xi_\phi}$ is the {\it $\JJ$-Hamiltonian vector
field}  on $Y$ defined by the condition that $\widehat{\xi_\phi}$
is   holomorphic and coincides with $\xi_{\phi}$
on  holomorphic  functions.  We write $\L_{\eta}$ for the Lie derivative
operator (acting  on holomorphic half-forms)
with respect to a holomorphic vector field $\eta$.   

Differentiating the $K$-action on $Y$ we get an infinitesimal
holomorphic vector field action 
$$\k\to\Vect^{\hol}\;Y,\qquad x\mapsto\eta^x\eqno(2.6)$$
Then $\eta^x=\widehat{\xi_{\phi^x}}$  for $x\in\kR$ and so 
$$\Q(\phi^x)=-i\L_{\eta^x},\qquad\hbox{for }x\in\kR\eqno(2.7)$$ 

The problem, since our polarization $\JJ$ is only $\KR$-invariant, is to
quantize the  remaining functions $\phi^v$, $v\in\pR$ corresponding
to the second piece in the Cartan decomposition (2.2).

A key aspect of our program for quantization of real nilpotent orbits
(see [B1-3]) is that we regard $\rho$ as the Hamiltonian function on
$\OR$. This generalizes the case of the harmonic oscillator
discussed in \S1 where the Hamiltonian is the total energy. 
It may seem strange  that the oscillator energy Hamiltonian is
homogeneous quadratic while our function  $\rho$  is homogeneous 
linear. However  the oscillator phase space
$\R^{2n}-\{0\}$ arises as the {\it double} cover of a real nilpotent
orbit. In that case, our  linear potential function $\rho$ does indeed
pull back to a quadratic function   on $\R^{2n}-\{0\}$, and it is  easy
to check that we recover the classical energy 
$p_1^2+q_1^2+\cdots+p_n^2+q_n^2$ (see [B3]).

In physical terms, the Hamiltonian governs the time evolution of the
classical system. The quantum mechanical problem is
to find the eigenvalues and eigenstates of  the   operator
quantizing the Hamiltonian. 

Thus we now demand that 
quantization should not only promote  the symmetry functions
$\phi^w$ to operators, but  should  also promote $\rho$   to
an operator.  In fact the Hamiltonian flow of $\rho$ preserves $\JJ$
and is periodic; we call this the $KV$ (Kronheimer-Vergne)   
$S^1$-action on $\OR$ ([B1]).  Under $\V$, the $KV$    
$S^1$-action corresponds to the 
circle part of the Euler $\C^*$-action on $Y$. It follows  that 
the $\JJ$-Hamiltonian vector field of $\rho$ is $iE$. Hence
$$\Q(\rho)=-i\L_{iE}=\L_E\eqno(2.8)$$

Let $\Om$ be the canonical holomorphic symplectic form on $T^*Y$.
Then $\Om$ defines a   Poisson bracket on the algebra of 
holomorphic functions on $T^*Y$, and also on the field of
meromorphic functions.

A main result  of  [B1]  is   to realize the holomorphic
cotangent bundle $(T^*Y,\Om)$ as a symplectic  complexification of
$\OR$.  To do this, we push forward $\rho$ to a smooth function
$\rho_Y$  on $Y$ so that $\rho=\rho_Y\circ\V$.   Next we construct
the following real $1$-form $\be$ on $Y$ 
$$\be=-{i\over 2}(\del-\odel)\rho_Y\eqno(2.9)$$
Then $\be$  defines a smooth section of the cotangent bundle
$T^*Y\to Y$.

\theo{2.1[B1]}{The composition
$$b:\OR\;\mapright{\V}\;Y\;\mapright{\be} \; T^*Y\eqno(2.10)$$
embeds $\OR$ as a totally real symplectic submanifold of $T^*Y$.
In particular, $b^*(\Re\Om)=\sig$ and $b^*(\Im\Om)=0$.}

Now, given a function $\phi$ on $\OR$ which we wish to quantize, we 
can ask if $\phi$ extends to a holomorphic function $\Phi$ on $T^*Y$.
(Such an extension, if it exists, is necessarily unique.)
If so, then $\Phi$ is our candidate for the symbol of  $\Q(\phi)$.

This philosophy is consistent with what we already found in 
(2.7) and (2.8).  Indeed we can define the holomorphic symbols,
where $x\in\k$,
$$\Phi^x=\symbol\eta^x  \AND \la=\symbol E\eqno(2.11)$$

Our convention for symbols is specified by the following formula in
holomorphic Darboux coordinates:
$$\symbol \quad f(z_0,\dots,z_m)
{\pa^{k_0+\cdots+ k_m}\phantom{xx}
\over \pa z_0^{k_0}\cdots\pa z_m^{k_m}}=
 f(z_0,\dots,z_m)i^{k_0+\cdots+ k_m}\ze_0^{k_0}\cdots\ze_m^{k_m}  $$

It is easy to check ([B1,3]) 
\cor{2.2}{
\item{\rm(i)} The Kaehler potential $\rho$ on $\OR$ extends
uniquely to a holomorphic function on $T^*Y$.  Precisely,
$\rho$   extends to $\la=\symbol \Q(\rho)$.
\item{\rm(ii)} For $x\in\kR$,  $\phi^x$ extends
uniquely to a holomorphic function on $T^*Y$. Precisely,
$\phi^x$ extends to $\Phi^x=\symbol  \Q(\phi^x)$.  \mypar }
\noindent In effect, $\be$ was engineered to make (i) true. 

The passage from functions on $\OR$ to holomorphic functions on
$T^*Y$ preserves Poisson brackets. I.e., 
if $\Phi_1$ and $\Phi_2$ are respectively the
holomorphic extensions  of two real functions $\phi_1$ and $\phi_2$
on $\OR$ then
$\{\Phi_1,\Phi_2\}_{\Om}$ is the holomorphic extension
of $\{\phi_1,\phi_2\}_{\sig}$,
where the subscripts indicate the symplectic form defining the
Poisson brackets. This follows easily from Theorem 2.1.

Thus  if all the Hamiltonian functions $\phi^w$, $w\in\gR$, extend
holomorphically  from $\OR$ to $T^*Y$, then this
in effect converts our  quantization problem on $\OR$ into a
holomorphic  quantization problem on $T^*Y$.

This brings us to the question as to whether the Hamiltonian functions
$\phi^v$, $v\in\pR$, extend to holomorphic functions on $T^*Y$.
The general answer is no.
However, the better question is whether the $\phi^v$ extend to
{\it  meromorphic} functions  $\Phi^v$ on $T^*Y$.
In [B3], we show that the answer is yes in every case, at least if we
allow  the $\Phi^v$ to lie in a finite   extension of  the field of 
meromorphic  functions  on $T^*Y$. This relies on the powerful
result of Biquard [Bi1, Bi2] that the homogeneous hyperkaehler
potential on a complex nilpotent orbit is always a positive Nash
function.

In fact, the symbols  that arise here are
all rational functions on $T^*Y$ (or at least regular functions on an
\'etale cover of a Zariski open set of $T^*Y$) in the sense of algebraic
geometry.  The holomorphic symplectic form
$\Om$ is manifestly  algebraic and so
$\Om$ defines a Poisson bracket on the algebra $R(T^*Y)$ of
algebraic holomorphic functions on $T^*Y$ and also on the field
$\C(T^*Y)$ of rational functions on $T^*Y$.

In the next section, we explain in detail how this works for the
smallest orbits.

\Sec{\S3.  Pseudo-Differential Symbol Realization  of $\gR$ for $\OR$
Strongly Minimal.} 

Each nilpotent orbit $\OR\subset\gR$ lies in a unique complex
adjoint orbit $O\subset\g$.  Then $O$ is a complex nilpotent orbit
(i.e., $O$ consists of nilpotent elements in $\g$). We call $O$ the
{\it complexification} of $\OR$. The nilpotent elements  in   $\g$ are 
characterized by the property that their adjoint orbits are stable
under the scaling action of   $\C^*$.

We assume from now on that the complex Lie algebra $\g$ is simple.
Let $G$ be the adjoint group of $\g$. Then $G$ is a 
connected complex semisimple 
algebraic group with Lie algebra $\g$ and  complex conjugation on 
$\g$ defines  a complex conjugation map $g\mapsto\og$ on $G$.
Let  $(\,,\,)_{\g}$ be the complex Killing form of $\g$.
We often identify $\g$ with  $\g^*$ and 
$\k$ with $\k^*$ by means of  $(\,,\,)_{\g}$.

Recall from \S2 that $\KR\subset\GR$ is a fixed maximal compact 
subgroup with complexification $K$. We have natural maps
$\GR\to G$ and $K\to G$  and both maps have finite kernel.
 
Since $\g$ is a simple Lie algebra, the adjoint representation of $G$ on
$\g$ is irreducible. The orbit $\Omin$ of highest weight vectors is then 
nilpotent, as it is the orbit of a highest root vector.
Moreover, $\Omin$ is minimal among all non-zero nilpotent orbits   in the
sense that it lies in the closure of every non-zero nilpotent   orbit.
It follows that $\Omin$ is the {\it unique} (non-zero) minimal nilpotent
orbit.

We will call a real nilpotent orbit 
$\OR$  {\it strongly minimal} if the complexification of $\OR$ is
$\Omin$.  In  Theorem 4.1  below  we recall from [B-K5] the
classification of strongly minimal real nilpotent orbits.
For $\OR$ strongly minimal, formulas  for  $\V$ and
$\rho$ are easy to write down because  the action of $\KR\times\RP$ is
transitive on
$\OR$.  (However, for general $\OR$ the action is not transitive, and
working out $\V$ and $\rho$  is a hard open problem.)
  
As $\g$ is simple, 
there are just two possibilities for the  center of $\KR$: either 
(i) $\Cent \KR$ is a circle subgroup  or  (ii)  $\Cent \KR$ is finite. 
These cases correspond exactly to the nature of the irreducible 
symmetric space $\GR/\KR$, so that $\GR/\KR$ is
Hermitian in (i) and non-Hermitian in (ii). Accordingly,  
we call the  complex symmetric pair  $(\g,\k)$ 
Hermitian or non-Hermitian.

For each $v\in\p$, let $f_v$ be the linear function on $\p$ defined by
$f_v(u)=(v,u)_{\g}$. Then by restriction to $Y$ we get a
$K$-equivariant complex linear map 
$$\p\to R(Y),\qquad v\mapsto f_v\eqno(3.1)$$ 
Every  algebraic holomorphic function on $Y$ defines an 
algebraic holomorphic  function on
$T^*Y$ by pullback through the projection $T^*Y\to Y$.

From now on in \S3, we assume that $\OR$ is strongly minimal and
the center of  $\KR$ is finite.

\theo{3.1}{ Let $v\in\pR$ and $x\in\kR$.
\item{\rm(i)} Recall the embedding $b:\OR\to T^*Y$ from  {\rm(2.10)}.
Each function $\phi^v$ on $\OR$ extends uniquely to a  
rational  function $\Phi^v$ on $T^*Y$.  Set $\Phi^{x+v}=\Phi^x+\Phi^v$
where    $\Phi^x$ was defined  in {\rm (2.11)}.
The resulting linear map
$$\gR\to\C(T^*Y),\qquad w\mapsto\Phi^w\eqno(3.2)$$
is a $1$-to-$1$ real Lie algebra homomorphism with respect to the
Poisson bracket on $\C(T^*Y)$ defined by $\Om$.
\item{\rm(ii)} Each rational function $\Phi^v$   is everywhere defined
on the Zariski open dense complex algebraic submanifold
$$M=\{m\in T^*Y\,|\, \la(m)\neq 0\}\eqno(3.3)$$ 
so that $\Phi^v$ is algebraic holomorphic on $M$.
\item{\rm(iii)} We have  
$$\Phi^v=f_v+g_v\eqno(3.4)$$
where   $g_v$ is an algebraic holomorphic function on $M$
which is homogeneous of degree
$2$ with respect to the  Euler $\C^*$-action on the fibers
of the the projection $M\hookrightarrow T^*Y\to Y$.
\mypar }
 
\proof{This is proven in a more general setting in [A-B1].\qed}

We write $R(X)$ for the algebra of algebraic holomorphic functions on a
complex algebraic variety $X$.
Recall from (2.11) that $\la\in R(T^*Y)$ is the symbol of the Euler 
vector field.

\lem{3.2}{We have $R(M)=R(T^*Y)[\la\i]$.}
\proof{This follows easily since $M$ is the  complement of  the
irreducible divisor $(\la=0)$ in the smooth (and hence normal)
variety $T^*Y$.\qed }

It is  natural now to extend (3.2) $\C$-linearly so that
$\Phi^{x+iy}=\Phi^x+i\Phi^y$ for $x,y\in\gR$. This is consistent with 
(2.11). We have the complexified Cartan decomposition 
$$\g=\k\oplus\p\eqno(3.5)$$

\cor{3.3}{The map {\rm(3.2)}  extends to     a $1$-to-$1$
complex Lie algebra homomorphism
$$\g\to R(T^*Y)[\la\i],\qquad z\mapsto \Phi^z\eqno(3.6)$$
Then for $v\in\p$ we have again the same formula {\rm(3.4)}.}

The significance of Corollary 3.3 is that we can regard functions in
$R(T^*Y)[\la\i]$ as ``pseudo-differential" symbols; cf. \S6.

We have  now, in Theorem 3.1 and
Corollary 3.3, transformed our original problem of quantizing the 
functions $\phi^w$, $w\in\gR$, on $\OR$ into the problem of 
quantizing the rational functions $\Phi^w$, $w\in\gR$, on $T^*Y$.  We
mandate
$$\Q(\Phi^w)=\Q(\phi^w)$$
The new problem lies in the holomorphic symplectic  category: the 
problem is to quantize each $\Phi^w$ into a self-adjoint operator 
$\Q(\Phi^w)$ on a Hilbert space consisting of holomorphic sections of a
holomorphic half-form bundle on  $Y$.  

The  advantage  of the new problem is that $\Phi^w$ is already a 
symbol, and so we  can try to quantize it by constructing reasonable
quotients of differential operators with symbol $\Phi^w$. 
We emphasize that (3.4) says that  $\Phi^v$, $v\in\pR$, is 
{\it not} a principal symbol, but instead is a sum of two principal symbols
$f_v$ and $g_v$. We will get around this by a  naive trick: we will
quantize $f_v$ and $g_v$ separately and then add the answers.
 
Since $f_v$ is just a holomorphic function on $Y$, we mandate that the
quantization of $f_v$ is $\Q(f_v)=f_v$, i.e.,  $\Q(f_v)$ is the operator
defined by multiplication by $f_v$.

The aim of the  rest of this section is to state a formula 
for the symbols $g_v$.
We want  to express $g_v$   in terms 
of the basic symbols $f_v$, $v\in\p$, $\Phi^x$, $x\in\k$, and
$\la$  since we already know how to quantize these symbols.
To work this out, we construct a set of local (\'etale) coordinates on
$T^*Y$ consisting of   basic symbols  in Lemma 3.5 below.

We  begin by setting up  some of the 
Lie theoretic structure associated to $\OR$ following  [B-K4,\S2].  We  
 will make use of this throughout the paper.   
We note that the  discussion of $\Omin$ in [B-K4,\S2] was in the same
generality we have here, and it was only from \S3 onwards in that
paper that the work  specialized to the three cases where  $\GR$ is
split of type $E_6,E_7$ or $E_8$.

To begin with   we have
$$\OR=\Omin\cap\gR=\GR\cdot(e+ih+\oe)\AND 
Y=\Omin\cap\p=K\cdot e\eqno(3.7)$$ 
where $(e,h,\oe)$ are  chosen as in \S2. Then  
$$\s=\C e\oplus\C h\oplus\C\oe\eqno(3.8)$$ 
is the corresponding $\sl(2,\C)$-subalgebra. We assume from now
on that $(\,,\,)_\g$ is rescaled so that $(e,\oe)_\g=1$.

The   action of $\adh$ on $\g$ is diagonalizable with 
spectrum $\{\pm 2,\pm 1, 0\}$ so that we have the $5$-grading
$$\g=\g_{2}\oplus\g_{1}\oplus\g_0\oplus\g_{-1}\oplus\g_{-2}
\eqno(3.9)$$
where the subscripts indicate the  corresponding eigenvalues.
Then $\g_s=\k_s\oplus\p_s$ and 
$$\k=\k_{1}\oplus\k_0\oplus\k_{-1}\AND 
\p=\p_{2}\oplus\p_{1}\oplus\p_0\oplus\p_{-1}\oplus\p_{-2}
\eqno(3.10)$$
Clearly then   $\k_{\pm 1}$ and $\p_{\pm 1}$  are   abelian Lie 
subalgebras. We recall from [B-K4,\S2.2.-2.4]

\lem{3.4}{The spaces $\g_{\pm 2}$ are $1$-dimensional with
$$\g_2=\p_2=\C e\AND  \g_{-2}=\p_{-2}=\C\oe\eqno(3.11)$$
We have $\dimC\p_s=\dimC\p_{-s}$ and $\dimC\k_s=\dimC\k_{-s}$.
The Lie bracket defines a perfect pairing $\k_1\times\p_1\to\C e$. 
Thus we may define
$$m=\dimC\p_{\pm1}=\dimC\k_{\pm 1}\eqno(3.12)$$
\mypar
The subspace $\g_2\oplus\g_1$ is a $(2m+1)$-dimensional Heisenberg
Lie algebra   with center $\g_2$. 
We have $\g=\g^\oe\oplus\C h\oplus\g_{2}\oplus\g_{1}$. 
Consequently
$$\dimC\Omin=2m+2\eqno(3.13)$$}

In particular 
$$\dimC~Y={1\over 2}\dimC\Omin=m+1\eqno(3.14)$$

Now we take a basis   $v_1,\dots,v_m$  of $\p_1$. 
We put $v_0=e$. The corresponding regular functions on $Y$
defined by (3.1) are
$$\fo=f_{v_0}, \qquad \fa=f_{v_1}, \quad\dots,\quad\fm=f_{v_m}
\eqno(3.15)$$ 
These form a system  of local coordinates on $Y$
by ([B-K3, Prop. 5.2]). In fact we get an isomorphism of varieties  
$$\Yo\to\C^*\times\C^m, \qquad
y\mapsto (\fo(y),\fa(y),\dots,\fm(y))\eqno(3.16)$$
where $\Yo\subset Y$ is the open set  given by
$$\Yo=(\fo\neq 0)\eqno(3.17)$$
In our local coordinates we have
$$\eta^x=\sum_{k=0}^m(\eta^xf_k)\pdb{f_k}\eqno(3.18)$$

Let $\{x_1,\dots,x_m\}$ be the basis of  $\k_1$ such that
$[x_i,v_j]=\delta_{ij}e$. Then   in terms of our local coordinates 
$\fo,\fa,\dots,\fm$   the expressions  for  our vector fields $E$ and 
$\eta^{x_i}$, $i=1,\dots,m$ are
$$ E=\sum_{i=0}^m f_i\pdb{f_i}\AND
\eta^{x_i}=\fo\pdb{f_i} \eqno(3.19)$$
It  follows easily from these formulas that 
\lem{3.5}{The $2m+2$ functions 
$f_0,\fa,\dots,\fm,\la,\Phi^{x_1},\dots,\Phi^{x_m}$ form a system of
local \'etale coordinates on $T^*Y$.}

The single function $g_{v_0}$ determines all  the functions $g_v$
because of the $K$-action; indeed,
$g_{[x,v]}=\{\Phi^x,g_v\}_{\Om}$ for $x\in\k$. To state our formula for
$g_{v_0}$ we need one more ingredient, the polynomial function
function $P$ defined below.

The vector fields $\eta^x$, $x\in\k$ define a natural complex 
algebra homomorphism from the universal enveloping algebra     
$\U(\k)$ to the algebra $\D(Y)$ of  algebraic holomorphic
differential operators on $Y$. So   in particular  we get a
representation
$$\pi_K:\U(\k)\to\End\, R(Y)\eqno(3.20)$$
On the symbol level, (3.20) corresponds to the graded Poisson
algebra homomorphism
$$\Phi_K:S(\k)\to R(T^*Y)\eqno(3.21)$$
defined by $\Phi_K(x)=\Phi^x$ for $x\in\k$.

The adjoint action of $\g$ defines a a complex algebra
homomorphism $\ad: \U(\k)\to\End\p$, 
$Q\mapsto\ad\,Q=\ad_Q$.
Let $P$ be the polynomial function on $\k_{-1}$  defined by
$${1\over 4!}\ad_y^4(e)=P(y)e\eqno(3.22)$$
where $y\in\k_{-1}$. Then $P$ is  homogeneous of degree $4$.

We have a perfect pairing 
$$\k_{1}\times\k_{-1}\to\C\eqno(3.23)$$ defined
by the Killing form $(\,,\,)_\g$ as in [B-K4,\S2.5].
This gives an identification of  $S(\k_1)$ with the algebra of polynomial
functions on $\k_{-1}$. This identification places   $P\in S^4(\k_1)$ 
so that  we may write
$$P=P(x_1,\dots,x_m)\AND
\Phi_KP=P(\Phi^{x_1},\dots,\Phi^{x_m})\eqno(3.24)$$

\theo{3.6}{The function $g_{v}$ in Theorem {\rm 3.1(iii)}
is given  for $v=v_0=e$  by  
$$g_{v_0}=-{1\over\la^2}{\Phi_KP\over f_0}
\eqno(3.25)$$}
\proof{A   more general result is proven  in [A-B2].\qed}

In the next section, we set up the Jordan algebra machinery which 
gives us a useful and computable way to understand the polynomial $P$.
In \S5, we already use this machinery to classify half-form bundles 
on $Y$. The reader eager to see how we   quantize the symbols $g_v$
and then the symbols $\Phi^w$ can skip ahead to Lemma 5.3 
and Proposition 5.5 and then to \S6.  

\Sec{\S4.  Complex Minimal Nilpotent Orbits and Jordan Algebras.}

Our first aim in this section is to classify the real simple Lie algebras 
$\gR$ which possess a strongly minimal real nilpotent orbit. 
This amounts to classifying $\gR$ such that
$\Omin$ has real points because of (3.7).
This classification, recalled in Theorem 4.1 below, uses the
geometry of $\Omin$.
 
Any complex nilpotent orbit $O$, and so in particular $\Omin$,
is a quasi-affine smooth locally closed complex algebraic subvariety
in $\g$. This follows since the adjoint action of $G$ on $\g$ is complex
algebraic. Furthermore, $O$  is  an
algebraic holomorphic symplectic manifold   with respect to its
$G$-invariant  holomorphic KKS symplectic form 
$\Sig$ (cf. [B-K1]).    The $G$-action on $O$ is Hamiltonian with
holomorphic moment map given by the embedding  $O\subset\g$. 

Let $\mu:\g\to\k$ be the projection defined by (3.5). Then
the composite map  $\mu: \Omin\to\g\to\k$ is  the moment map for
the  Hamiltonian $K$-action on $\Omin$.  
Let $\cal N(\k)$ be the cone of nilpotent elements in $\k$.

\theo{4.1[B-K5]}{
The following conditions are equivalent:
\mypar {\rm (i)} $\Omin\cap\gR$ is empty
\mypar {\rm (ii)} $K$ has a Zariski open  orbit on $\Omin$
\mypar {\rm (iii)} $\mu(\Omin)\subset{\cal N}(\k)$ 
\mypar
\noindent and  imply that the principal isotropy group of $K$ on 
$\Omin$   is  $\Ks$ where $\s$ was defined in {\rm(3.8}.
\mypar
 The   complete list of all complex symmetric pairs $(\g,\k)$ 
{\rm (}with $\g$   simple {\rm )} which satisfy {\rm(i)-(iii)} is:
\mypar {\rm (a)} $(\sl(2n,\C),\sp(2n,\C))$, where  $n\ge 2$ 
\mypar {\rm (b)} $(\so(p+1,\C),\so(p,\C))$, where $p\ge 3$
\mypar {\rm (c)} $(\sp(2p+2q,\C),\sp(2p,\C)+\sp(2q,\C))$,  
where $p,q\ge 1$
\mypar  {\rm (d)} $(F_4,\so(9,\C))$
\mypar  {\rm (e)} $(E_6,F_4)$
\mypar 
\noindent Each pair $(\g,\k)$ in this list is non-Hermitian.}

From the point of view of representation theory, the condition
that $\Omin\cap\gR$ is non-empty is very natural. To explain this, we
recall the theory of the associated variety.

Suppose    $\pi_o:\GR\to\Unit\H$ is an irreducible unitary
representation.  Let $H\subset\H$  be  the space of $\KR$-finite 
vectors with its natural $(\g,K)$-module structure; $H$ is then the 
Harish-Chandra  module of the representation.  Differentiation  of the
group representation gives a Lie algebra representation of $\gR$ on 
$H$ and so a representation  $\wt{\pi}:\Ug\to\End H$ of the universal
enveloping algebra. The annihilator $\I$ is then the primitive  ideal
attached to
$\pi_o$. The graded ideal $\gr\I$ cuts out a closed complex algebraic 
subvariety  ${\cal V}(\gr\I)\subset\g^*\simeq\g$ called the associated
variety of  $\I$.   Since
$\pi_o$ admits a central character, it   follows  that
${\cal V}(\gr\I)$ is  a union of complex nilpotent orbits. 
A basic result (due independently to   Borho and J.L. Brylinski, 
to Ginzburg, and to Joseph)  is that ${\cal V}(\gr\I)$ is in fact the 
closure  of a {\it single} nilpotent orbit, which is then called  the
associated complex  nilpotent orbit of  $\H$ and $H$.
 
The following observation is an easy consequence of the theory of
associated varieties. For instance, it is a  corollary of 
[Vo2, Theorem 8.4].
\lem{4.2}{Suppose  $\GR$ admits an irreducible unitary
representation with associated complex nilpotent orbit $O\subset\g$.
Then $O\cap\gR$ is non-empty, i.e., $O$ has a real form with
respect to $\gR$.}

We will call an irreducible  unitary representation
$\pi_o:\GR\to\Unit\H$   {\it  minimal}  if its associated nilpotent orbit 
is $\Omin$ and also the   image of $\wt{\pi}:\Ug\to\End H$ has no
zero-divisors , i.e., the annihilator $\I$ of $\wt{\pi}$ is completely
prime.  For $\g$ not of type $A_n$, $\pi_o$ is minimal  if and only if
$\I$ is the Joseph ideal. So Lemma 4.2 says that a necessary (but not
sufficient) geometric  requirement  for $\GR$ to admit a minimal
representation is that $\Omin\cap\gR$ is non-empty.

From now on in this paper, we assume that $\OR$ is strongly minimal 
and   $(\g,\k)$ is non-hermitian.  For convenience, we also take
$\GR$ to be simply-connected. There is no problem in this as 
the universal cover has finite center.

We freely identify $\OR$ with $Y$  via the Vergne diffeomorphism (2.3).
Using Lemma 3.4 and (3.14), we find the dimension of $Y$ is given by:
$$\matrix{\g& \sl(n,\C)&\so(n,\C)&\sp(2n,\C)&
G_2&F_4&E_6&E_7&E_8\cr
\dimC~Y&n-1&n-3&n&3&8&11&17&29} \eqno(4.1)$$

Next we want to develop the Jordan theory interpretation of the 
polynomial $P$ defined in (3.22). We   find in Proposition 4.4 below a
Jordan structure on the  space $\k_{-1}$ from (3.10). 

There is a natural  symmetry group acting on the space $\k_{-1}$,
namely the isotropy group $K_0=K^h$ for the adjoint action of $K$. 
We use this symmetry throughout the paper.
 $K_0$ is a closed reductive complex algebraic subgroup
of $K$ with Lie algebra $\k_0$. Also $K_0$ is connected; this follows
immediately from the fact that the adjoint orbit $K\cdot h$ is
simply-connected.  Basic
constructions like  (3.10), (3.16), and (3.17) break the $K$-symmetry 
but not  the $K_0$-symmetry. 

In particular, $K_0$ acts on 
$\p_2$ by a (non-trivial) character 
$$\chi:K_0\to\C^*\eqno(4.2)$$ 
so that $a\cdot e=\chi(a)e$ for $a\in K_0$.   Let 
$$K_0'=\hbox{kernel of   }\chi=\Ks\eqno(4.3)$$
and let $\k_{0}'=\ks$  be the Lie algebra of   $K_{0}'$.
We get an  orthogonal decomposition 
$$\k_0=\k_{0}'\oplus\C h\eqno(4.4)$$
 
Since $\k_{1}$ is abelian, we have a natural identification
$$\U(\k_1)=S(\k_1)\eqno(4.5)$$
So in particular, $P$ defines an element of $\U(\k_1)$.
We recall  from  [B-K4,\S2.2-2.6]:
\prop{4.3}{ The polynomial $P\in S^4(\k_1)$ defined in {\rm(3.22)} is
semi-invariant under $K_0$ and   transforms  by the character
$\chi^2$. Moreover, $P$ is, up to scaling, the unique 
$K_0$-semi-invariant polynomial in $S(\k_1)$  such that
$$\ad_P(\C\oe)=\C e\eqno(4.6)$$}

We recall some work from  [B-K4, \S2.6-8]. We found with Kostant 
a nilpotent  element $\ek\in\k_1$  such that  
$$\l=\C h\oplus\C\ek\oplus\C\oek\eqno(4.7)$$
is a complex Lie subalgebra in $\k$ isomorphic to $\sl(2,\C)$
and $(2h,-\ek,\oek)$ is an S-triple basis of  $\l$.
We normalized the choice of $\ek$ so that
${1\over 4!}(\ad\,{\ek})^4(\oe)=e$ and 
${1\over 4!}(\ad\,{\oek})^4(e)=\oe$. Hence 
$$P(\oek)=1\eqno(4.8)$$
We note that the nilpotents $e\in\p$ and $\ek\in\k$ were called,
respectively, $z$ and $e$ in   [B-K4].

Then we showed  that $(\k,\k_0)$ is a Hermitian symmetric pair
of tube type with rank $q$ where $q\le 4$.
In addition, the pair $(\k_0,\kl)$  
is  a complex symmetric pair so that we have a  complex Cartan
decomposition $\k_0=\kl\oplus\r$ where $[\r,\r]\subset\kl$.

Consequently, elaborating on [B-K4, Proposition 2.8] we get
\prop{4.4}{The Tits-Kantor-Koecher construction   gives  $\k_{-1}$  
the structure of a  complex semisimple Jordan algebra   with
$\Kl$-invariant Jordan product defined by 
$$[x,\oek]\circ[y,\oek]=[x,[y,\oek]]\eqno(4.9)$$ where 
$x,y\in\r$. The Jordan  identity element is $\oek$. The Jordan
algebra degree of  $\k_{-1}$ is 
$$\deg \k_{-1}  = q=\hbox{rank } (\k,\k_0)\le 4\eqno(4.10)$$ }

The T-K-K theory identifies $\k_{-1}$ as the complexification
$\J_{\C}$ of a  real Euclidean   Jordan algebra $\J$.
The book [F-K] is an excellent  reference for  the theory of  
real Euclidean and complex semisimple Jordan algebras.

Next we write out the decomposition of $\k_{-1}$ into a direct sum of 
complex simple Jordan subalgebras:
$$\k_{-1}=\j_{[1]}\oplus\cdots\oplus\j_{[\ell]}\eqno(4.11)$$
Then each space $\j_{[n]}$ carries an irreducible  representation
of  $\k_0$.   Let $q_n$ be the degree of  $\j_{[n]}$; then
$$q_1+\cdots +q_{\ell}=q\eqno(4.12)$$
Let $P_{[n]}$ be the Jordan norm of  $\j_{[n]}$; then $P_{[n]}$ has    
degree   $q_n$.

\prop{4.5}{The polynomial $P\in S^4(\k_1)$, constructed in {\rm(3.22)}, 
 as a function on $\k_{-1}$,  is uniquely expressible as a monomial 
$$P=P_{[1]}^{w_1}\cdots P_{[\ell]}^{w_\ell}\eqno(4.13)$$
in the Jordan norms $P_{[n]}$ of the simple components $\j_{[n]}$
of   $\k_{-1}$. Every exponent $w_1,\dots,w_{\ell}$ is positive.}
\proof{In the  Tits-Kantor-Koecher construction,
the simple  Lie components of $\k$ correspond to the simple
Jordan components of   $\k_{-1}$. I.e.,
$\j_{[n]}=\k_{-1}\cap\r_{[n]}$  where 
$\k=\r_{[1]}\oplus\cdots\oplus\r_{[\ell]}$ is the decomposition of $\k$  
into complex simple Lie subalgebras. Then we get also the decomposition 
$\k_0=\t_{[1]}\oplus\cdots\oplus\t_{[\ell]}$ into a direct sum of
subalgebras where $\t_{[n]}=\k_{0}\cap\r_{[n]}$.   
Each pair $(\r_{[n]},\t_{[n]})$  is complex symmetric of tube type.
Then any  polynomial 
function on $\j_{[n]}$  semi-invariant under $\t_{[n]}$ is a polynomial in
$P_{[n]}$; see, e.g. [K-S, Th. 0]. It follows that $P$ is of the form
$P=cP_{[1]}^{w_1}\cdots P_{[\ell]}^{w_\ell}$  for some scalar $c$.
But $c=1$ since $P(\oek)=P_{[1]}(\oek)=\cdots=P_{[\ell]}(\oek)=1$.
Finally the fact that $h$ has non-zero projection to each component 
$\r_{[1]},\dots,\r_{[\ell]}$ implies that each $w_1,\dots,w_{\ell}$
is non-zero.  \qed}

Notice that (4.13) {\it defines}  the exponents $w_1,\dots,w_{\ell}$ and
gives the numerical equality
$$q_1w_1+\cdots q_{\ell}w_{\ell}=4\eqno(4.14)$$

In [B-K4], we quantized with Kostant the real form $\OR$ of $\Omin$
in the  three cases where $\cal J_{\C}$  is a simple complex
Jordan algebra of degree $4$ so that $\ell=1$ and $P=P_{[1]}$.  In this
paper we  treat the general  case where $P$ may factor non-trivially.

We proceed in the rest of  this section to make an explicit
list of the Jordan algebras occurring here. 
The associations we get   between exceptional 
Lie algebras and Jordan algebras are in many cases already familiar
from the constructions  discovered by Tits, Kantor, Koecher,  and  
Allison  and Falkner  to produce    exceptional  Lie algebras out of   
Jordan algebras.

In our  tables  we adopt the following conventions. We write
$\so_p$, $\sp_{2p}$, $\sl_p$, $G_2$, $F_4$, $E_6$, 
$E_7$, $E_8$ for the corresponding  complex Lie algebras.
Also $S_{o}^n\C^p$ denotes the irreducible representation 
of $\so(p,\C)$ satisfying $S_{o}^n\C^p+S^{n-2}\C^p\simeq
S^{n}\C^p$ while $\wedge_{o}^n\C^{2p}$ denotes the
irreducible representation of  $\sp(2p,\C)$ satisfying 
$\wedge_{o}^n\C^{2p}+\wedge^{n-2}\C^{2p}
\simeq\wedge^{n}\C^{2p}$ where $n\le p$.

A complete list of all
non-isomorphic formally real simple Jordan
algebras of degree $\le 4$ follows immediately from the
Tits-Kantor-Koecher theory and the known list of all
irreducible Hermitian symmetric tube domains
(see [H], p. 528, Example 4 and \S6.4, pp. 518-520).
We give this list in  Table 4.6 together with the 
corresponding pair $(\k,\k_0)$. 
The number $d$ arises in the following way.
The restricted root system for the pair
$(\k,\k_0)$ is of type $C_q$ where $q$ is the degree of $\J$
and then the  long roots have multiplicity 1 while the 
short roots all have common multiplicity $d$.

In Table 4.6,  $\R^{p}$, $p\ge 1$, is the $p$-dimensional real Jordan 
algebra associated to the  Euclidean norm and $\Herm(n,\Bbb F)$ is
the real Jordan algebra of $n\times n$ hermitian matrices over 
$\Bbb F$ where $\Bbb H$ and $\Bbb O$ denote the quaternions and
the octonions (the Cayley numbers) respectively. Then 
$\Herm(3,\Bbb O)$ is the exceptional $27$-dimensional Jordan
algebra while all the others in Table 4.6 are special (i.e., arise in the
standard way from  associative algebras).
The last column in Table 4.6 gives a name to  the Jordan norm of
$\J$.
\vskip 2pc
\centerline{\bf Table 4.6.  All Simple Euclidean  Real Jordan
Algebras 
$\J$ of rank $\le 4$}
\vskip .6pc
{\settabs\+aaaaaaaaaaaaaaaaaaaaaa&ddddddddddd&bbbbb&bbbbb&
bbbbbbbbbbb&ccccccc&ddddddd\cr
\+  \phantom{XXX}$\J$ & {\hskip -5pt}$\dimR\J$ & $d$ &  $\k$ &
 $\k_0$ & {\hskip -6pt}$\deg \J$ &{\hskip -6pt}Norm\cr
\vskip -6pt
\+{\hbox to 29pc{\hrulefill }}\cr
\+ $\J(1)=\R$ & $1$ & $0$ &$\sl_2$ & $\so_2$ & $1$ & $P_1$ \cr
\+ $\J(2;p)=\R^{p-2}$,$\,\,p\ge 5$ & $p-2$ &{\hskip -5pt}$p-4$ 
& $\so_{p}$ &$\so_{p-2}\oplus\so_2$ & $2$ & $P_{2;p}$ \cr 
\+ $\J(3,\Bbb R) = \Herm(3,\Bbb R)$ & $6=3+3d$ & $1$
& $\sp_{6}$ & $\gl_3$ & $3$ & $P_{3;\R}$ &\cr 
\+ $\J(3,\Bbb C) = \Herm(3,\Bbb C)$ & $9=3+3d$
&$2$ &$\sl_{6}$ & $\s(\gl_3\oplus\gl_3)$ & $3$ & $P_{3;\C}$ &\cr 
\+ $\J(3,\Bbb H) =
\Herm(3,\Bbb H)$ & $15=3+3d$ &$4$ &$\so_{12}$ & $\gl_6$ 
& $3$ & $P_{3;\Bbb H}$ &\cr 
\+ $\J(3,\Bbb O) = \Herm(3,\Bbb O)$ & $27=3+3d$ &$8$ &$E_7$
 & $E_6\oplus\so_2$ & $3$ & $P_{3;\Bbb O}$ &\cr 
\+ $\J(4,\Bbb R) = \Herm(4,\Bbb R)$ & $10=4+6d$ &$1$ 
&$\sp_{8}$ &
$\gl_4$ & $4$ & $P_{4;\R}$ &\cr \+ 
$\J(4,\Bbb C) = \Herm(4,\Bbb C)$ & $16=4+6d$ &$2$
&$\sl_{8}$ & $\s(\gl_4\oplus\gl_4)$ & $4$ & $P_{4;\C}$ &\cr 
\+ $\J(4,\Bbb H) = \Herm(4,\Bbb
H)$ & $28=4+6d$ &$4$ &$\so_{16}$ & $\gl_8$ & $4$ & 
$P_{4;\Bbb H}$ &\cr} 
\vskip 2pc

In Table 4.7 we list all pairs $(\g,\k)$ occurring here 
(i.e.,  non-hermitian  complex symmetric pairs $(\g,\k)$
with $\Omin\cap\gR\neq\emptyset$ where $\g$ is simple)
together with the Jordan algebra $\J$  arising from the
T-K-K theory and the polynomial function $P\in S^4(\k_1)$ 
on $\cal J_\C$  written  as the product of the Jordan norms.

\vskip 2pc
\centerline{\bf Table 4.7    Non-Hermitian 
 pairs $(\g,\k)$ with $\Omin\cap\gR\neq\emptyset$ }
\vskip .6pc 
{\settabs\+aaaaaaaaaaaaaaaaaaaaaaaaaaaaa&qqqqqqq
&bbbbbbbbbbb
&kkkkkkkkkk&ppppppppppp&ggggg\cr
\+  \phantom{XXX}$\J$ & $q$ & $P$ & $\k$ & $\p$ & $\g$ \cr
\vskip -6pt
\+{\hbox to 32pc{\hrulefill }}\cr
\+$\J(4,\Bbb R)$  & $4$&  $P_{4,\Bbb R}$ & $\sp_8$  & 
$\wedge_{o}^4\C^8$ & $E_6$ \cr
\+$\J(4,\Bbb C)$  & $4$& $P_{4,\Bbb C}$  & $\sl_8$  & 
$\wedge^4\C^8$ & $E_7$ \cr
\+$\J(4,\Bbb H)$   & $4$&  $P_{4,\Bbb H}$ & $\so_{16}$  
& $\C^{128}$ & $E_8$  \cr
\+$\J(3,\Bbb R)\oplus\J(1)$  & $4$& $P_{3,\Bbb R}P'_1$ & 
$\sp_6\oplus\sl_2$  & $\wedge_{o}^3\C^6\otimes\C^2$ & $F_4$  \cr
\+$\J(3,\Bbb C)\oplus\J(1)$  & $4$& $P_{3,\Bbb C}P'_1$ & 
$\sl_6\oplus\sl_2$  & $\wedge^3\C^6\otimes\C^2$ & $E_6$  \cr
\+$\J(3,\Bbb H)\oplus\J(1)$  &  $4$&$P_{3,\Bbb H}P'_1$  & 
$\so_{12}\oplus\sl_2$  & $\C^{32}\otimes\C^2$ & $E_7$ \cr
\+$\J(3,\Bbb O)\oplus\J(1)$  &$4$&  $P_{3,\Bbb O}P'_1$  & 
$E_7\oplus\sl_2$  & $\C^{56}\otimes\C^2$ & $E_8$ \cr
\+$\J(1)\oplus\J'(1)$  & $2$& $P_{1}^3{P}'_1$  & $\sl_2\oplus\sl_2$ 
 & $S^3\C^2\otimes \C^2$ & $G_2$   \cr
\+$\J(2;p)\oplus\J'(2;q)$ & $4$& $P_{2;p}P'_{2;q}$  & 
$\so_{p}\oplus\so_q$  & $\C^p\otimes\C^q$ & $\so_{p+q}$   \cr
\+$\J(2;p)\oplus\J(1)\oplus\J'(1)$  &  $4$&$P_{2;p}P_1P'_1$ & 
$\so_{p}\oplus\so_4$  & $\C^p\otimes\C^4$ & $\so_{p+4}$   \cr
\+$\J(1)\oplus\J'(1)\oplus\J''(1)\oplus\J'''(1)$ 
&$4$&  $P_1P'_1P''_1P'''_1$ 
& $\sl_2^{\oplus 4}$  & $(\C^2)^{\otimes4}$ & $\so_8$ \cr
\+$\J(2;p)\oplus\J(1)$  & $3$& $P_{2;p}P_{1}^2$ & $\so_{p}\oplus\so_3$  
& $\C^p\otimes\C^3$ & $\so_{p+3}$  \cr
\+$\J(1)\oplus\J'(1)\oplus\J''(1)$  &$3$& $P_{1}^2P'_1P''_1$ & 
$\so_3\oplus\so_4$ & $\C^3\otimes\C^4$ & $\so_7$  \cr
\+$\J(2;p)$     &   $2$&   $P_{2;p}^2$  & $\so_{p}$  & 
$S_{o}^2\C^p$ & $\sl_p$       \cr
\+$\J(1)\oplus\J'(1)$  &  $2$&$P_{1}^2{P'_{1}}^2$ & $\so_3\oplus\so_3$ 
& $\C^3\otimes\C^3$ & $\so_6$ \cr
\+$\J(1)$  & $1$& $P_{1}^4$ & $\sl_2$  & $S^4\C^2$ & $\sl_3$  \cr
\+  $p,q\ge 5$ throughout the table &&&&\cr}
\vskip 2pc

Comparing Tables 4.6 and 4.7, we find 
\prop{4.8}{There is a bijection between 
{\rm (i)} the pairs  $(\g,\k)$ in Table {\rm 4.7} and 
{\rm (ii)} the triples $(\J, P)$ where  $\J$
is a  Euclidean real Jordan algebra and $P$  is a  monomial in the  
Jordan norms $P_{[n]}$ of the simple components of $\J$ such 
that each $P_{[n]}$ occurs at least once in $P$ and 
$P$ has total degree {\rm 4}.} 

Notice that the condition that $\J$ has degree $\le 4$ is 
necessary but not sufficient for $\J$ to occur here (indeed
$\J$ cannot be a Jordan algebra of rank 3) and that the same
$\J$ can give rise to different polynomials $P$ and hence
different $\g$ (this occurs for $\g=G_2$ and $\g=\so_6$).
 
From Table 4.7, we get  in Table 4.9 a list of  the real semisimple
groups  $\GR$ occurring here.

\vskip 2pc
\centerline{\bf Table 4.9.}\nobreak
\vskip 1pc
{\settabs\+GGGGGGG&kkkkkkkkkkk&
llllllllllll&dddddd&mmmmmmmm&tttttttttttttttt\cr
\+$\GR$&$\k$&rank&$d$& $m$ \cr \vskip -6pt
\+{\hbox to 25pc{\hrulefill }}\cr
\vskip 10pt
\+$E_{6(6)}$ & $\sp_8$ & $6$ &$1$ &$4+6d=10$ \cr
\+$E_{7(7)}$ & $\sl_8$ & $7$ & $2$ &$4+6d=16$\cr
\+$E_{8(8)}$ & $\so_{16}$ &$8$ & $4$&$ 4+6d=28$\cr
\+$F_{4(4)}$ &$\sp_6\oplus\sl_2$ & $4$ &$1$&$4+3d=7$\cr
\+$E_{6(2)}$ &$\sl_6\oplus\sl_2$ &$4$
&$2$&$4+3d=10$\cr
\+$E_{7(-5)}$ &$\so_{12}\oplus\sl_2$ & 
$4$&$4$& $4+3d=16$\cr
\+$E_{8(-24)}$ &$E_7\oplus\sl_2$ & $4$&$8$ &$4+3d=28$\cr
\+$G_{2(2)} $ & $\sl_2\oplus\sl_2$ & 
$2$ & ${2\over 3}$ &$(4+3d)/3=2$&\cr
\+$\wt{SO}(p,q)$ & $\so_p\oplus\so_q$& 
$p$  & &$p+q-4$ & $3\le p\le q$\cr
\+$\wt{SL}(n,\R)$ &$\so_n$& $n-1$  &  &$n-2$& $3\le n$  \cr}
\vskip 2pc

In the listing of the exceptional real groups, the subscripted  number
in parentheses is equal to $\dimR\pR-\dimR\kR$ and serves
to distinguish between  simply-connected real forms having the same
complexified Lie algebra.   In  the first three cases
$d$ is the ``correct" parameter for the corresponding simple
Jordan algebra $\k_{-1}$ while in the next five cases $d$ is a
fictitious  parameter which we make up as it gives
consistent formulas in Tables 4.9 and  6.9.

\Sec{\S5.  Holomorphic Half-Form Bundles on   $\OR$.}

In this section, we construct and classify all holomorphic 
half-form bundles $\NH$ over $\OR$ equipped with its instanton
Kaehler structure $\JJ$ from \S2.  Right away, we identify
$(\OR,\JJ)$ with the complex cone $Y$ in $\p$ by means of the Vergne
diffeomorphism (2.3). 

We are assuming, throughout the rest of the  
paper, that   $\OR$ is strongly minimal and   
$(\g,\k)$ is non-hermitian.  Then $\k$ is a semisimple Lie algebra and 
$\p$ is irreducible as a  representation of $K$. The spaces $\OR$ and $Y$
are given by (3.7).  In particular
$Y$ is the conical $K$-orbit of highest weight vectors in $\p$. 
The cases occurring here were classified in Table 4.7.

Each holomorphic half-form bundle $\NH$ over $Y$
is  automatically homogeneous under $K$ (see Lemma 5.2). The space 
$H=\Gamma(Y,\NH)$ of global algebraic holomorphic sections  breaks up
under the action of
$K$ into a multiplicity free ladder decomposition which we analyze in 
Lemmas 5.2 and 5.3. 

In Proposition 5.5 we determine the spectrum of the operator   $E'$     
on $H$ given by the  Lie   derivative of the holomorphic Euler vector  field
$E$.  We find  that $E'$ is diagonalizable with    {\it positive} spectrum.
This result about $E'$ is
crucial since it allows us to {\it invert} $E'$ and $E'+1$ in \S6  in order to
quantize the symbol $g_{v_0}$ from Theorem 3.6.

We regard $E'$ as the ``energy" operator on the space 
$H$ of quantization. Indeed, $E'=\Q(\la)=\Q(\rho)$ by Corollary 2.2(i).
I.e.,  $E'$ is the quantization of our chosen Hamiltonian $\rho$
(see discussion before (2.8)).  Since $\rho$ is positive everywhere 
on $\OR$, the positivity of $E'$ is exactly what we expect from the
quantum  theory.

The Euler $\C^*$-action on $Y$ (see \S2) defines  a complex 
algebra grading  
$$R(Y)=\bigoplus_{p\in\ZP}R_p(Y)\eqno(5.1)$$
where
$$R_p(Y)=\{f\in R(Y)\,|\, Ef=pf\}\eqno(5.2)$$ 
Here $\ZP$ denotes the set of non-negative integers.
A priori, the grading in (5.1) extends over all integers, but 
since $Y$ is the orbit of highest weight vectors in $\p$, it follows that
the pullback map $S(\p^*)\to R(Y)$ on functions is surjective.
So $R_p(Y)=0$ for $p$ negative.
 
Our first aim is to compute the fundamental group of $Y$. To do this, we 
use fact that $Y$ is a homogeneous space of $K$:  
$$Y\simeq{K/K^e}\eqno(5.3)$$
Since $\pi_1(K)=0$, it follows that $\pi_1(Y)$ is isomorphic to the
component group of $K^e$.  

Now  let  $Q\subset{K}$ be  the (closed) subgroup which 
preserves the line $\C e$. Clearly $K^e$ lies in $Q$ as the kernel  of
the action of $Q$ on $\C e$; we put $Q'= K^e$.
The quotient ${\Bbb P}(Y)$ of $Y$ by the 
Euler $\C^*$-action  identifies with $K/Q$. Since $Y$ is an orbit of 
highest weight vectors,  ${\Bbb P}(Y)$ is a (generalized) flag variety of
$K$.  So we have a $K$-equivariant  principal $\C^*$-bundle
$$Y\,\to\, {\Bbb P}(Y)\,\simeq\,K/Q\eqno(5.4)$$

It  follows  easily  that $Q$ is the connected subgroup of $K$ with Lie
algebra $\q=\k_0\oplus\k_1$. Thus we get the Levi decomposition
$$Q=K_0\ltimes K_1\eqno(5.5)$$
where   $K_1=\exp\k_1$ is the connected unipotent subgroup of $K$ with 
Lie  algebra $\k_1$.

\lem{5.1}{$Y$ is simply-connected except if  $\g$ is of type  $A_n$.
If $\g=\sl(n,\C)$,  then the fundamental group of $Y$ is 
$\pi_1(Y)\simeq \Z_2$ if $n\ge 4$ or
$\pi_1(Y)\simeq \Z_4$ if $n=3$.}
\proof{The discussion above gives $K^e=K_0'\ltimes K_1$
where $K_0'$ was defined in (4.3).
The exponential map  $\exp:\k_1\to K_1$ is an
isomorphism.   Hence, the component groups of $K^e$ and $K_0'$ 
identify. So $\pi_1(Y)$ is isomorphic to the component
group of $K_0'=\Ker\,\chi$ where  $\chi$ is the $K_0$-weight of $e$.
\mypar
Proposition 4.3 says that  $\chi$  is the square root
of the weight of  $P$.   (The square root is unique since $K_0$ is
connected.) Thus the product decomposition (4.13) gives
$\chi^2=\chi_{1}^{2w_1}\cdots\chi_{\ell}^{2w_\ell}$ 
where $\chi_{n}^2$ is the weight of  the Jordan norm $P_{[n]}$. So
$$\chi=\chi_{1}^{w_1}\cdots\chi_{\ell}^{w_\ell}\eqno(5.6)$$
\mypar
The well-known theory of the $K_0$-action on $S(\k_1)$  
(see [B-K3,  Theorem 3.4 and Corollary 3.6])  says that
$\chi_{1},\dots,\chi_{\ell}$  are primitive characters of $K_0$  and
form a basis of the character group. The kernel of a primitive character
is connected.  It follows that the component group of  $\Ker\,\chi$ is
isomorphic to $\Z/g\Z$ where $g=\gcd\{{w_1,\dots,w_{\ell}}\}$.
Now the Lemma follows immediately from Table 4.7 where $P$ is 
given explicitly in the third column as a monomial in
$P_{[1]},\dots,P_{[\ell]}$ so that the exponents $w_1,\dots,w_\ell$ can be
read off  in each case. \qed}
 
It  is often convenient to label irreducible representations of $K$ by
their highest weight in the sense of the Cartan-Weyl highest weight
theory. This requires that we fix a choice of Cartan subalgebra $\h$ in
$\k$  together with a Borel subalgebra $\b$ in $\k$. The
set $\cal R^+$ of non-zero weights of $\h$ on $\b$ is the set of                
{\it positive roots}.  Then $\cal R=\cal R^+\cup -\cal R^+$ is the set of 
all non-zero weights of $\h$ on $\g$. Each $\alpha\in\cal R$ is called
a       {\it root} and  the corresponding weight space $\k^{\alpha}$ is
$1$-dimensional and is called the $\alpha$-root space.
Then 
$$\k=\h\oplus\m\oplus\m^-\eqno(5.7)$$ where
$\m\subset\b$ is the span of the  positive root spaces and $\m^-$ is 
the span of the negative root spaces.

Cartan-Weyl highest weight theory says that in each
finite-dimensional irreducible $\k$-representation $V$, the subspace
$V^\m$ of  all vectors annihilated by the action of $\m$ is
$1$-dimensional. Clearly then $V^\m$ is a weight space of $\h$ of
some weight $\mu$. Moreover $V^\m$ turns out to be the full
$\mu$-weight space $V^{\mu}$ in $V$. Then $\mu$ is the so-called 
{\it highest weight} in $V$,  any (non-zero) vector in $V^\mu$ is
called a  {\it highest weight vector}, and we write $V=V_{\mu}$.
  
We   choose $(\h,\b)$ so that  $h\in\h$ and  $e$ is a highest
weight vector in $\p$.  We also require that $\k_1$ lies in $\m$,
$\h$ is complex conjugation stable
and  $\h$ is  stable under  the complex Cartan
involution of the symmetric pair $(\k_0,\kl)$.  It is easy to
meet these conditions. Then  
$$\p\simeq V_{\psi}\eqno(5.8)$$
where $\psi$ is the weight of $e$ so that $\psi$ is  the highest 
weight  of $\p$. Notice that $\psi$ is just the restriction to $\h$ 
of  the $\k_0$-weight $\psi=\d\chi\in\k_{0}^*$.
 
From now on in this section, we work in the category of
algebraic holomorphic complex line bundles $L$ over a smooth
(irreducible) complex algebraic manifold $X$, which will soon be
specialized to $X=Y$. $R(X)$ and  $\Gamma(X,L)$ denote  the regular (i.e.,
algebraic holomorphic) functions on $X$ and sections of $L$.

Next  we discuss some general  notions about half-form bundles.
Let $L\to X$ be  a complex line bundle. A {\it square root} of  $L$  is
a  pair $(C,\alpha)$ where $C$ is a complex line bundle over $X$ and 
$\alpha:C^{\otimes 2}\to L$ is a bundle isomorphism.  
Notice that $\alpha$ gives an $R(X)$-module map
$$\alpha_X:\Gamma(X,C)^{\otimes 2}\to\Gamma(X,L)$$ 
so that the product of two sections of $C$ defines a section of $L$.

Two square roots $(C,\alpha)$ and $(C',\alpha')$ are isomorphic if 
there  exists a bundle isomorphism $\beta:C\to C'$ such that 
$\alpha=\alpha'\circ\beta^{\otimes 2}$ where 
$\beta^{\otimes 2}:C^{\otimes 2}\to (C')^{\otimes 2}$ is the obvious 
map. In counting or classifying half-form bundles, we will always 
work up to isomorphism.

It is easy to check that if  $s\in\Gamma(X,L)$ is a non-zero section then
there exists, up to isomorphism, at most one square root $(C,\alpha)$
of $L$ such that $s$ is the square of some section   of $C$.  In practice, 
we  suppress the isomorphism $\alpha$ from the notation.  
Notice that any line bundle is a square root of its square.

We will use the notation and terminology from [B-K3] on 
algebraic holomorphic differential 
operators and their symbols. If $\eta\in\D^1(X,L)$, i.e., $\eta$ is an 
order $1$ differential operator on sections of $L$, then  $\eta$
determines an  order $1$ differential operator on 
sections of any square root $C$ of
$L$  in the following way: there exists a unique operator
$\eta^\sharp\in\D^1(X,C)$  such that $\eta(s^2)=2s\eta^\sharp(s)$ for all
$s\in\Gamma(X,C)$.  Then the  symbols  of $\eta$ and $\eta^\sharp$
coincide.  We write $\eta$ for $\eta^\sharp$ when the meaning is clear.

The next result follows easily  using general facts about homogeneous
line bundles for the first part and  
[B-K3, Lem. 2.9 and Appendix \S{A}.12] for the second part. 
In particular, the proof of the second part uses the Borel-Weil theorem
on ${\Bbb P}(Y)$. 

 Recall that $f_0\in R(Y)$ was defined in (3.15). Using the
terminology of this section we see that $f_0$ is a highest weight
vector in $R_1(Y)$ of weight $\psi$ and $f_0$ is $Q$-semi-invariant.
 
\lem{5.2}{Suppose $C$ is a square root of a $K$-homogeneous 
line bundle on $Y$. Then $C$ has {\rm(}uniquely{\rm)} the structure of a
$K$-homogeneous line bundle.  The space $H=\Gamma(Y,C)$ of 
global sections is non-zero. The differential of the corresponding 
$K$-action on $H$   is  a representation of $\k$ on $H$ by  
differential operators of order   $1$ and is compatible with tensor
product of bundles. \mypar
The $\k$-representation on $H$ is completely reducible and
multiplicity-free.  A vector in $H$ is a highest weight vector for the
$\k$-action if and only if it is $Q$-semi-invariant.
The space $H^{\m}$  of highest weight vectors has a basis of the form
$\{\fo^ps_0\,|\, p\in\ZP\}$  where $s_0\in H$ is uniquely 
determined  up to scaling.  This gives a ladder  decomposition
$$H\simeq\bigoplus_{p\in\ZP} V_{\nu+p\psi}\eqno(5.9)$$ of $H$ as a
$\k$-representation where $\nu$ is the
$\h$-weight of $s_0$. }
 
We write $H_{[\nu+p\psi]}$ for the subspace of $H$ which carries the
$\k$-representation  $V_{\nu+p\psi}$. The representation $V_\nu$ is
the minimal $\k$-type in $H$. We may  also call 
its carrier space $H_{[\nu]}$ the minimal $\k$-type in $H$.

A   square root $\NH$ of  the canonical line bundle
$\NN=\wedge^n(T^*X)$ on $X$, where $n=\dimC X$,  
is called an algebraic holomorphic 
{\it half-form} line bundle on $X$.  Sections of
$\NH$ are called {\it half-forms}. Let $\L_{\xi}$ denote the Lie derivative
of   a  vector field  $\xi$ on $X$ so that in particular we get an
operator $\L_{\xi}\in\D^1(X,\NN)$.   The  Lie derivative on half-forms is
the  operator   $\L_{\xi}\in\D^1(X,\NH)$ given by
$$\L_{\xi}(s^2)=2s\L_{\xi}(s)\eqno(5.10)$$ 
In this way we get a representation
$$\Vect(X)\to\D(X,\NH)\to\End H,\hskip 2pc
\xi\mapsto\L_{\xi}\eqno(5.11)$$
where $\Vect(X)$ is the Lie algebra of algebraic holomorphic vector
fields on $X$. Then on the symbol level
$$\symbol \L_{\xi} \,= \,\symbol \xi\,\in R_{[1]}(T^*X)\eqno(5.12)$$
where $R_{[p]}(T^*X)$ denotes the space of regular functions on
the cotangent bundle $T^*X$ which are homogeneous of degree $p$
on the fibers of the natural projection $T^*X\to X$.

Now the Lie derivative of the Euler vector field $E$ on $Y$ is the 
operator
$$E'=\LE\eqno(5.13)$$   
The following result is easy to verify and defines
the minimal $E'$-eigenvalue $r_0$.

\lem{5.3}{Suppose $\NH$ is  a half-form bundle on $Y$. Then
the $K$-action  on $\NH$ gives a natural representation of $K$ on 
the space $H=\Gamma(Y,\NH)$ of global sections. Differentiating  this
gives the representation 
$$\pi'_K:\k\to\D(Y,\NH)\to\End H,\hskip 2pc
x\mapsto\L_{\eta^x}\eqno(5.14)$$
\mypar
$E'$ is diagonalizable on $H$ with spectrum $r_0+\ZP$ where
$r_0\in\half\Z$.  This defines $r_0$.
Thus we have the eigenspace decomposition
$$H=\bigoplus_{p\in\ZP}H_{r_0+p}\eqno(5.15)$$ 
where $H_q$ is  the  $q$-eigenspace of $E'$ in H. The action of $E'$ 
on $H$ commutes with the  $K$-action and  the eigenspaces   
$H_{r_0+p}$ are
the irreducible $K$-submodules in $H$. Furthermore $H_{r_0+p}$
carries the representation $V_{\nu+p\psi}$ so that
$$H_{r_0+p}=H_{[\nu+p\psi]}\eqno(5.16)$$   
Let $s_0\in H_{r_0}$ be a non-zero $Q$-semi-invariant section.
Then  $\{\fo^ps_0\in H_{r_0+p}\,|\, p\in\ZP\}$ is a   complete set
of linearly independent $Q$-semi-invariant sections in $H$.} 

The Hilbert space of  our quantization of $\OR$ will be a  certain
completion of  $H=\Gamma(Y,\NH)$.
We think of $E'$ as the  energy operator  and call
$H_{r_0}$ the  {\it vacuum space} in $H$. The  sections in
$H_{r_0}$ are the {\it vacuum vectors} in $H$.
The   vector $s_0$ chosen in $H_{r_0}$
is unique up to scaling.

A {\it  spherical} $K$-representation is one that
contains a non-zero $K$-invariant vector, which is then  called the
spherical vector.  It follows from Lemma 5.3 that $H$ is  
$K$-spherical if  and only if  $H_{r_0}\simeq\C$, so if and only if
$s_0$ is $K$-invariant.

\rem{5.4}{Our discussion of square root bundles and half-form
bundles  generalizes in the  obvious way to $n$th  root bundles and
$n$th roots of the canonical bundle.  
Then all the results in Lemmas 5.2 and  5.3 go over to the case
where we replace $\NH$ by an $n$th root of the canonical
bundle, the only change being that then $r_0\in{1\over n}\Z$ instead
of   $r_0\in\half\Z$.} 

Our key result on half-form bundles is
\prop{5.5}{Suppose $\NH$ is  a half-form bundle on $Y$. Then
the minimal eigenvalue $r_0$ of $E'$ on $H=\Gamma(Y,\NH)$
is positive. Thus 
$E'$ has positive spectrum $r_0+\ZP$ on $H$.}

For the proof, we need to construct a concrete holomorphic 
$(m+1)$-form $\Lambda$ on   $Y$. We will use this form $\Lambda$ 
throughout \S5 and \S7..  We construct $\Lambda$ out of
the functions $\fo,\fa\dots,\fm$ on $Y$ we defined in  (3.15) 
in the following way:
$$\Lambda=\d \fo\wedge\d\fa\wedge\cdots\wedge
\d\fm\in\Gamma(Y,\NN)\eqno(5.17)$$
Since $\fo,\fa\dots,\fm$ were coordinates  on the open set $\Yo$ 
defined  in (3.17), it follows that  $\Lambda$  is nowhere vanishing on
$\Yo$.  Furthermore $\Lambda$ is $Q$-semi-invariant and 
$$\ze =  \hbox{ $Q$-weight  of $\Lambda$}\eqno(5.18)$$
is the character by which $Q$ acts on $\p_2\otimes\wedge^m\p_1$.
\prooff{of Proposition 5.5}{We will show that  the spectrum of
$\L_{E}$ on $\Gamma(Y,\NN)$ is positive. This  implies the 
positivity of the spectrum of $E'$ on $H$.
\mypar  
It follows by  Lemma 5.2 that $\Gamma(Y,\NN)$ is
a multiplicity-free $K$-module of ladder type and the set
$\{\fo^{p-\be}\Lambda\,|\, p\in\ZP\}$ is a complete set of linearly
independent  highest weight vectors, where $\be$ is some  
non-negative integer.  Let
$$\Lambda_0=\fo^{-\be}\Lambda\eqno(5.19)$$ 
Since  $E(f_v)=f_v$ for all $v\in\p$
and $\L_{\eta^h}(f_v)=kf_v$ if $v\in\p_k$,  we find
$$\eqalign{\L_{E}(\Lambda_0)=&t\Lambda_0\quad\hbox{with}\quad
t=-\be+1+m\cr
\L_{\eta^h}(\Lambda_0)=&j\Lambda_0 \quad\hbox{with}\quad
j=-2\be+2+m \cr}\eqno(5.20)$$  
But  $j\ge 0$. This follows from the representation theory of the Lie
algebra  $\l$ defined in  (4.7)  as $\L_{\eta^{x}}(\Lambda_0)=0$ for
all $x\in\k_1$ and so in particular for $x=\ek$. Therefore 
$\Lambda_0$ is the highest weight vector  of a finite-dimensional 
irreducible  $\l$-representation.  
So $2t=j+m$ is positive  and hence  $t$ is positive.
This gives the result as  $t+\ZP$ is the spectrum of $E'$ on 
$\Gamma(Y,\NN)$.  \qed}

\prop{5.6}{\item{\rm (i)} If $\g$ is not of type $A_n$, then up to 
isomorphism $Y$ admits at most  one half-form bundle. 
\item{\rm (ii)} If   $\g=\sl(n,\C)$ {\rm(}$n\ge 3${\rm)}, then either 
$Y$ admits no  half-form bundle or $Y$ admits exactly  two 
non-isomorphic half-form bundles. 
\item{\rm (iii)} If  $\NH$ is a half-form bundle on $Y$ with
$H=\Gamma(Y,\NH)$ and  $s_0\in H_{r_0}$   is a non-zero
$Q$-semi-invariant  vector,  then up to scaling  
$$\hbox{either }\qquad s_{0}^2=\Lambda_0\qquad\hbox{or}\qquad
 s_{0}^2=f_0\Lambda_0 \eqno(5.21)$$
Moreover the two possibilities in {\rm(5.21)} classify half-form
bundles on $Y$ up to isomorphism.\mypar}
\proof{From formal properties of square root bundles, it follows that
if the canonical bundle $\NN$ on $Y$ admits any square root, then the 
set of all square roots  up to isomorphism is parameterized 
by the order $2$ characters of $\pi_1(Y)$. Thus by Lemma 5.1, if  $Y$
admits a half-form bundle then it is unique  with the exception of    the
cases $\g=\sl(n,\C)$  where there would be  two half-form bundles.
This proves (i) and (ii).
\mypar Now suppose we are  given $\NH$ and $s_0$.
Recall $\Lambda_0$ from (5.19).
The ladder structure on $\Gamma(Y,\NN)$  (Lemma 5.2) implies that
(up to scaling) $s_{0}^2=\fo^p\Lambda_0$ where $p\in\ZP$.
If $p\ge 2$ then we consider the local section 
$\fo\i s_0\in\Gamma(\Yo,\NH)$. The square of $\fo\i s_0$ is equal to 
$\fo^{p-2}\Lambda_0$ and this has no poles on $Y$. It follows that 
$\fo\i s_0$ has no poles on 
$Y$ and so $\fo\i s_0\in\Gamma(Y,\NH)$. But $\fo\i s_0$ has 
$E'$-degree equal to $r_0-1$ and so
this contradicts the minimality of $r_0$. Thus  $p=0$ or $p=1$.
Furthermore, $\Lambda_0$ and   $\fo\Lambda_0$
cannot both be squares of sections of the same bundle $\NH$
since, by (5.1), $\fo\in R_1(Y)$ is not a square in $R(Y)$.
\qed}

\prop{5.7}{\item{\rm (i)} If $\GR=\wt{SO}(p,q)$ where $p+q$ is odd 
and $p,q\ge 4$, then $Y$  admits no half-form bundle.
\item{\rm (ii)} If $\GR=\wt{SL}(p,\R)$ where $p$ is odd and $5\le p$,
then $Y$  admits no half-form bundle.
\item{\rm (iii)} If $\GR=\wt{SL}(n,\R)$ where $n\ge 4$ is even,
then $Y$ admits exactly two half-form  bundles $\NH_+$ and $\NH_-$.
These may be characterized by the conditions: $\Lambda_0$ is the
square of a section of $\NH_+$ and $\fo\Lambda_0$ is the square of a
section of $\NH_-$.
\item{\rm (iv)} In all other cases,  $Y$ admits a unique half-form 
bundle $\NH$ and $\Lambda_0$ is the square of a section of $\NH$.
\mypar
\noindent All the half-form bundles occurring here are listed in Table 
{\rm 6.9} along with  their minimal $E'$-eigenvalue $r_0$ and the
$\k$-type $V_{\nu}$ of  the vacuum space $H_{r_0}$.}
\proof{ By (5.19) the form $\Lambda_0$ is 
$K_0$-semi-invariant of weight $\zeta_0=\chi^{-\be}\zeta$.   It
follows  easily that $Y$ admits a half-form bundle with
$s_{0}^2=\Lambda_0$ (respectively  $s_{0}^2=\fo\Lambda_0$)
if and only if  $\zeta_0$ (respectively $\chi\zeta_0$)
is the square of a character  $\tau$ of $K_0$. Then
$\tau$ is   unique and $\d\tau=\nu$ is the highest weight
of $H_{r_0}$.
\mypar
Thus we need to compute the  character $\zeta_0$ and look for 
square roots in the character group of  $K_0$. 
The Jordan structure on $\k_{-1}$ gives us a convenient way 
to  do the calculations.   Indeed we found   in the proof of
Lemma 5.1 (i) the basis $\chi_1,\dots,\chi_{\ell}$  of the character group
of $K_0$  and (ii) the formula (5.6).
\mypar
Now  by (5.17) $\Lambda$ transforms in the $1$-dimensional 
$K_0$-representation 
$$\p_2\otimes\wedge^{m}\p_1~\simeq~
\chi^{m+1}\otimes \wedge^{m}\k_{-1}\eqno(5.22)$$  But 
$\wedge^{m}\k_{-1}$ is the tensor product of the  top exterior 
powers of the spaces $\j_{[1]},\dots,\j_{[\ell]}$. The 
weight of $K_0$ on the top exterior power of $\j_{[n]}$  
is  $\chi^{-u_n}_n$ where  $u_n=2+d_n(q_n-1)$.  
Here $d_n$ is the  root multiplicity parameter
of $\j_{[n]}$  given in Table 4.6 and the formula
for $u_n$ is immediate from the description of the 
restricted root system of $(\k_0,\kl)$. Thus the weight of $\Lambda$ is 
by (5.6)
$$\zeta=\chi^{m+1}\chi_{1}^{-u_1}\cdots\chi_{\ell}^{-u_\ell}
=\chi_{1}^{(m+1)w_1-u_1}\cdots\chi_{\ell}^{(m+1)w_\ell-u_\ell}
\eqno(5.23)$$
\mypar
We   can decide if a section $s=\fo^p\Lambda$ over $\Yo$  extends to $Y$ 
just by examining its  weight $\chi_1^{t_1}\cdots\chi_\ell^{t_\ell}$.
Indeed $s$ extends to $Y$ if and only if $t_1,\dots,t_\ell\ge 0$.
This  follows easily by  using the Borel-Weil theorem  on ${\Bbb P}(Y)$
(in its geometric form involving  the orders of poles along irreducible 
divisors in the complement of the big cell)   
to compute $\Gamma(Y,\NN)$. Thus, going back to the definition of
$\Lambda_0=\fo^{-\be}\Lambda$ in (5.19), we see that $\be$
is  the largest non-negative integer such that all the numbers
$-\be w_n+(m+1)w_n-u_n$ are non-negative. To simplify, we put
 $\al=m+1-\be$. Then the weight  of  $\Lambda_0$ is
$$\ze_0=\chi_{1}^{\al w_1-u_1}\cdots\chi_{\ell}^{\al w_\ell-u_\ell}
\eqno(5.24)$$
where  $\al$ is the smallest positive integer such that $\al w_n\ge u_n$  
for all $n$.
 \mypar
It is now easy to go through Table 4.7,  calculate $\zeta_0$
in each case, and see if  $\zeta_0$ and/or $\chi\zeta_0$ admits
a square root.  The cases group together naturally into families.
In the  first three cases  in Table 4.7 we have
$P=P_{4,\Bbb R},P_{4,\Bbb C},P_{4,\Bbb H}$ so that
$\ell=1$, $w_1=1$, $u_1=2+3d$. Then $\al=2+3d$ and
$\zeta_0=1$.  Thus we
get a half-form bundle  (unique as $\g\neq\sl(p,\C)$)
with $s_{0}^2=\Lambda_0$ and  $H_{r_0}\simeq\C$. 
\mypar
In  the four cases where $P=P_{3,\Bbb F}P'_1$, we have
$\ell=2$, $(w_1,w_2)=(1,1)$ and $(u_1,u_2)=(2+2d,2)$ where $d=d_1$.
Then $\al=2+2d$ and $\zeta_0=\chi_{2}^{2d}$ which is a   square.
So we get a unique half-form bundle with $s_{0}^2=\Lambda_0$ and 
$H_{r_0}\simeq\C\otimes S^d\C^2$. The $G_2$
case $P=P_{1}^3P'_{1}$ is similar with
$(w_1,w_2)=(3,1)$ and $(u_1,u_2)=(1,1)$. Then $\al =2$ and 
$\zeta_0=\chi_{1}^4$ which is a square.  This gives a 
unique half-form bundle with $s_{0}^2=\Lambda_0$ and we find
$H_{r_0}\simeq S^2\C^2\otimes\C$ since $\chi=\chi_{1}^2\chi_2$
is the highest weight of $\p=S^3\C^2\otimes\C^2$.
\mypar
We can treat all the  cases where $\gR=\so(p,q)$, $3\le p\le q$,
simultaneously  with $\ell=2$ as long as we formally set
$\J(2;4)=\J(1)\oplus\J(1)$, $P_{2;4}=P_1P'_1$,
and $\J(2;3)=\J(1)$, $P_{2;3}=P_{1}^2$ with $S^{k/2}\C^3=S^k\C^2$
as $\so(3)$-representations. This follows easily and 
we   get  $(w_1,w_2)=(1,1)$ and  $(u_1,u_2)=(p-2,q-2)$.
Then $\al =q-2$ and  $\zeta_0=\chi_{1}^{q-p}$.
If $q-p$ is even, then $\zeta_0$ is a square, and we
get a half-form bundle with $s_{0}^2=\Lambda_0$ and 
$H_{r_0}\simeq S^{(q-p)/2}\C^2\otimes\C$ since
$\chi=\chi_1\chi_2$ is the highest weight of 
$\p\simeq\C^p\otimes\C^q$. This is the unique half-form bundle 
unless  $(p,q)=(3,3)$ in which case
we get a second half-form bundle with  $s_{0}^2=\fo\Lambda_0$
and $H_{r_0}^{\triangle 2}\simeq\p$ 
(where $\triangle$ denotes Cartan product)   so that 
$H_{r_0}\simeq\C^2\otimes\C^2$.
Indeed $\so(3,3)=\sl(4,\R)$ and this is the only time when  our
$\so(p,q)$ cases and  $\sl(n,\R)$ cases coincide. 
Now if  $q-p$ is odd with  $p>3$, 
then neither $\zeta_0$ nor $\chi\zeta_0=(\chi_{1}^{q-p+1},\chi_2)$ 
are squares.  However, if    $q-p$ is odd with  $p=3$, then $\zeta_0$,
but not $\chi\zeta_0$, is a square.  Thus we get one half-form bundle
with $s_{0}^2=\Lambda_0$ and 
$H_{r_0}\simeq S^{(q-3)/2}(\C^3)\otimes\C=S^{q-3}\C^2\otimes\C$.
\mypar
Finally we consider the $\gR=\sl(p,\R)$ cases. If $p\ge 5$, then
$P=P^2_{2;p}$,  $\ell=1$, $w_1=2$, $u_1=p-2$. So if $p$ is even  
then $\al=(p-2)/2$, $\zeta_0$ is trivial and we  get a half-form bundle
with $s_{0}^2=\Lambda_0$ and  $H_{r_0}\simeq\C$. The second
half-form bundle has  $s_{0}^2=\fo\Lambda_0$ and
$H_{r_0}^{\triangle 2}\simeq\p$ so that $H_{r_0}\simeq\C^p$.
If $p$ is odd, then there is no half-form bundle.
We already did the case $p=4$. If $p=3$ then
$P=P_{1}^4$, $\ell=1$, $w_1=4$, $u_1=2$  and so $\al=1$
and $\zeta_0=\chi_{1}^2$. 
We  get one  half-form bundle with
$s_{0}^2=\Lambda_0$ and $H_{r_0}\simeq\C^2$, and a second
with  $s_{0}^2=\fo\Lambda_0$ and   $H_{r_0}\simeq S^3\C^2$.
\mypar
We have now proven everything except for the values of $r_0$.
But (5.19)   and (5.20) 
give $E'\Lambda_0=(-\be+m+1)\Lambda_0=\al\Lambda_0$.
So if  $s_{0}^2=\Lambda_0$ then $r_0=\al/2$ while
if  $s_{0}^2=\fo\Lambda_0$ then $r_0=(\al+1)/2$.
We have computed the parameter $\al$ in each case above, and this
produces the values of $r_0$ in  Table 6.9. 
\qed}

\rem{5.8}{Case (i)  is the {\it Howe-Vogan counterexample}.   Howe
proved  that these groups $SO(p,q)$ admit  no minimal unitary
representation  and then  Vogan ([Vo1]) extended this to the
simply-connected covering groups.}

The isomorphism (3.16) implies in particular that $R(\Yo)$ is the
localization of $R(Y)$ at $\fo$ so that
$$R(\Yo)=R(Y)[\fo\i]=\C[\fo,\fa,\dots,\fm][\fo\i]\eqno(5.25)$$
It is easy to  prove 
\prop{5.9}{Suppose $C$ is a $K$-homogeneous line bundle
on $Y$ and $s\in\Gamma(Y,C)$ is a non-zero $Q$-semi-invariant vector. 
Then $s$ is nowhere vanishing on
$\Yo$.  Consequently, since $\Yo$ is affine,  the space of sections
$\Gamma(\Yo,C)$ is a cyclic $R(\Yo)$-module generated by  $s$ so 
that $$\Gamma(\Yo,C)=\Gamma(Y,C)[\fo\i]=R(\Yo)s\eqno(5.26)$$
\mypar In particular then, in Proposition {\rm 5.6},
$H=\Gamma(\Yo,\NH)$ is  a cyclic $R(\Yo)$-module generated by any
section  $\fo^ps_0$.}
This is a key result as it enables us to analyze $H$ in a uniform  
manner in \S7 below regardless of whether $H$ is $K$-spherical or
not.  In fact, we  further simplify our work in \S7 by making the 
following observation, which obviates the need to consider
separately the two possibilities in (5.21).

The regular function $\fo\in R_1(Y)$ is not a square in $R(\Yo)$
(because of (5.1) and (5.25))  and is nowhere vanishing on $\Yo$ by 
the definition of $\Yo$. Thus we may  construct a non-trivial two-fold
\'etale covering  
$$\wt{Y^o}\to Y^o\eqno(5.27)$$ by ``extracting  a square root of 
$\fo$".  Then  $\wt{Y^o}$, like $\Yo$,   is again a smooth 
affine complex algebraic variety and has a unique $K_0$-action
such that the cover $\wt{Y^o}\to Y^o$ is $K_0$-equivariant.
The pull back of $\NH$ through the covering is  a 
$K_0$-homogeneous bundle on  $\wt{Y^o}$ which we again
call $\NH$. Now Proposition 5.9 gives

\cor{5.10}{The space
$\Gamma(\wt{\Yo},\NH)$ of   algebraic holomorphic sections  is a
cyclic 
$R(\wt{\Yo})$-module generated by $\sL$ so that
$$\Gamma(\wt{\Yo},\NH)=
\C[\fo^{1\over 2},\fo^{-{1\over 2}},\fa,\dots,\fm]\sL\eqno(5.28)$$}

\Sec{\S6.  Quantization of  $\OR$.}

In this section we construct explicitly our quantization of $\OR$.
This is purely ``from scratch"; we assume no a priori information on the 
existence of any quantizations or unitary representations.
In the next section we prove that the constructions of this section 
``work", i.e., we prove Theorem 6.3 and 6.8.

The first step of our quantization, carried out in \S2 and \S3,  was to
transform the quantization problem on $\OR$ into a quantization
problem on $T^*Y$. We replaced (by holomorphic extension) each
function $\phi^w$, $w\in\gR$, by a rational  (pseudo-differential)
symbol $\Phi^w\in R(T^*Y)[\la\i]$. Then, after complexification, 
we ended up
in (3.6) with a realization of $\g$ as a complex  Lie algebra of rational
symbols,
$$\Phi^z\in R(T^*Y)[\la\i],\qquad  z\in\g\eqno(6.1)$$

Our aim now is to quantize the symbols $\Phi^z$, $z\in\g$,  into
operators $\Q(\Phi^z)$ on a Hilbert space $\H$ which
is a completion of $H=\Gamma(Y,\NH)$ for some half-form bundle 
$\NH$ over $Y$. As usual, we freely identify $(\OR,\JJ)$ with $Y$ via the
Vergne diffeomorphism (2.3). 
We require our operators satisfy certain explicit and implicit axioms.
This solves our quantization problem on
$\OR$ as we set $\Q(\phi^w)=\Q(\Phi^w)$, $w\in\gR$. 

In \S5 we have already  quantized the symbol $\la$,   corresponding to   
our chosen Hamiltonian $\rho$ (cf. Corollary 2.2(i)), into the operator 
$E'$ on half-forms.

We require that the operators $\Q(\Phi^w)$, $w\in\gR$, be 
self-adjoint, or equivalently, that the operators $\Q(\Phi^z)$ satisfy
$$\Q(\Phi^z)^{\dagger}=\Q(\Phi^\oz),\qquad  z\in\g\eqno(6.2)$$ 
Of course at this point, $H$ carries no preferred positive definite
Hermitian inner product. So we will have to construct this along the 
way. Our operators will not be defined everywhere on $\H$, but they
will all contain $H$ in their domain. 

We further require that the Dirac commutation relations (2.4) be
satisfied. If we put
$$\pi^z=i\Q(\Phi^z),\qquad  z\in\g\eqno(6.3)$$ then
these relations amount to  
$$[\pi^z,\pi^{z'}]=\pi^{[z,z']}\eqno(6.4)$$ 
for all $z,z'\in\g$. I.e., the map 
$\pi:\g\to\End H$, $z\mapsto\pi^z$,  must be a complex Lie
algebra homomorphism.

In order to satisfy the implicit axioms, we require that the 
``symbol"  of  $\Q(\Phi^z)$ is just $\Phi^z$.  Here our definition of 
``symbol" is not precise, as we found in Theorem 3.1(iii) that the symbols
$\Phi^v$, $v\in\p$ are not homogeneous. However, we will get around 
this by dealing individually with the homogeneous pieces.

With this in mind, we mandate  
$$\pi^x=i\Q(\Phi^x)=\L_{\eta^x},\qquad x\in\k\eqno(6.5)$$
So  $\pi^x$ is the Lie derivative on half-forms of the algebraic
holomorphic vector field $\eta^x$ defined in (2.6) by differentiating 
the $K$-action on $Y$.  Thus, just as we would expect, $\pi^x$ 
corresponds to the natural
$K$-action on $H$, i.e., $\pi^x=\pi'_K(x)$ in the notation of  (5.14).

Next, guided by the complex Cartan decomposition (3.5),  we need to
quantize the symbols $\Phi^v$, $v\in\p$. In (3.4) we found these
break into a sum $f_v+g_v$ of two homogeneous pieces. 
We now  mandate  
$$\pi^v=i\Q(\Phi^v)=if_v+iT_v,\qquad v\in\p\eqno(6.6)$$
where $T_v=\Q(g_v)$ is some ``nice" quantization of the 
homogeneous degree $2$ rational symbol  $g_v$ from (3.4). 

This leaves the  problem of how to construct  $T_v$.
Of course we want in the end for $\pi$ to be a complex Lie algebra
homomorphism. So we certainly want 
$[\pi^x,\pi^v]=\pi^{[x,v]}$ and this  implies   
$$[\pi^x,T_v]=T_{[x,v]}\eqno(6.7)$$
Hence the fact that $\p$ is irreducible as a $\k$-representation  
insures that $T_{v_0}$ and the  $\pi^x$,
$x\in\k$, already determine all operators $T_v$, $v\in\p$.

To construct $T_{v_0}$, we break down the symbol $g_e=g_{v_0}$
computed in Theorem 3.6 by (3.25) into its elementary factors.
On the face of it, it seems hard to imagine  what to do with
the factor $f_0\i(\Phi_KP)$ appearing in (3.4). While $\Phi_KP$ is the
symbol of  the perfectly nice order $4$
differential operator $\pi'_K(P)$ on 
sections of  $\NH$  because of  (4.5) and (5.14),  the quotient
$f_0\i(\Phi_KP)$ a priori only defines a differential operator on 
sections of the bundle $\NH$ restricted to the open set $Y^o$ from
(3.17).  Fortunately,  the result in  [B-K3, Ths. 3.10 and  4.5] (which
applies more generally to any homogeneous line bundle over $Y$)
tells us

\theo{6.1}{Suppose $\NH$ is a half-form bundle on $Y$ and 
$H=\Gamma(Y,\NH)$. Then the operators  $\fo$ and $\pi'_K P$ on $H$
commute and have the same image. Hence the formula
$$D_e={1\over \fo}(\pi'_K P)\eqno(6.8)$$ 
defines an operator on $H$. It follows that $D_e$ is
an algebraic  differential operator of order $4$ on sections of $\NH$.
\mypar
The assignment $e\mapsto D_e$ extends naturally and uniquely to
a complex linear $1$-to-$1$ $K$-equivariant map
$$\p\to\D(Y,\NH),\qquad  v\mapsto D_v\eqno(6.9)$$   
of $\p$ into the 
algebra $\D(Y,\NH)$ of algebraic differential operators on 
sections of $\NH$. 
Then, for each non-zero $v\in\p$, $D_v$ has  order $4$; also  
$D_v$ has  degree $-1$, i.e., $[E',D_v]=-D_v$.
\mypar
The subalgebra $\A\subset\D(Y,\NH)$ generated by 
the $D_v$, $v\in\p$,
is  abelian, isomorphic to $R(Y)$ and graded by 
$\A=\oplus_{p\ge 0}\A_{-p}$ where 
$\A_{-p}=\{D\in\A\,|\,[E',D]=-pD\}$.
Putting $D_v=D_{f_v}$ for  $v\in\p$, we get a graded 
$K$-equivariant complex  algebra isomorphism  
$$R(Y)\to\A,\qquad  f\mapsto D_f\eqno(6.10)$$
\mypar
There is a unique $\KR$-invariant  positive-definite  Hermitian inner
product $B_o$ on $H$  such that  $B_o(s_0,s_0)=1$ 
\(when we fix a choice of $s_0$ in Lemma {\rm 5.3}\)  and 
the operators $f_v$ and $D_{\ov}$ are adjoint with respect to 
$B_o$ for all $v\in\p$. Then the  grading {\rm(5.15)} is a  
$B_o$-orthogonal decomposition.}

The expression for $D_e$ in terms of our local coordinates (3.15) 
on $Y$  is, in the notation of  (3.24),
$$D_e=\fo^3
P\big(\L_{\pdb{\fa}},\dots,\L_{\pdb{\fm}}\big)\eqno(6.11)$$

To complete our quantization of  the rational symbol $g_e$, we  need
to quantize the factor $-\la^{-2}$ in (3.25). It is natural to try to
quantize
$-\la^{-2}$ into the operator $(E'+a)\i(E'+b)\i$ where $a$ and $b$ are
some constants to be determined.  Of course  $a$ and $b$ must be
chosen so that neither $-a$ nor $-b$ belongs to the spectrum of $E'$. 

In fact it turns out that 
exactly one choice, namely $a=0$ and  $b=1$, satisfies our 
requirement that the resulting operators defined by (6.5) and (6.6) 
satisfy the bracket relations of $\g$. In fact,   just the one relation
$[\pi^e,\pi^{\oe}]=\pi^{h}$ mandates that
$$-{1\over\la^2}\hbox{ quantizes to }  {1\over E'(E'+1)}\eqno(6.12)$$
(We prove in \S7 that this choice  works.)
We emphasize that the operators $E'$ and $E'+1$ are invertible since 
the spectrum of $E'$ on $H$ is positive by Proposition 5.5.

Thus, rather than putting 
$T_v={1\over  (E'+a)(E'+b)}D_v$ and solving for $a$ and $b$ in \S7
we simply define
$$T_v=\Q(g_v)={1\over  E'(E'+1)}D_v\eqno(6.13)$$
In [B2, proof of Theorem 4.2],  we   wrote out   the
latter procedure of determining $a$ and $b$ from the bracket relation 
for the  case $\gR=\sl(3,\R)$.

These operators $T_v$ are no longer differential operators, 
but they share many of the same properties. To explain this, we 
introduce the notion of $\k$-finite endomorphism.

We have a natural representation of $\k$ on $\End H$ defined by
$x* D=[\L_{\eta^x},D]$.  Then $D\in\End H$ is called {\it $\k$-finite}
if $D$ generates a finite-dimensional representation of $\k$ inside 
$\End H$. 
The space $\End_{\k-fin} H$ of all $\k$-finite endomorphisms of $H$ is 
a complex subalgebra of $\End H$. Let  
$$\End_{[p]}H\subset\End_{\k-fin}H\eqno(6.14)$$  
be the $p$-eigenspace of  $\ad\,E'$.   We have a  natural graded
$K$-equivariant complex algebra inclusion (see [B-K4, A.6, A.12]) 
$\D(Y,\NH)\subset \End_{\k-fin}H$.

Now Theorem 6.1 easily gives (cf. proof of [B-K4, Th. 3.4])
\cor{6.2}{The operator $T_v$,  $v\in\p$, lies $\End_{[-1]}H$.
Thus we get a $K$-equivariant complex linear  map
$$T:\p\to\End_{[-1]}H,\qquad v\mapsto T_v\eqno(6.15)$$
The pseudo-differential operators $T_v$, $v\in\p$,  
commute and generate a graded abelian  $K$-stable subalgebra  
$\T=\oplus_{p\ge 0}\T_{-p}$
of  $\End_{\k-fin}H$  where $\T_{-p}=\T\cap\End_{[-p]}H$.
We then get a graded $K$-equivariant complex  algebra 
isomorphism 
$$R(Y)\to\T, \qquad f\mapsto T_f\eqno(6.16)$$
where $T_{f_v}=T_v$ for $v\in\p$.
There is a $\KR$-invariant  positive-definite 
Hermitian inner product  $B$ on $H$, such that  $B(s_0,s_0)=1$ 
\(when we fix a choice of $s_0$ in Lemma {\rm 5.3}\)  and 
the operators $f_v$ and $T_{\ov}$ are adjoint with respect to 
$B$ for all $v\in\p$.   Then the  grading {\rm(5.15)} is a  
$B$-orthogonal decomposition.}

Taking inventory of our operators, we see that
$\pi^x$, $x\in\k$, and $f_v$, $T_v$, $v\in\p$ are each graded
operators on $H$ of degrees $0$,$1$, and $-1$ respectively. I.e.,
we have  
$$\pi^x:H_t\to H_t,\qquad f_v:H_t\to H_{t+1},\qquad
T_v:H_t\to H_{t-1}\eqno(6.17)$$

Now we can state
\theo{6.3}{Suppose $\NH$ is a half-form bundle on $Y$ and 
$H=\Gamma(Y,\NH)$.  Let  
$$\pi:\g\to\End_{\k-fin}H,\qquad z\mapsto\pi^z\eqno(6.18)$$
be the complex linear  map defined by
{\rm(6.5)},  {\rm(6.6)}, {\rm(6.13)} so that
$$\eqalign{\pi^x&=i\Q(\Phi^x)=\,\L_{\eta^x}
\hskip 49pt\hbox{if  $x\in\k$},\cr  
\pi^v&=i\Q(\Phi^v)=\, if_v+iT_v \hskip 26pt\hbox{if
$v\in\p$}}\eqno(6.19)$$    Then, except in the one  case where
$\GR=\wt{SL}(3,\R)$, $r_0=1$, and
$H_1\simeq\C^4$, the map  $\pi$ is a complex Lie algebra 
homomorphism 
so that $\pi$ is a representation of $\g$ by global algebraic 
pseudo-differential operators on sections of $\NH$.} 
\proof{As in [B-K4,\S6],  the problem reduces  to proving the  single
bracket relation of operators on $H$
$$[\pi^e,\pi^{\oe}]=\pi^h\eqno(6.20)$$ 
because of   [B-K4, Lem. 3.6].   Since the operators satisfy 
$$[f_v,f_{v'}]=[T_v,T_{v'}]=0\eqno(6.21)$$
 for all $v,v'\in\p$, we
get $[\pi^e,\pi^{\oe}]=[f_\oe,T_{e}]-[f_e,T\soe]$. Thus (6.20)  
amounts  to the relation
$$[f_\oe,T_{e}]-[f_e,T\soe]=\L_{\eta^h}\eqno(6.22)$$
We   prove (6.22) in \S7. We also show that in the  $SL(3,\R)$ case
we omitted, $\pi$  fails to  be a Lie algebra homomorphism. \qed} 

For the rest  of this section, we assume that we are in the situation of 
Theorem {\rm 6.3} with the one bad case excluded so that
$\pi$ is a Lie algebra homomorphism. 
Let $$\wt{\pi}:\U(\g)\to\End_{\k-fin} H\eqno(6.23)$$
be the complex algebra homomorphism  defined by $\pi$. Let
$\E$ be the image of $\wt{\pi}$ and let $\I\subset\U(\g)$ be the 
kernel of $\wt{\pi}$.
Then we have  a natural complex algebra isomorphism 
$$\U(\g)/\I\simeq\E\eqno(6.24)$$
Let  $S^{[p]}(\g)$ be the $p$th Cartan power
of the adjoint representation of $\g$.

\theo{6.4}{The representation $\pi$ of $\g$ on $H$ is irreducible. 
The algebra homomorphism $\wt{\pi}$ is surjective so that
$$\E=\End_{\k-fin} H\eqno(6.25)$$
The algebra $\E$  has no zero-divisors. Thus 
the annihilator $\I$ of $\pi$ is a completely prime primitive ideal
in $\U(\g)$.
\mypar
We have    $\E^p/\E^{p-1}\simeq S^{[p]}(\g)$ as  $\g$-modules
and so there is a   multiplicity free  $\g$-module decomposition
$$\E\simeq\oplus_{p\in\ZP}S^{[p]}(\g)\eqno(6.26)$$
The associated graded ideal $\gr\I\subset S(\g)$
is the prime ideal defining the closure of ${\Omin}$.
Thus  the associated graded map to $\wt{\pi}$  gives  a graded 
algebra isomorphism  $$R(\Omin)\iso \gr\E\eqno(6.27)$$}
\proof{The proofs in [B-K4,\S5] of the corresponding results  
go through verbatim (by design) into this more general setting. The 
only change needed is that we replace 
the line ``Let $s_0=1\in H_{r_0}$" in [B-K4] by 
``Let $s_0\in H_{r_0}$ be a non-zero highest weight vector for the 
$\k$-action." \qed}

If $\g\neq\sl(n,\C)$, then Joseph ([J2]) proved  that  $\U(\g)$ 
contains a unique
completely prime primitive ideal with associated nilpotent orbit
$\Omin$ (cf. \S4). This is called the Joseph ideal. Thus we get
\cor{6.5}{If $\g\neq\sl(n,\C)$ then $\I$ is the Joseph ideal.}
We remark that this gives a new proof of Garfinkle's ([G]) result 
that the associated graded ideal of the Joseph ideal is prime.

We now fix  a non-zero vacuum vector $s_0\in H_{r_0}$ 
which is  $Q$-semi-invariant, or equivalently, a  highest  weight vector   
for the $\k$-action.
 
\theo{6.6}{$H$ admits a  unique $\gR$-invariant positive-definite
Hermitian   inner product   $\<s|s'\>$
such that $\<s_0|s_0\>=1$. This    coincides with
the inner product $B$ found in Corollary {\rm  6.2} so that
$$B(s,s')=\<s|s'\>\eqno(6.28)$$   
\mypar  Consequently, the representation of $\gR$ on $H$ given by 
$\pi$  integrates uniquely to give a  unitary representation
$$\pi_{o}:{\GR}\to\Unit\,\H\eqno(6.29)$$
on the Hilbert space  $\H$ obtained by completing $H$ with respect to
$B$. Then $H$ is the space of $\KR$-finite vectors in $\H$.}
\proof{This follows by the proof of 
[B-K4, Th. 5.2], using  the same modification described in
Theorem 6.4.   See, e.g., [W, \S6.A.4] for a proof of  Harish-Chandra's
theorem that the $\gR$-action on an admissible finitely generated 
$(\g,K)$-module $S$ endowed with a
$\gR$-invariant positive-definite Hermitian inner product  
integrates to a unitary representation of ${\GR}$ on the
Hilbert space completion of $S$.  We apply this with
$H=S$. Indeed $H$ is irreducible by Theorem 6.4 and so    generated
by any non-zero vector.  Also  $H$  is admissible (i.e., all
$K$-multiplicities are finite) since $H$ is  in  fact multiplicity free  
by   Lemma 5.3.\qed}

We will write   $\<s|s\>=\ns{s}$.
Theorems 6.4 and 6.6  give, in the language of \S4,   the representation
theoretic result  
\cor{6.7}{$\pi_o$ is a minimal unitary representation of 
${\GR}$ and $H$ is its associated 
Harish-Chandra $(\g,K)$-module.}

\theo{6.8}{ There exist  positive real numbers $a$ and $b$ 
{\rm(}depending  on  $\GR$ and $\NH${\rm)} such that 
$$\biggns{f_{0}^ns_0\over n!}=~{(a)_n(b)_n\over n!(r_0+1)_n}
\eqno(6.30)$$
for all $n\in\ZP$. Moreover $a$ and $b$ are unique up to ordering 
and satisfy
$$a+b=r_0+1+X_0\eqno(6.31)$$
where $X_0$ is the eigenvalue of $\L_{\eta^h}$ on $s_0$.
We compute $a$ and $b$ below in Table {\rm 6.9}. 
\mypar
The values $\ns{f_{0}^ns_0}$ and $\KR$-invariance uniquely 
determine the  inner product $B$ on $H$ because of the ladder
decomposition {\rm(5.15)}. }
Here we are using the hypergeometric function notation
$(a)_n=a(a+1)\cdots(a+n-1)$.
\proof{The adjoint of multiplication by $f_0=f_e$ is $T_\oe$
by Corollary 6.2. We find
$T\soe(f_{0}^ks_0)=\ga_kf_{0}^{k-1}s_0$ for $k\in\ZP$ where
$\ga_k$ is a scalar and $\ga_0=0$. This follows by $E'$-degree and 
weight as  in  [B-K4, proof of Th. 5.2]. Indeed, 
$T\soe(f_{0}^ks_0)$ lies in 
$H_{r_0+k-1}$ and has $\h$-weight $\nu+(k-1)\psi$.
But  by Lemma 5.3,  $f_{0}^{k-1}s_0$ is a highest weight vector
of weight  $\nu+(k-1)\psi$ in the irreducible $\k$-representation
$H_{r_0+k-1}$. Thus $T\soe(f_{0}^ks_0)$ is a multiple of 
$f_{0}^{k-1}s_0$. We now find
$$\ns{f_{0}^ns_0}=\<s_0\,|\, T\soe^n(f_{0}^ns_0)\>=
\ga_1\cdots\ga_n\<s_0\,|\, s_0\>\eqno(6.32)$$
We evaluate the RHS of  (6.32) in \S7 below and obtain (6.30) and (6.31).
The final statement that these values determine $B$  is immediate 
from Lemma 5.3 -- in particular the ladder decomposition of $H$ is
multiplicity-free. \qed}

\centerline{\bf Table 6.9.}\nobreak
\vskip 1pc
{\settabs\+ccccccc&kkkkkkkkkkk&
sssssssss&VVVVVVVVVV&rrrrrrrrrrrrrrrr&aaaaaaaaaaaaa&
bbbbbbbbbbb&\cr
\+Case&$\GR$& $s_{0}^2$ & $V_{\nu}\simeq H_{r_0}$
& $\quad r_0$ & $\quad a$ & $\quad b$\cr \vskip -6pt
\+{\hbox to 31pc{\hrulefill }}\cr
\vskip 10pt
\+(i)& $E_{6(6)} $  &   $\Lambda_{0}$ 
& $\C$  & $1+{3\over 2}d={5\over 2}$ &
$1+{1\over 2}d={3\over 2}$  &$1+d={2}$ \cr
\+(ii)& $E_{7(7)}$  &   $\Lambda_{0}$
& $\C$ &$1+{3\over 2}d=4$ & 
$1+{1\over 2}d=2$  &$1+d=3$\cr
\+(iii)& $E_{8(8)}$   & $\Lambda_{0}$& $\C$ 
& $1+{3\over 2}d=7$ &$1+{1\over 2}d=3$  &$1+d=5$\cr
\+(iv)& $F_{4(4)}$   &  $\Lambda_{0}$&
$\C\otimes S^1\C^2$ & $1+d=2$ & $1+{1\over 2}d={3\over 2}$ 
&$1+d={2}$\cr
\+(v)& $E_{6(2)}$  & $\Lambda_{0}$& 
$\C\otimes S^2\C^2$ & $1+d=3$&$1+{1\over 2}d=2$  &$1+d=3$\cr
\+(vi)& $E_{7(-5)}$   & $\Lambda_{0}$& 
$\C\otimes S^4\C^2$& $1+d=5$&$1+{1\over 2}d=3$  &$1+d=5$\cr
\+(vii)& $E_{8(-24)}$  & $\Lambda_{0}$& 
$\C\otimes S^8\C^2$
& $1+d=9$ &$1+{1\over 2}d=5$  &$1+d=9$\cr
\+(viii)& $G_{2(2)}$   & $\Lambda_{0}$ 
&$S^2\C^2\otimes\C$ & $1$ & $1+{1\over 2}d={4\over 3}$  
&$1+d={5\over 3}$\cr
\+(ix)&$\wt{SO}(p,q)$  & $\Lambda_{0}$& 
$S_o^{(q-p)/2}\C^p\otimes\C$ & ${q-2\over 2}$ & 
${q-2\over 2}$ &  ${q-p+2\over 2}$ &\cr
\+ & \qquad $3\le p\le q$,   $p+q$ is even\cr \vskip .5pc
\+(x) & $\wt{SO}(3,q)$  & $\Lambda_{0}$ & 
$S^{q-3}\C^2\otimes\C$ & ${q-2\over 2}$ &
${q-2\over 2}$ &  ${q-1\over 2}$ \cr 
\+ & \qquad $4\le q$,  $q$ is even\cr \vskip .5pc
\+(xi) & $\wt{SO}(3,3)$ & $\fo\Lambda_{0}$ & 
$\C^2\otimes\C^2$ & ${1 }$ & ${3\over 2}$ & ${3\over 2}$ \cr
\+(xii)&$\wt{SO}(p,q)$  & none & $*$ & $*$ & $*$ & $*$ \cr      
\+ & \qquad $4\le p\le q$, $p+q$ is odd\cr \vskip .5pc
\+(xiii)&$\wt{SL}(n,\R)$  & $\Lambda_{0}$& $\C$ &${n-2\over 4}$
&${1\over 2}$ &  ${n\over 4}$ \cr
\+ & \qquad $4\le n$,  $n$ is even \cr \vskip .5pc
\+(xiv)&$\wt{SL}(n,\R)$ &  
$\fo\Lambda_{0}$& $\C^n$ &${n\over 4}$ &${3\over 2}$ & 
${n+2\over 4}$ \cr
\+ & \qquad $4\le n$, $n$ is even  \cr \vskip .5pc
\+(xv)&$\wt{SL}(n,\R)$  & none  &   $*$ & $*$ & $*$& $*$   \cr
\+ & \qquad $5\le n$,    $n$ is odd\cr \vskip .5pc
\+(xvi)&$\wt{SL}(3,\R)$  &   
$\Lambda_{0}$ & $\C^2$ & ${1\over 2}$  & ${3\over 4}$ & 
${5\over 4}$ \cr
\+(xvii) &$\wt{SL}(3,\R)$  &  
$\fo\Lambda_{0}$ & $S^3\C^2$ &$1$ & $*$ & $*$ \cr}
\vskip 2pc
In Table 6.9, the symbol $*$ means that there is no entry
because   some  aspect of the construction has failed.
In Cases (xii) and (xv) there is no half-form bundle, while
in Case (xvii) the operators $\pi^z$  fail to satisfy the 
bracket relations of $\sl(3,\R)$.

The final result we present here about our minimal representations
is the computation of a matrix coefficient on the one parameter
subgroup generated by $x=e+\oe\in\pR$. Indeed,  the same 
arguments used in  [B-K4, Th. 6.6] go though in this setting  to give
\theo{6.10}{We have, for $t\in\R$,
$$\<(\exp~tx)\cdot s_0|s_0\>=
\phantom{l}_2F_{1}(a,b; 1+r_0;-\sinh^2t)\eqno(6.33)$$} 
 
\Sec{\S7.  Differential Operators on Half-forms  
and the Generalized Capelli Identity.}

The purpose of this section is to complete the proofs of the results
in \S6, i.e., to prove  Theorems  6.3 and 6.8.
We already reduced Theorem 6.3 to  the operator relation (6.22) on
$H$. 

To begin the proof of (6.22), we observe that the two operators  
$[f_\oe,T_{e}]-[f_e,T\soe]$ and $\L_{\eta^h}$
appearing on the LHS and the RHS of  (6.22) are both
$K_0$-invariant. The first idea of the proof is to exploit  this
observation. Since $K_0$ is reductive, it follows that $H$ is 
completely reducible as
$K_0$-representation. So we can fix a direct sum  decomposition
$H=\oplus_{\al}H_{\al}$  
where each subspace $H_{\al}$ carries an irreducible 
$K_0$-representation. Now to prove the operator relation (6.22) 
holds  on $H_{\al}$  it suffices to prove that (6.22) holds on just one
non-zero section $s^\al$ in 
$H_{\al}$. 

There is a natural method to pick a section  $s^\al$ from $H_{\al}$,  
unique  up to scaling. This uses structure of $H_{\al}$ as a
$\k_0$-representation.  The method is to pick $s^\al$ to be a lowest
weight vector for $\k_0$.
Here we appeal again to the Cartan-Weyl theory recalled 
 in \S5 (where  we applied it to $H$ considered as a 
$\k$-representation).  This time we use lowest weights rather than
highest weights just for convenience.

To get the notion of lowest weight  for $\k_0$, we take the  triangular  
decomposition
$$\k_0=\h\oplus\m_0\oplus\m_{0}^-\eqno(7.1)$$
induced by (5.7).  Then $\h$ is a Cartan subalgebra of $\k_0$ and 
$\b_0=\h\oplus\m_0$ is a Borel subalgebra. So now a 
{\it lowest weight vector} in $H_{\al}$ for $\k_0$  is a vector in the 
$1$-dimensional space  $H_{\al}\mm$.
Then 
$$H\mm=\bigoplus_{\al}H_{\al}\mm\eqno(7.2)$$
is the space of  (i.e., spanned by) all lowest weight vectors in $H$ for 
$\k_0$. 

So proving (6.22) reduces to verifying it on each lowest weight 
vector  $s^\al\in H_{\al}\mm$. By $K_0$-invariance and Schur's
Lemma , the vectors
$([f_\oe,T_{e}]-[f_e,T\soe])(s^{\al})$ and $\L_{\eta^h}(s^{\al})$
again lie in $H_{\al}\mm$. Consequently, for $s=s^{\al}$ we have
$$\eqalign{\L_{\eta^h}(s)=&Xs\cr
([f_\oe,T_{e}]-[f_e,T\soe])(s)=&Rs}\eqno(7.3)$$
where $X$ and $R$ are scalars  depending on $\al$. So   proving 
(6.22) reduces to showing that, for each $\al$,  the scalars $X$ and $R$
coincide. 
 
The second, more serious,  idea is to work out the first idea using the 
Jordan
algebra structure on $\k_{-1}$ from \S4. It turns out that the Jordan
structure gives us (i) a nice way to write down a basis of  $H\mm$
and (ii)  a means to compute $R$ in (7.3) in the form of the 
generalized Capelli Identity of Kostant and Sahi ([K-S]).
The computation of $X$ follows easily from (i).  Then, with everything
computed, we see manifestly that $R=X$.

The rest of this section is devoted to working this out.  
To start off, we embed $H$ in a larger space $\Hs$
which is easier to work with.  We choose the natural embedding
$$H\subset\Hs=\Gamma(\wt{\Yo},\NH)\eqno(7.4)$$
We constructed  the $2$-fold  covering $\wt{\Yo}$  of $\Yo$ in (5.27)
and then the pullback bundle (again called) $\NH$.
Clearly $H$ sits inside $\Hs$ as the space of  sections
which descend to $\Yo$ (i.e., are $(\Z/2\Z)$-invariant) and then
extend to all of $Y$.

In Corollary 5.10 we got an nice description in (5.28) of $\Hs$. 
Since $f_0$ and $\sL$ are $K_0$-semi-invariant, we might as well
rewrite (5.28) as
$$\Hs=\left(\C[\fa,\dots,\fm]\otimes
\C[\fo^{1\over 2},\fo^{-\half}]\right)\sL\eqno(7.5)$$
This makes it clear that a decomposition of 
the polynomial algebra  $\C[\fa,\dots,\fm]$ into irreducible
$K_0$-representations   will produce a decomposition 
$\Hs=\oplus_{\al}\Hs_{\al}$   into irreducible
$K_0$-representations.  In particular we have
$$(\Hs)\mm=\left(\C[\fa,\dots,\fm]^{\m_0^-}\otimes
\C[\fo^{1\over 2},\fo^{-\half}]\right)\sL\eqno(7.6)$$
We will deal with the problem of locating $H\mm$
inside $(\Hs)\mm$ when the time comes.

Now we can bring in  the Jordan algebra $\k_{-1}$.
We have a vector space isomorphism
$\k_{-1}\to\p_{1}$, $y\mapsto[y,e]$; cf. Lemma 3.4. This induces a 
graded complex algebra isomorphism
$$S(\k_{-1})\to \C[\fa,\dots,\fm],\qquad
g\mapsto\hat{g}\eqno(7.7)$$
defined in degree $1$ by $\hat{y}=f_{[y,e]}$ for $y\in\k_{-1}$
where $f_v$ was defined in (3.1).
This isomorphism has weight $\chi^p$ in degree $p$ under the  action
of $K_0$.  Hence (7.7) is $K_0'$-equivariant  and so  gives by
restriction   a  graded complex algebra isomorphism
$$S(\k_{-1})\mm\to \C[\fa,\dots,\fm]\mm,\qquad
g\mapsto\hat{g}\eqno(7.8)$$

Recall from (4.12) that $q_n$ is the degree of the Jordan algebra 
$\j_{[n]}$. We put $c_1=0$ and
$$c_n=q_1+\cdots+ q_{n-1},\qquad n=2,\dots,\ell\eqno(7.9)$$

We have the well-known result 
(see   [B-K3,  Theorem 3.4 and Corollary 3.6])

\lem{7.1}{The natural representation of $K_0$ on $S(\k_{-1})$ is 
completely reducible and multiplicity free.
\mypar
The  ring of lowest weight vectors in
$S(\k_{-1})$ for the $K_0$-action is  a polynomial algebra in $q$
independent graded generators so that
$$S(\k_{-1})\mm=\C[N_1,\dots,N_q]\eqno(7.10)$$ 
The polynomials $N_1,\dots,N_q$ are  uniquely determined   
by the conditions {\rm(i)}  for   $1\le j\le q_n$,  we have
 $N_{c_n+j}\in S^j(\j_{[n]})$ and  {\rm(ii)}   
$N_1(\ek)=\cdots=N_q(\ek)=1$.}

Now combining  the work we have done thus far in this section we 
obtain
\prop{7.2}{The ring of lowest weight vectors in $\C[\fa,\dots,\fm]$   
for the $K_0$-action is  a polynomial algebra in the $q$ independent 
graded generators $\hat{N}_1,\dots,\hat{N}_q$ so that
$$\C[f_1,\dots,f_m]^{\m_{0}^{-}}=\C[\hat{N}_1,\dots,\hat{N}_q]
\eqno(7.11)$$ 
Then $\deg\hat{N}_{c_n+j}=j$ for    $1\le j\le  q_n$.
\mypar
Suppose $\NH$ is a half-form bundle on $Y$ and $H=\Gamma(Y,\NH)$.
Then the natural $K_0$-action on $H$ is completely reducible
and has  a basis of lowest weight vectors of the form
$$s=\fo^p\hat{N}^{t_1}_1\cdots\hat{N}^{t_q}_q\sL\eqno(7.12)$$
where   $p\in\half\Z$ and  $t_1,\dots,t_q\in\ZP$. Here we regard $s$ as
a section in the larger space $\Hs$. The section $s$ in {\rm(7.12)}  
determines  $p,t_1,\dots,t_q$ uniquely.
\mypar
The $K_0$-representation on each $E'$-eigenspace $H_{r_0+n}$
is multiplicity free.} 
 
If we range  over all tuples  $p,t_1,\dots,t_q$ in (7.12), then we
obtain a   basis of   $(\Hs)\mm$.   
We could   explain how to decide when the corresponding
section lies in $H$, but we omit this as it is not necessary in our work
below. (Only the  partial answer  we give later is needed.)

Now that we have obtained a nice basis of $H\mm$ as promised, we 
need to start computing the eigenvalues of our various operators on
the section $s$ in (7.12).  We set $\bt=(t_1,\dots,t_q)$. We put
$$N^{\bt}=N^{t_1}_1\cdots N^{t_q}_q\AND
z=\deg N^{\bt}=\sum_{i=1}^qt_i\deg
N_i\eqno(7.13)$$   Recall  from (3.14) that $m=\dimC Y -1$.

\lem{7.3}{Suppose $s\in H$ is of the form {\rm(7.12)} Then 
$$\eqalign{E's=&rs\,\qquad\hbox{where}
\qquad r=p+z+{m+1\over 2}\cr
\L_{\eta^h}s=&Xs\qquad\hbox{where}\qquad
X=2p+z+{m+2\over 2}}\eqno(7.14)$$}
\proof{We have
$$E'(s)=(E\fo^p\hat{N}^{\bt})\sL
+\fo^p\hat{N}^{\bt}(E'\sL)=(p+z)s+{m+1\over 2}s\eqno(7.15)$$ 
using (7.12) and  the  formulas (3.19) and (5.17) for  $E$ and $\Lambda$
in terms of  our local coordinates. Next
$$\eta^h=2\fo\pdb{\fo}+\sum_{i=1}^m f_i\pdb{f_i}\eqno(7.16)$$
Then computing the same way as in (7.15) we find 
$$\L_{\eta^h}(s)=(\eta^h\fo^p\hat{N}^{\bt})\sL
+\fo^p\hat{N}^{\bt}(\L_{\eta^h}\sL)=
(2p+z)s+{m+2\over 2}s\eqno(7.17)$$
\qed}

Now we have come to the tricky part, computing the eigenvalue $R$
of $[f_\oe,T_{e}]-[f_e,T\soe]$ on each section $s$ from (7.12). 
Expanding out we get
$$\eqalign{&[f\soe,T_{e}]-[f_e,T\soe]=
f\soe T_{e}-T_{e}f\soe-f_eT\soe+T\soe f_e\cr
&={1\over (E'-1)E'}f\soe D_e -{1\over E'(E'+1)}D_ef\soe
-{1\over (E'-1)E'}f_e D\soe +{1\over E'(E'+1)}D\soe f_e}\eqno(7.18)$$  

We will use the generalized Capelli Identity of Kostant and Sahi ([K-S])
to compute the eigenvalues
of  $f\soe D_e$,  $D_ef\soe$, $f_eD\soe$, $D_\oe f_e$, , and  on $s$.
The idea is to transform this computation   into a 
computation on $S(\k_{-1})$. This works   because we 
will write everything in terms of  our local coordinates 
$\fo,\fa,\dots,\fm$ and use   (7.7).

The first thing we will compute is $f\soe D_e(s)$.
We already have the expression (6.11) for $D_e$ in terms of our local
coordinates. So we need the  expression for the function $f\soe$.
In analogy to  (4.13)  we set
$$N=N_{q_1}^{w_1} N_{q_1+q_2}^{w_2}
\cdots N_{q_1+\cdots+q_{\ell}}^{w_\ell} 
\eqno(7.19)$$

Complex conjugation on $\g$ preserves $\Omin$ and $\p$ and so
preserves $Y$ by (3.7). This induces naturally a complex conjugation
map $f\mapsto\of$ on $R(Y)$. This is then a 
$\C$-anti-linear  real algebra involution. 

\prop{7.4}{The unique expression for $f\soe$ in terms of our local
coordinates $\fo,\fa,\dots,\fm$ on $Y$ is 
$$f_{\oe}=\ofo={\hat{N}\over\fo^3}\eqno(7.20)$$}
\proof{This is proven  in the same way
 as  in [B-K4, Prop. 4.3].\qed}

To state and prove our  computation of $f\soe D_e(s)$, we need to
encode the monomial $N^{\bt}$ into a new $q$-vector, namely
the {\it multi-degree}  $\mdeg(N^{\bt})$. We define this by
$$\mdeg (N_{c_n+j})=(0,\dots,0,2,\dots,2,0,\dots,0),\qquad
1\le j\le q_{n}\eqno(7.21)$$ where
there are   $c_n$ zeroes follows by $j$ twos followed by zeroes, and
$$\mdeg(N^{\bt})=\sum_{i=1}^qt_i\mdeg(N_i)\eqno(7.22)$$
where addition of vectors is component-wise.  If
$\mu=(\mu_1,\dots,\mu_q)=\mdeg(N^{\bt})$ then
$$\mu_1+\cdots+\mu_q=2\deg N^{\bt}=2z\eqno(7.23)$$ 

Let $\delta$ be the $q$-vector such that 
$\delta_{c_n+j}=d_n(q_n-j)$ for $1\le j\le q_n$
where we recall the root multiplicity numbers $d_n$ from the 
proof of  Proposition 5.7. So
$$\eqalign{\delta&=(\delta_1,\dots,\delta_q)\cr
&=(d_1(q_1-1),\dots,d_1,0,d_2(q_2-1),\dots,d_2,0,\dots
d_{\ell}(q_{\ell}-1),\dots,d_{\ell},0)}\eqno(7.24)$$
Finally we define the $q$-vector
$$\bv=(v_1,\dots,v_q)=(w_1,\dots,w_1,\dots,w_{\ell},\dots,w_{\ell})
\eqno(7.25)$$
where each $w_n$ occurs $q_n$ times.

From now on, we assume that   $s\in H$ is of the form 
$s=f_0^p\hat{N}^{\bt}\sL$ of {\rm(7.12)}  and 
$\mu=\mdeg({N}^{\bt})$.
 
\prop{7.5}{We have
$$f\soe D_e(s)=
s\prod_{i=1}^q\prod_{j=0}^{v_i-1}  C_{i,j}(\mu)\eqno(7.26)$$
where $C_{i,j}(\mu)$ is the  Capelli multiplier given by
$$C_{i,j}(\mu)={\mu_i+\delta_i-2j\over 2v_i} \eqno(7.27)$$}
\proof{First (6.11) and (7.20) give
$$f\soe D_e=\hat{N}P(\L_{\pa_1},\dots,\L_{\pa_m})\eqno(7.28)$$
where $\pa_k=\pdb{f_k}$ for $k=1,\dots,m$. 
To compute $f\soe D_e(s)$ we first compute $\L_{\pa_k}(s)$. 
But  $\pa_k(f_0)=0$  and also taking the Lie derivative of (5.17) we 
get  $\L_{\pa_k}(\sL)=0$.   So we get
$$\L_{\pa_k}(s)=\pa_k(f_0^p\hat{N}^{\bt})\sL
+f_0^p\hat{N}^{\bt}(\L_{\pa_k}\sL)=f_0^p(\pa_k\hat{N}^{\bt})\sL
\eqno(7.29)$$
Then (7.28) gives
$$f\soe D_e(s)=f_0^{p}\hat{N}
\left(P({\pa_1},\dots,{\pa_m})(\hat{N}^{\bt})\right)\sL
\eqno(7.30)$$ 
\mypar
Next we introduce the  graded complex algebra isomorphism
$$S(\k_1)\to\C\left[\pa_1,\dots,\pa_m\right],
\qquad B\mapsto\pa_{B}\eqno(7.31)$$
defined in degree $1$ by  $\pa_x=\fo\i\eta^x$,
$x\in\k_1$,   so that, by (3.19), $\pa_{x_k}=\pdb{f_k}$.
Then we can rewrite (7.30) as
$$f\soe D_e(s)=f_0^{p}\hat{N}\left(\pa_P(\hat{N}^{\bt})\right)\sL
\eqno(7.32)$$ 
\mypar
It follows easily, as in [B-K4, (4.4.4)],  that
$$\hat{N}\pa_P(\hat{N}^{\bt})=
\widehat{N\pa_P(N^{\bt})}\eqno(7.33)$$
Indeed, it suffices to check that $\pa_x\hat{y}=\pa_xy$
where $x\in\k_1$ and $y\in\k_{-1}$. We find
$\pa_x\hat{y}=f_e\i f_{\psi([x,y])e}=\psi([x,y])=(x,y)_{\g}=\pa_xy$.
\mypar
Now  we can compute $N\pa_P(N^{\bt})$  using the generalized  Capelli
Identity of Kostant and Sahi [K-S]. This is similar to 
[B-K4, proof of Th. 4.4], but we are in a more general situation here
where  the Jordan algebra $\k_{-1}$ is not necessarily simple and 
also the multiplicities $w_1,\dots,w_{\ell}$ may be non-trivial.
The point is that $N\pa_P(N^{\bt})$ breaks into a product with one
factor for each simple component $\j_{[n]}$ of $\k_{-1}$. The $n$th
factor is  
$$N_{c_n+q_n}^{w_n}\pa_{P_{[n]}^{w_n}}
\left(N_{c_n+1}^{t_{c_n+1}}\cdots N_{c_n+q_n}^{t_{c_n+q_n}}\right)   
\eqno(7.34)$$
The Capelli Identity for $\j_{n}$ says that the operator
$N_{c_n+q_n}^{w_n}\pa_{P_{[n]}^{w_n}}$ acts on it argument in (7.34)
by a scalar and computes that scalar.  
\mypar
Putting all the factors together  we obtain
$$N\pa_P(N^{\bt})=N^{\bt}\prod_{i=1}^q
\prod_{j=0}^{v_i-1}C_{i,j}(\mu)\eqno(7.35)$$
which then gives (7.26) because of (7.32) and (7.33).
\mypar
There is one subtle point here: the
appearance of the factor $v_i\i$ in $C_{i,j}(\mu)$. The 
corresponding Capelli multiplier from [K-S] is just
$\half(\mu_i+\delta_i-2j)$.   However the factor $v_i\i$
arises because of the way we have paired $\k_{-1}$ with $\k_{1}$.  
Let $\ek=\ek^1+\cdots+\ek^{\ell}$ be the decomposition of $\ek$  
corresponding to (4.11) so that $\ek^n$ is  the Jordan identity
element in $\j_{[n]}$; similarly  we get
$\oek=\oek^1+\cdots+\oek^{\ell}$.  Then we easily find
$$(\ek,\oek)_{\g}=\sum_{n=1}^{\ell} (\ek^n,\oek^n)_{\g}\AND
 (\ek^n,\oek^n)_{\g}=q_nw_n\eqno(7.36)$$
This fits with (4.14) since $(\ek,\oek)_{\g}=4$ just as in [B-K4, proof of
Theorem 4.4.)  So we get $\pa_{\ek^n}(\oek^n)=q_nw_n$.
However the normalization from [K-S] is that 
$\pa_{\ek^n}(\oek^n)=q_n$. The ratio  $w_n$ then appears in the
denominator of our Capelli multiplier.  \qed} 

Let $\Xi$ be the set of  ordered pairs $(i,j)$ occurring in (7.26).
The cardinality of   $\Xi$ is $q_1w_1+\dots+q_nw_n=4$ by (4.14).

\cor{7.6}{We have
$$D_ef_{\oe}(s)=s\prod_{(i,j)\in\Xi}\big(C_{i,j}(\mu)+1\big)
\eqno(7.37)$$ }
\proof{As in the last proof we  find, cf. (7.28), 
$$ D_ef\soe=P(\L_{\pa_1},\dots,\L_{\pa_m})\hat{N}\eqno(7.38)$$
and so, as in  (7.32) and (7.33), 
$$D_ef\soe (s)=f_0^{p}\widehat{\left(\pa_P(N{N}^{\bt})\right)}\sL
\eqno(7.39)$$
Then we find
$$\pa_P(N{N}^{\bt})=\prod_{(i,j)\in\Xi}C_{i,j}(\mu')\eqno(7.40)$$
where $\mu'=\mdeg(N{N}^{\bt})$. But   $\mu'=2\bv+\mu$ and so
$C_{i,j}(\mu')=C_{i,j}(\mu)+1$.  Now we get (7.37).
\qed}

Next we want compute  $f_e D_\oe(s)$.
We solve this as in [B-K4] by introducing, in the next result, an
involution $\theta$ of $H$.   We construct this using the group element 
$$\theta_o={\bf\exp}~{\pi\over 2}(\ek-\oek)\in
K\eqno(7.41)$$ 
where ${\bf\exp}:\k\to{K}$ is the exponential  map.  
The same arguments used in  [B-K4, Lem. 4.6, Prop.
4.6 and Th. 4.7] give   

\lem{7.7}{The action of $\th_o$ on $\p$  preserves $Y$ and defines a
graded complex algebra involution $\th$ of $R(Y)$ which commutes
with  complex conjugation.  We have $\fo^\th=\pm\ofo$. 
\mypar
The natural action of $\th_o$ on the $K$-homogeneous half-form
bundle  $\NH$ over $Y$ defines a complex linear involution  
$\th:H\to H$. Then $\th:H\to H$ preserves the grading {\rm(5.15)} and is 
compatible with the
$R(Y)$-module structure so that $(fs)^\th=f^\th s^\th$. 
\mypar
$\th$ permutes the simple $\k_0$-submodules in  $H$ and 
moves lowest weight  vectors to highest weight vectors. For any $s\in
H$ we have
$$f_eD_\oe(s)=(f_\oe D_e(s^\th))^\th\eqno(7.42)$$ }
 
We can now prove 
\prop{7.8}{We have
$$f_eD_\oe(s)=s\prod_{(i,j)\in\Xi}\big(r-1-C_{i,j}(\mu)\big)
\eqno(7.43)$$}
\proof{Lemma 7.7 reduces the problem to computing $f_\oe
D_e(s^\th)$.   Now   $s$ is a lowest weight  vector in some  simple
$\k_0$-submodule
$F$ in $H$, and so  Lemma 7.7 implies that  $s^\th$ is a highest  weight 
vector in the simple $\k_0$-submodule $F^\th$ in $H$. Then there is a 
lowest weight vector $s^*\in F^\th$ (unique up to scaling). We can 
write $s^*=\fo^{a}\hat{u}\sL$ where
$u\in (S^{b}(\k_{-1}))\mm$.  Let $\mu^*=\mdeg(u)$. Then 
$$f_\oe D_e(s^\th)=s^\th\prod_{(i,j)\in\Xi}C_{i,j}(\mu^*)\eqno(7.44)$$
\mypar
We claim that there is an involution $(i,j)\mapsto(i^*,j^*)$ on
the set $\Xi$  such that
$$C_{i^*,j^*}(\mu^*)=r-1-C_{i,j}(\mu)\eqno(7.45)$$
This, because of Lemma 7.7, gives (7.43). 
\mypar
The construction of the involution  requires several calculations. To
carry  these out, we  bring to the forefront the theory of
weights associated to our triple $(\k,\k_0,\kl)$
from \S4.   Indeed,  our choice of $(\h,\b)$ in \S5 was compatible
with complex conjugation and 
the complex Cartan decomposition $\k_0=\kl\oplus\r$ so that
we have a complex conjugation stable splitting  $\h=\a\oplus\t$  where
$\a\subset\r$ is  a maximal abelian subalgebra and $\t\subset\kl$. 
\mypar
Now  $\dimC\a=q$ and 
$\a^*$ has a unique basis $\vep_1,\dots,\vep_q$ such that the
$\a$-weight of $N_{c_n+j}$, where $1\le j\le q_n$, is 
$-2(\vep_{c_n+1}+\dots+\vep_{c_n+j})$.
The weights $\ep_i$ are pure imaginary and the  action of  $\th_o$ 
gives the complex  involution $\th$ of  $\k_0$ with fixed algebra $\kl$
and $(-1)$-eigenspace $\r$. In particular,  $\th_o$ acts as $-1$ on 
$\a$.  It follows that
$$\sig=\hbox{$\a$-weight  of  $s$}\quad\Rightarrow\quad
-\sig=\hbox{$\a$-weight  of  $s^\th$}   \eqno(7.46)$$ 
\mypar
Let $\nu\mapsto\nu^\tau$ be the involution of $\a^*$ which exchanges
the highest weight of a simple $\k_0$-submodule in $S(\k_{-1})$
with the lowest weight.  Then 
$$-\sig^\tau=\hbox{$\a$-weight  of  $s^*$}\eqno(7.47)$$
In terms of our basis of $\a^*$ we have
$\vep_{c_n+i}^\tau=\vep_{c_n+q_n-i+1}$ for $1\le i\le q_n$. 
The $\a$-weight of  $f_0=f_e$ is $\psi|_\a$, which we
again call $\psi$. From now on, we identify a $q$-vector
$\nu=(\nu_1,\dots,\nu_q)$ with the $\a$-weight 
$\nu=\sum_{i=1}^q\nu_i\vep_i$.
\mypar
We now define our involution on $\Xi$ by  
$$\vep_{i^*}=\vep^{\tau}_i\AND
j^*=v_i-j-1\eqno(7.48)$$
(This  automorphism of $\Xi$ is 
can happen to be the  identity, so by involution we mean just
that the order divides $2$.)  Then $\mu_i^\tau=\mu_{i^*}$.
\mypar
Let us put $g=N^{\bt}$ so that 
$s=f_0^p\hat{g}\sL$ ,  $\mu=\mdeg(g)$ and $g\in(S^z(\k_1))\mm$. 
Then the  $\a$-weight  of $\hat{g}$ is  $z\psi-\mu$.   Also
$$\la=\hbox{$\a$-weight of}~\sL={m+1\over 2}\psi-{\kap\over 2}
\eqno(7.49)$$ 
where $-\kap$ is the weight  of $\a$ on
$\wedge^m\k_{-1}$. The sum of the $\a$-weights of
$f_0^p$, $\hat{g}$ and $\sL$ is   
$$\sig=(p+z)\psi-\mu+\la=r\psi-\mu-{\kap\over 2}\eqno(7.50)$$
\mypar
Now  we can compute $\mu^*$.  The sections
$s=\fo^p\hat{g}\sL$ and $s^*=\fo^a\hat{u}\sL$
have the same eigenvalues  under $E'$ and $\L_{\eta^h}$ which
means that $p+z=a+b$ and $2p+z=2a+b$. Hence $p=a$ and $z=b$.
So the $\a$-weight of $\hat{u}$ is $z\psi-\mu^*$ and we get
 $$-\sig^\tau=\hbox{$\a$-weight  of  $s^*$}=
(p+z)\psi-\mu^*+\la\eqno(7.51)$$ 
Applying $\tau$ to   (7.49) and subtracting  (7.50)  we  get
$$\mu^*=2\sig^\tau+\mu^\tau=2r\psi-\mu^\tau-\kap\eqno(7.52)$$
\mypar
To obtain a proof of (7.45), we write out (7.52) in terms of 
components. The $\a$-weight $\psi$ has the same components 
as the vector $\bv$ in (7.25).   A key observation is that the
components   of $\kap$ are $\kap_i=2+\delta_i+\delta_{i}^\tau$.  
To see this,  we start from the fact that 
$\kap$ is the sum of  the $\a$-weights 
$\kap^{[n]}$ of the  top exterior powers of  the spaces  
$\j_{[n]}$.  The weights of $\a$ on $\j_{[n]}$ are
precisely  the weights  $2\vep_i$ and $\vep_i+\vep_j$ where
$c_n+1\le i<j\le c_n+q_n$. So
$\kap^{[n]}=2\sum_i\vep_i+d_n\sum_{i<j}(\vep_i+\vep_j)$ and
this gives our  formula for $\kap_i$. Now (7.52) gives
$$\mu^{*}_{i}+\delta_{i}=2rv_i-2-(\mu^{\tau}_i+\delta^{\tau}_i)
\eqno(7.53)$$
\mypar
Now subtracting $2j$ from both sides of (7.53) and using (7.48) we 
get
$$\mu^{*}_{i}+\delta_{i}-2j=2(r-1)v_i-
(\mu_{i^{*}}+\delta_{i^{*}}-2j^{*})\eqno(7.54)$$
Dividing through by $2v_i$ (notice $v_i=v_{i^*}$) we get  (7.45).\qed} 
Arguing as in the proof of Corollary 7.6 we get
\cor{7.9}{We  have 
$$D_\oe f_e(s)=s\prod_{(i,j)\in\Xi}\big(r-C_{i,j}(\mu)\big)
\eqno(7.55)$$}

We  can now compute  the scalar $R$ from (7.3); we already
computed $X$ in (7.14).  Starting with (7.18) and plugging in
(7.26), (7.37),  (7.43) and (7.55) we get
$$R={\prod C_{i,j}(\mu) \over (r-1)r}
-{\prod (C_{i,j}(\mu)+1)\over r(r+1)}
-{\prod (r-1-C_{i,j}(\mu))\over (r-1)r}
+{\prod (r-C_{i,j}(\mu))\over r(r+1)}
\eqno(7.56)$$
This is only valid when $r\neq 1$ (as we know $r>0$).

Fortunately, the expression for $R$ in (7.56) simplifies greatly. 
We can apply the following  formal identity given 
[B-K4, Lem. 4.8]. Put
$$J(a_i;b)=J(a_0,a_1,a_2,a_3;b)={a_0a_1a_2a_3
\over b(b+1)}\eqno(7.57)$$
and $a'_{m}=b-a_m$ where $b,a_0,a_1,a_2,a_3$ are five 
indeterminates. Then
$$\eqalign{&J(a_i;b)-J(a'_{i};b)-J(a_i+1;b+1)+J(a'_{i}+1;b+1)
\cr&=2b-(a_0+a_1+a_2+a_3)}\eqno(7.58)$$ 
Applying this with $b=r-1$ and $a_1,a_2,a_3,a_4$ set equal to
the four Capelli multipliers $C_{i,j}(\mu)$ (taken in any order), we get
$$ R=2r-2-\sum_{(i,j)\in\Xi} C_{i,j}(\mu)\eqno(7.59)$$
To prove (6.22), we need to show $R=X$. We compute   
$$\eqalign{\sum C_{i,j}(\mu)
=&\sum_{i=1}^q\sum_{j=0}^{v_i-1}{\mu_i+\delta_i-2j\over 2v_i}
=\sum_{i=1}^q{\mu_i+\delta_i-v_i+1\over 2}\cr
=&z+{m\over 2}-2=2r-2-X}\eqno(7.60)$$
The second to  last equality follows as $\sum_{i=1}^q\mu_i=2z$, 
$\sum_{i=1}^q\delta_i=\sum_{n=1}^ld_n(q_n-1)q_n/2=m-q$
and $\sum_{i=1}^q v_i=\sum_{n=1}^{\ell}q_nw_n=4$, while the last
follows by (7.14). So  $R=X$.   Thus (6.22)  holds on $s$ if $s\notin H_1$. 
This proves Theorem 6.3 in all cases where $r=1$ never occurs in the
spectrum  of  $E'$ on $H$. 
\mypar
Now suppose that $r=1$ does occur and $s\in H_1$. Then  
Proposition  5.5 implies  that $r_0=1$ so that $H_1$ is the vacuum space.  
Thus by  $E'$-degree $D_es=D_\oe s=0$. Then (7.26) implies that 
$\prod  C_{i,j}(\mu)=0$ and so  $C_{i,j}(\mu)=0$ for some $(i,j)$.
But also (7.56) collapses to
$$\eqalign{R&=\half\prod(1-C_{i,j}(\mu))-\half\prod(C_{i,j}(\mu)+1)
\cr&=-\sum_{(i,j)} C_{i,j}(\mu)-\sum_{(i,j)\neq(i',j')\neq(i'',j'')}
C_{i,j}(\mu)C_{i',j'}(\mu)C_{i'',j''}(\mu)}\eqno(7.61)$$ 
But (7.60) gives $\sum C_{i,j}(\mu)=-X$.  Hence
$R=X$ if and only if the third elementary symmetric 
function of the four numbers $C_{i,j}(\mu)$ is zero. But we already 
know at least one $C_{i,j}(\mu)$ vanishes. So $R=X$ if and only if 
at least  two of the four numbers $C_{i,j}(\mu)$  are zero.
\mypar
At this point we observe that $r=r_0=1$ implies something very
particular about the form of   $s$: in our normal form
$s=\fo^p\hat{N}^{\bt}\sL$ we have 
$${N}^{\bt}=N_{q_1}^{w_1-u_1}N_{q_1+q_2}^{w_2-u_2}\cdots
N_{q_1+\cdots+q_\ell}^{w_\ell-u_\ell}\eqno(7.62)$$ 
where $0\le u_n\le w_n$. Now it 
follows  from (7.27) that  the list of  four numbers $C_{i,j}(\mu)$ has at
least  $\ell$ zeroes. Thus we are left with the case where $\ell=1$. 

Suppose $\ell=1$.  Then  $m\le 4$.   To see this, we consider
a   highest weight vector  $s_1\in H_1$ for the $\k$-action.
Then $s_1=\fo^{-j}\sL$ for some $j\in\half\Z$. The eigenvalue
of $\L_{\eta^h}$ on $s_1$ is $-2j+(m+2)/2$  and must be non-negative.
The eigenvalue of $E'$ on $s_1$ is $-j+(m+1)/2$ and is equal to $1$.
But then $(m+2)/4\ge j=(m+1)/2-1$ and so $4\ge m$.

Now  looking at  Table 4.6, we see that  $\ell=1$ and $m\le 4$ only if
$\gR=\sl(p,\R)$ where $p=6,5,$ or $3$. 
We rule out $p=5$ because of Proposition 5.7(ii).
For $p=6$ we have $r_0=1$ when $s_{0}^2=\Lambda_0$ and
$H_1\simeq\C$. Then $s=s_0$ and $\mu=(\mu_1,\mu_2)=(0,0)$.
The multipliers $C_{i,j}(\mu)$ are 
${\mu_1+p-4\over 2},{\mu_2\over 2},{\mu_1+p-6\over 2},
{\mu_2-2\over  2}$ and this list  has  two zeroes as required.
For $p=3$, we have  $r_0=1$ when $s_{0}^2=\fo\Lambda_0$ and
$H_1\simeq S^3\C^2$.   Then $s$ is one of four vectors
with $\mu=(\mu_1)$ where $\mu_1=0,2,4,$ or $6$. 
The multipliers $C_{i,j}(\mu)$ are 
${\mu_1\over 4},{\mu_1-2\over 4},{\mu_1-4\over 4},
{\mu_1-6\over  4}$  and so we never get two zeroes in this list.
Thus this one case fails to produce a representation.
This concludes the proof of  Theorem 6.3. 
\vskip 1pc
Next  we  finish the proof  of  Theorem 6.8.
We started this in \S6,  and left off at (6.32) where we needed
to compute the numbers $\ga_k$ defined by
$T_{\oe}(f_{0}^ks_0)=\ga_kf_{0}^{k-1}s_0$.
But now we can compute the $\ga_k$ because of 
Proposition 7.8. Indeed, let $s=f_{0}^ks_0$; then
$r=r_0+k$ and $\mu=\mdeg(s)=0$. Now   (7.4.3) gives
$$T\soe(f_{0}^ks_0)=f_{0}^{k-1}s_0
{\prod_{(i,j)\in\Xi}\big(r_0+k-1-C_{i,j}(0)\big)\over 
(r_0+k-1)(r_0+k)}\eqno(7.63)$$

We can simplify the factor  $\ga_k$ appearing in (7.63) by
computing the four numbers $C_{i,j}(0)$.
We have $T_e(s_0)=T_{\oe}(s_0)=0$ since $s_0$ is a vacuum vector.
Hence, if   $r_0\neq 1$ then the list of four
multipliers $C_{i,j}(0)$ contains zero and  $r_0-1$. Moreover
we just showed that  if $r_0=1$ then
the list contains zero  with multiplicity at least two.
Thus, regardless of the value of $r_0$, we can write the list of four
multipliers $C_{i,j}(0)$ as $0,r_0-1,r_0-a,r_0-b$ where $a$ and $b$ are
unknown. Then (7.63) gives 
$$\ga_k= {k(k-1+a)(k-1+b)\over (r_0+k)}\eqno(7.64)$$
Consequently
$$\ga_1\cdots\ga_n={n!(a)_n(b)_n\over (r_0+1)_n}\eqno(7.65)$$
Because of  (6.32) this gives (6.30). In fact we have gotten   the 
more precise information
\prop{7.10}{We have equalities of multi-sets:
$$\{C_{i,j}(0)\}_{(i,j)\in\Xi}=
\left\{{\delta_i-2j\over 2v_i}\right\}_{(i,j)\in\Xi}=
\{0,r_0-1,r_0-a,r_0-b\}\eqno(7.66)$$
In this way  $P$ and $r_0$  determine uniquely the  numbers $a$ and
$b$ appearing in Theorem   {\rm 6.8}. } 
 
Notice that the choice of $\gR$ determines $P$ while the choice of
$\NH$ determines $r_0$.

Proposition 7.10 implies $\sum C_{i,j}(0)=3r_0-1-a-b$.
But also  $\sum C_{i,j}(0)=2r_0-2-X_0$ by (7.60). Comparing, we get
(6.31).

Finally we can use Proposition 7.10 to   compute the numbers
$a$ and $b$ in Table 6.9. We  then observe that $a$ and $b$ are always
positive. But we also have a nice  theoretical proof of the positivity.

By (7.66), the four numbers $r_0-C_{i,j}(0)$ are  $1,r_0,a,b$. 
Using first (7.45) and then (7.27)  we can write
$$r_0-C_{i,j}(0)=1+C_{i^*,j^*}(\al)=
{\al_{i^*}+\delta_{i^*}+2v_{i^*}-2j^*\over 2v_{i^*}}
\eqno(7.67)$$
where $\al=0^*$.  The last expression in (7.67) is positive since
$\al_{i^*},\delta_{i^*}\ge 0$ and $v_{i^*}>j^*$.
This concludes the proof of Theorem 6.8.

\Sec{\S 8. The Reproducing Kernel of $\H$.}

The aim of this section is to show that the 
Hilbert spaces $\H$ carrying the unitary
irreducible representations   constructed in  \S6 each admit 
a reproducing kernel $\K$ and $\K$ is a holomorphic
half-form on $\YoY$. It follows then that $\H$ consists entirely of
holomorphic half-forms on  $Y$.

We work in the setting of Theorems 6.3, 6.6 and 6.8.  (So the
one case $\GR=SL(3,\R)$, $r_0=1$,
$H_1\simeq\C^4$ is excluded as this case did not quantize.) 

To begin with, we explain how the notion of reproducing kernel 
applies here. Our Hilbert space $\H$ is the completion of 
$H=\Gamma(Y,\NH)$.
Therefore, using the grading  (5.15) of $H$, we may
regard a section  $s\in\H$ as a formal sum 
$$s=\sum_{n\in\ZP}s_n\eqno(8.1)$$
where $s_n\in H_{r_0+n}$.  Then $s$ is a holomorphic section of
$\NH$ if and only if the series in (8.1) converges locally uniformly.

The complex conjugate  space $\oH$ identifies naturally with
the space $\Gamma(\oY,\oNH )$
of algebraic holomorphic sections of $\oNH$ over $\oY$.
Here $\oNH$ is the complex conjugate line over the complex  
conjugate algebraic manifold $\oY$.    So we get an identification
$$H\otimes\oH=\Gamma(\YoY,\NH\otimes\oNH)\eqno(8.2)$$

A {\it reproducing kernel} for $\H$ is a   section
$\K$ of $\NH\otimes\oNH$ over $\YoY$
such that  for each $v\in Y$,  the formula
$$\K_v(u)=\K(u,\ov)\eqno(8.3)$$
defines a section  $\K_v\in\H\otimes\oNH_{\ov}$ 
and  we have the ``reproducing" property for all $s\in\H$
$$s(v)=\<s|\K_v\>\eqno(8.4)$$ 
This makes sense as both sides of (8.4) define vectors in 
$\NH_v$.

$\H$ admits a reproducing kernel if and only
if the evaluation map $s\mapsto s(v)$
is continuous on $\H$ for every point $v\in Y$.
A reproducing kernel on $\H$, if it exists,
is unique and is computed by 
$$\K =\sum_{k}g_k\otimes\ovl{g_k}\eqno(8.5)$$
where $\{g_k\}$ is  any orthonormal basis of $\H$.
See , e.g., [F-K, IX,\S2] for the case of Hilbert spaces of
holomorphic functions.

Each space  $H_{r_0+n}$, $n\in\ZP$, is finite-dimensional and 
so admits a reproducing kernel 
$\Pi_n\in\Gamma(\YoY,\NH\otimes\ovl{\NH}\,)$. The
reproducing kernel
$\K$ of $\H$ exists if and only if $\K=\sum_{n\in\ZP}\Pi_n$,
i.e., if and only if the series $\sum_{n\in\ZP}\Pi_n$ converges.

We have a $K$-invariant function $T\in R(\YoY)$ defined by
$$T(u,\ov)=(u,\ov)_{\g}\eqno(8.6)$$
Then $u\mapsto T(u,\ou)$ is a positive  real function on $Y$.

\theo{8.1}{For any orthonormal basis $\{g_k\}$  of $H$,
the series in {\rm(8.5)} converges locally uniformly  and 
moreover we have the formula
$$\K=\phantom{i}_1F_2 (r_0+1;a,b;T)\Pi_0\eqno(8.7)$$ 
where $a$ and $b$ are as in Theorem {\rm6.8} and $\Pi_0$ is the 
reproducing kernel of $H_{r_0}$.
Consequently, $\K$ is a holomorphic section
$$\K\in\Gamma^{hol}(\YoY,\NH\otimes\oNH)\eqno(8.8)$$}
\proof{For each $n$, $\Pi_n$ is a $K$-invariant section of
$\NH\otimes\oNH$. This follows as 
the Hermitian inner product on $H_{r_0+n}$ is $\KR$-invariant.
Hence the quotient $\Pi_n/\Pi_0$ is a $K$-invariant rational
function on $\YoY$. We have   natural actions of $K\times\C^*$
on $Y$ and $\oY$ where $\C^*$ acts by the Euler scaling action.
\lem{8.2}{The product action  of $K\times\C^*$ on the
variety $\YoY$ has a unique Zariski dense orbit
$W$. The function $T$ separates the $K$-orbits in $W$. 
\mypar
Moreover any $K$-invariant rational function on $\YoY$ is a
polynomial in $T$ and $T\i$.}
\proof{The isotropy group of $K$ at $(e,\oe)$ is 
$K^e\cap K^\oe=\Ks=K_0'$. 
So the $K$-orbit of $(e,\oe)$ is isomorphic to $K/K_0'$ and hence
has codimension $1$ in $\YoY$ by Theorem 4.1 since 
$\dimC\YoY=2\dimC Y=\dimC O$. 
The $\C^*$-orbit of $(e,\oe)$ and the  $K$-orbit of $(e,\oe)$
meet in exactly two points, namely $\pm(e,\oe)$.
This follows from the easy fact that $(a\cdot
e,a\cdot\oe)=(te,t\oe)$ if and only if  $a\in K^h$ and
$\chi(a)=t=\chi(a)\i$.
\mypar
In particular then the $\C^*$-orbit  and the  
$K$-orbit are transverse at    $(e,\oe)$.
So by dimension, the orbit $W$ of $(e,\oe)$ under $K\times\C^*$
is Zariski dense in $\YoY$.
\mypar
Now the (punctured) line $\{(te,t\oe)\,|\,t\in\C^*\}$ meets all the
$K$-orbits in $W$ and the function $t^2$ separates out the
points lying in different $K$-orbits. But the function  $T$ is 
$K$-invariant and  satisfies $T(te,t\oe)=t^2$. So $T$ separates the
$K$-orbits and the last assertion of the Lemma follows easily.\qed}
Lemma 8.2 implies that $\Pi_n/\Pi_0$ is a polynomial in $T$ and
$T\i$. But also $\Pi_n/\Pi_0$ is bihomogeneous of degree $(n,n)$
under the scaling action of $\C^*\times\C^*$ on $\YoY$. Since
$T$ is bihomogeneous of degree $(1,1)$, it follows by
bihomogeneity that 
$$\Pi_n=p_nT^n\Pi_0\eqno(8.9)$$
for some scalar $p_n\in\C^*$.
\mypar
Our problem now is to compute the scalars $p_n$. We will do this by 
writing  out the ``leading terms" of $\Pi_n$, $T^n$ and $\Pi_0$.
We formulate a notion of leading term in the following way.
Suppose  $V$ is a highest weight representation of $K$ of weight 
$\kappa$ and $S\in V\otimes\oV$ is $K$-invariant. Then we can write
$S$ as a sum of  weight vectors $S_\al$ where the weight of each 
$S_\al$ is of the form $(\al,-\al)$.  Then we call the term $S_\kappa$
 of  weight  $(\kappa,-\kappa)$ the {\it leading term}. 
\mypar
We can identify $R(\YoY)=R(Y)\otimes R(\oY)$ and then we have the 
expansion
$$T=\sum_{k=1}^rf_{u_k}\otimes\ovl{f_{u_k}}\eqno(8.10)$$
where $u_1,\dots,u_r$ is any basis of $\p$ which is orthonormal with 
respect to the Hermitian inner product on $\p$ given by 
$(u_i|u_j)=(u_i,\ou_j)_\g$. Choosing $u_1,\dots,u_r$ to be an 
orthonormal basis by weight vectors, we find (recall   $f_e=f_0$ and
$(e,\oe)_\g=1$)
$$\hbox{leading term of }T=f_0\otimes\ovl{f_0}  \eqno(8.11)$$
and also
$$\hbox{leading term of }T^n=f_0^n\otimes\ovl{f_0^n}\eqno(8.12)$$
\mypar
Next we choose an orthonormal basis on $H_{r_0+n}$ consisting
of weight vectors. This   basis then contains the highest weight
vector $f^n_0s_0/\norm f^n_0s_0\norm$ and we find
$$\hbox{leading term of }\Pi_n=
{f_0^ns_0\otimes\ovl{f_0^ns_0}\over\ns{f^n_0s_0}}\eqno(8.13)$$ 
The leading term of a product is the product of the leading terms,  and
so equating leading terms in (8.9) we get
$${f_0^ns_0\otimes\ovl{f_0^ns_0}\over\ns{f^n_0s_0}}
=p_n\left(f_0^n\otimes\ovl{f_0^n}\right)(s_0\otimes\ovl{s_0})
\eqno(8.14)$$
since $\norm s_0\norm=1$.  So using (6.30) we find
$$p_n={1\over \ns{f^n_0s_0}}={(r_0+1)_n\over n!(a)_n(b)_n}
\eqno(8.15)$$
Thus
$$\sum_{n\in\ZP}\Pi_n=\sum_{n\in\ZP}
{(r_0+1)_n\over n!(a)_n(b)_n}T^n\Pi_0=
\phantom{i}_1F_2 (r_0+1;a,b;T)\Pi_0\eqno(8.16)$$ 
This proves (8.7).  The hypergeometric series here has infinite radius   
of convergence, and so $\phantom{i}_1F_2 (r_0+1;a,b;T)$ defines a
holomorphic function on $\YoY$. Thus $\K$ is a holomorphic  section
over $\YoY$. This concludes the proof of Theorem 8.1.
\qed }

Theorem 8.1 easily gives
\cor{8.3}{$\H$ consists entirely of holomorphic sections of
$\NH$ and $\K$ is the reproducing kernel of $\H$.}

\Sec{\S9. Examples of the Quantization.}

A feature of our results is that  we can construct the representation
$\pi$ in any model of $H$ so long as we are given both the
$K$-module structure  and the $R(Y)$-module structure on $H$.
In particular the half-forms can be  completely suppressed in the
model. We  illustrate this by 2 examples.
These cases  are particularly simple ones where
the polynomial $P$  factors  in (4.13) into a product of 4 linear terms.
\vskip 1pc

{\bf Example 9.1}
Let $\gR=\so(4,4)$; this is  Case (ix) in Table 6.9 with
$p=q=4$.  Then  $\k=\sl(2,\C)^{\oplus 4}$. As $K$-modules we have 
$H\simeq R(Y)\simeq\oplus_{n\ge 0}S^n(\C^2)^{\otimes 4}$.  
A model of $H$ is given in the following
way. Let $S$ be the polynomial ring in $8$ variables
$x_{p,i}$ where $p\in\{1,\dots, 4\}$ and $i\in\{1,2\}$. Then
$H$ is the subalgebra of $S$ generated
by the $16$ products $x_{1,i}x_{2,j}x_{3,k}x_{4,l}$
where $i,j,k,l\in\{1,2\}$ so that
$$H=\oplus_{n\ge 0}
\C_n[x_{1,1},x_{1,2}]\cdot\C_n[x_{2,1},x_{2,2}]\cdot
\C_n[x_{3,1},x_{3,2}]\cdot\C_n[x_{4,1},x_{4,2}]
\subset S$$
where $\C_n[u,v]$ is the space of degree $n$
polynomials in $u$ and $v$. Notice then that
$H$ is the space of invariants in $S$ under
a scaling action of $\C^*\times\C^*\times\C^*$.
Let $\beta$ be the differential operator
on $S$ given by
$$\beta=x_{1,1}\displaystyle{\pd\phantom{xX}\over\pd x_{1,1}}+
x_{1,2}{\pd\phantom{xX}\over\pd x_{1,2}}+1$$

Then  the following 28  pseudo-differential operators on $S$
preserve $H$ and satisfy the bracket relations of $\so(8,\C)$.
I.e., these 28 operators form a basis of a complex
Lie algebra $\g$ isomorphic to $\so(8,\C)$. 
\vskip 2pc
$x_{p,1}\displaystyle{\pd\phantom{xX}\over\pd x_{p,2}}$,  
\qquad
$x_{p,2}\displaystyle{\pd\phantom{xX}\over\pd x_{p,1}}$,  
\qquad 
$x_{p,1}\displaystyle{\pd\phantom{xX}\over\pd x_{p,1}}-
x_{p,2}\displaystyle{\pd\phantom{xX}\over\pd x_{p,2}}$ 
\hskip 10pt 
\hbox{ where }$p\in\{1,2,3,4\}$ 
\vskip 1pc
$x_{1,i}x_{2,j}x_{3,k}x_{4,l}-
\displaystyle{(-1)^{i+j+k+l}\over\beta(\beta+1)}
\displaystyle{\pd^4\over
\pd x_{1,i'}\pd x_{2,j'}\pd x_{3,k'}\pd x_{4,l'}}$
\vskip 1pc
\hskip 6pc
\hbox{ where }
$\{i,i'\}=\{j,j'\}=\{k,k'\}=\{l,l'\}=\{1,2\}$
\vskip 2pc
Then $r_0=a=b=1$ in Table 6.9 and so the $\gR$-invariant inner
product on  $H$ satisfies $$\biggns{{x_{p,i}^n\over n!}}
={(1)_n(1)_n\over n!(2)_n}={1\over (n+1)}$$
where $p\in\{1,\dots,4\}$ and $i\in\{1,2\}$.
This agrees with the result in  [K]. \vskip 1pc
{\bf Example  9.2.}
Let $\gR$ be of type $G_2$; this is Case (viii) in Table 6.9.
Let $S$ be the polynomial ring in $4$ variables
$u_1,u_2,x_1,x_2$ and let  $S'\simeq R(Y)$ be the subalgebra
generated by the $8$ products
$u_i^3x_j$ and $u_i^2u_{i'}x_j$, where 
$\{i,i'\}=\{1,2\}$ and $j\in\{1,2\}$.
A model of $H$  is the $S'$-submodule
$$H=\oplus_{n\ge 0}\C_{3n+2}[u_1,u_2]\cdot\C_n[x_1,x_2]
\subset S$$
Let $\beta$ be the differential operator on $S$ given by
$$\beta=x_1\displaystyle{\pd\phantom{X}\over\pd x_{2}}+
x_{2}\displaystyle{\pd\phantom{X}\over\pd x_{2}}+1$$

The following $14$ pseudo-differential operators on $S$ preserve
$H$ and satisfy the bracket relations of  $G_2$ so that
they form a basis of a complex simple Lie algebra 
of type $G_2$. 
\vskip 2pc \hskip 2pc
$u_1\displaystyle{\pd\phantom{X}\over\pd u_2}$,  \hskip 4pt
$u_2\displaystyle{\pd\phantom{X}\over\pd u_1}$,  \hskip 4pt
$u_1\displaystyle{\pd\phantom{X}\over\pd u_1}-
 u_2\displaystyle{\pd\phantom{X}\over\pd u_2}$ \vskip 1pc
\hskip 2pc
$x_1\displaystyle{\pd\phantom{X}\over\pd x_2}$,  \hskip 4pt
$x_2\displaystyle{\pd\phantom{X}\over\pd x_1}$,  \hskip 4pt
$x_1\displaystyle{\pd\phantom{X}\over\pd x_1}-
 x_2\displaystyle{\pd\phantom{X}\over\pd x_2}$\vskip 1pc
\hskip 2pc  $u_i^3x_j-
\displaystyle{(-1)^{i+j}\over 27\beta(\beta+1)}
\displaystyle{\pd^4\over\pd u_{i'}^3\pd x_{j'}}$
\hskip 1pc\hbox{ where } $\{i,i'\}=\{j,j'\}=\{1,2\}$
\vskip 1pc  \hskip 2pc  $u_i^2u_{i'}x_j-
\displaystyle{(-1)^{i+j}\over 27\beta(\beta+1)}
\displaystyle{\pd^4\over\pd u_{i'}^2\pd u_i\pd x_{j'}}$
\hskip 1pc\hbox{ where } $\{i,i'\}=\{j,j'\}=\{1,2\}$
\vskip 1pc
Then $r_0=1$,   $a=4/3$ and $b=5/3$ in Table 6.9 so that
the $\gR$-invariant inner product on $H$ satisfies
$$\biggns{u_i^{3n+2}x_{j}^n\over n!}
={(4/3)_n(5/3)_n\over n!(2)_n}={(3n+3)!\over 3^{3n}3!n!(n+1)!(n+1)!}$$
where $i,j\in\{1,2\}$. \vskip 2pc

\def\reff{\vskip 1.5pc{\centerline {\bf References}}
\vskip 4pt \noindent}
\reff
\medskip
\def\myitem{\noindent\hangindent=1.5em\hangafter=1} 

\myitem{[A-B1]}
 A. Astashkevich and  R. Brylinski 
 Cotangent bundle models of complexified small nilpotent  orbits,
preprint

\myitem{[A-B2]}
A. Astashkevich and  R. Brylinski, 
 Lie algebras of  exotic pseudo-differential symbols,
in preparation

\myitem{[A-M]}
R. Abraham and J. E. Marsden, Foundations of Mechanics,
Addison-Wesley, 2nd edition, updated printing (1985)

\myitem{[Bi1]}
O. Biquard,
Sur les \'equations de Nahm et la structure de Poisson des alg\`ebres de
Lie semi-simples complexes.
Math. Ann. 304, 253-276 (1996).

\myitem{[Bi2]}
O. Biquard,
Twisteurs des orbites coadjointes et m\'etriques
hyper-pseudok\"ahl\'eriennes.
Preprint, Ecole Polytechnique, 1997,
Bull. Soc. Math. France, to appear.

\myitem{[B1]} R. Brylinski, Instantons and Kaehler geometry of
nilpotent orbits, in  "Representation Theories and Algebraic Geometry", 
A. Broer  ed.,   Kluwer,  85-125, in press.

\myitem{[B2]}
    R. Brylinski, 
    Quantization of the 4-Dimensional Nilpotent Orbit of  $SL(3,\R)$,
    Can. J. Math. {\bf 49} (5), 1997, 916-943

\myitem{[B3]}
    R. Brylinski,  in  preparation

\myitem{[B-K1]} 
   R. Brylinski and B. Kostant,  Nilpotent orbits,   normality and
   Hamiltonian group actions, Jour.  Amer. Math. Soc.  (1994), 269-298.

\myitem{[B-K2]}
  R. Brylinski and B. Kostant,
  Minimal  representations 
  of $E_6,E_7$ and $E_8$ and the generalized  Capelli
  Identity,  {\it  Proc. Natl. Acad. Sci. USA}  {\bf 91} (1994),
  2469-2472, and {\it Minimal representations,  geometric 
  quantization and unitarity}, Proc. Natl. Acad. Sci. USA {\bf 91} (1994),  
  6026-6029.

\myitem{[B-K3]}
  R. Brylinski and B. Kostant,
   Differential operators on conical Lagrangian  manifolds,  
  in {\it Lie Theory and Geometry: in Honor of B. Kostant},  
  Progress in Mathematics, vol. 123,  Birkhauser, Boston,  (1994), 
  65-96.

\myitem{[B-K4]}
  R. Brylinski and B. Kostant,
  Lagrangian models of minimal representations of 
  $E_6$, $E_7$ and $E_8$,  in
  {\it Functional Analysis on the Eve of the 21th
  Century: Festschrift  in Honor of the Eightieth Birthday of  I.M.
  Gelfand}, 
  Progress in Math, vol.  131,  Birkhauser, Boston (1995),  13-63.     

\myitem{[B-K5]}
  R. Brylinski and B. Kostant, Geometry of the Spherical Moment Map 
  for 
  Complex Minimal Nilpotent Orbits (in preparation).

\myitem{[F-K]}
      J. Faraut and A. Koranyi, {\it Analysis on Symmetric Cones},
      Oxford University Press, Oxford, 1994

\myitem{[G]}  D. Garfinkle, {\it  A new 
  construction of the Joseph ideal},  MIT Doctoral Thesis, 1982.
  
\myitem{[H]}
  S. Helgason, {\it Differential Geometry, Lie Groups 
and Symmetric Spaces,}
  Academic Press, New York, 1978.

\myitem{[J1]}
   A. Joseph,  Minimal realizations and spectrum generating algebras, 
   {\it  Commun. Math. Phys.} {\bf 36} (1974), 325-338.

\myitem{[J2]}
   A. Joseph,  The minimal orbit in a simple Lie algebra and its associated
   maximal ideal, {\it Ann. Scient. Ec. Norm. Sup.} {\bf 9} (1976),
   1-30.

\myitem{[Ki]} A. Kirillov,  Geometric quantization, in {\it Dynamical Systems IV},
V.I. Arnold, S.P. Novikov, eds, Encyclopaedia of Mathematical Sciences,
Springer-Verlag, 1990, 137-172.

\myitem{[K]} 
   B. Kostant, The vanishing of scalar 
  curvature  and the minimal representation of SO(4,4),
  in {\it Operator Algebras, Unitary Representations, Enveloping 
   Algebras, and Invariant Theory} (A. Connes et al, ed),
   Birkh\"auser, 1990,  85-124.

\myitem{[K-S]} 
    B. Kostant and S. Sahi, The Capelli identity, tube
   domains and the generalized Laplace transform, Adv.
   Math. {\bf 87} (1991), 71-92.  

\myitem{[Kr]}
   P. B. Kronheimer, Instantons and the geometry of the 
   nilpotent
   variety, {\it Jour. Diff. Geom.} {\bf 32} (1990), 
   473-490.

\myitem{[Sek]}
   J. Sekiguchi,  Remarks on real nilpotent orbits of a symmetric
   pair, {\it J. Math. Soc. Japan} {\bf 39} (1987), 127-138.

\myitem{[T]}
     P. Torasso,  Methode des orbites de Kirillov-Duflo et
    representations minimales des groupes simples sur un corps
     local de caracteristique nulle,
     preprint March 1996

\myitem{[Ve]}
    M. Vergne,  Instantons et correspondance de Kostant-Sekiguchi,
     C.R. Acad. Sci. Paris {\bf 320} (1995), Serie 1, 901-906.

\myitem{[Vo1]}   
  D. A. Vogan,  Singular unitary
  representations, in { Non-Commutative Harmonic Analysis and
  Lie Groups}, J. Carmona and M. Vergne, eds,
  Springer Lecture Notes 880, Springer, Berlin, 1981, 506-535

\myitem{[Vo2]}
  D. A. Vogan, 
  Associated varieties and unipotent representations, W. Barker     
and P.   Sally eds,  in {\it Harmonic Analysis on Reductive Groups},
   Birkh\"auser, 1991,  315-388

\myitem{[W]}
  N. Wallach,
  {\it Real  Reductive Groups I},
   Academic Press, 1994

\vskip 1pc 
\noindent  Department of Mathematics,
Pennsylvania State  University, University Park,  PA 
16802 

\noindent  e-mail:  rkb@math.psu.edu;\quad
WorldWideWeb URL:    http://www.math.psu.edu/rkb/
\end